\documentclass[preprint,3p,sort&compress,final,times]{elsarticle}
\usepackage[overload]{empheq}
\usepackage{amssymb}
\usepackage{mathtools}
\usepackage{bm}
\usepackage[usenames,dvipsnames]{color}

\usepackage{tikz}
\usetikzlibrary{decorations.pathreplacing,decorations.markings}

\usepackage{lineno}

\newtheorem{remark}{Remark}

\newcommand{\innerprod}[2]{\left\langle #2\,,\, #1\right\rangle}

\newcommand{\figref}[1]{Figure~\ref{#1}} 
\newcommand{\secref}[1]{Section~\ref{#1}} 
\newcommand{\appendixref}[1]{\ref{#1}} 

\newcommand{\blue}[1]{{#1}}

\newcommand{\ex}{\boldsymbol{e}_{\lambda}}
\newcommand{\ey}{\boldsymbol{e}_{\varphi}}
\newcommand{\ez}{\boldsymbol{e}_{z}}

\newcommand{\ux}{u}
\newcommand{\uy}{v}
\newcommand{\uz}{w}
\newcommand{\upar}{\boldsymbol{u}_{\parallel}}
\newcommand{\uperp}{\boldsymbol{u}_{\perp}}

\newcommand{\degreepoly}{p}

\newcommand{\omegapar}{\boldsymbol{\omega}_{\parallel}}
\newcommand{\omegaperp}{\boldsymbol{\omega}_{\perp}}
\newcommand{\omegaparpar}{\boldsymbol{\omega}_{\parallel,\parallel}}
\newcommand{\omegaparperp}{\boldsymbol{\omega}_{\parallel,\perp}}

\newcommand{\gradpar}{\nabla_{\parallel}}
\newcommand{\gradperp}{\nabla_{\perp}}
\newcommand{\curlpar}{\nabla_{\parallel}\times}
\newcommand{\curlperp}{\nabla_{\perp}\times}

\newcommand{\incidencegrad}{\boldsymbol{\mathsf{E}}^{1,0}}
\newcommand{\incidencegradpar}{\boldsymbol{\mathsf{E}}^{1,0}_{\parallel}}
\newcommand{\incidencegradperp}{\boldsymbol{\mathsf{E}}^{1,0}_{\perp}}

\newcommand{\incidencecurl}{\boldsymbol{\mathsf{E}}^{2,1}}
\newcommand{\incidencecurlparpar}{\boldsymbol{\mathsf{E}}^{2,1}_{\parallel,\parallel}}
\newcommand{\incidencecurlparperp}{\boldsymbol{\mathsf{E}}^{2,1}_{\parallel,\perp}}
\newcommand{\incidencecurlperp}{\boldsymbol{\mathsf{E}}^{2,1}_{\perp}}

\newcommand{\incidencediv}{\boldsymbol{\mathsf{E}}^{3,2}}
\newcommand{\incidencedivpar}{\boldsymbol{\mathsf{E}}^{3,2}_{\parallel}}
\newcommand{\incidencedivperp}{\boldsymbol{\mathsf{E}}^{3,2}_{\perp}}

\newcommand{\pspace}{\mathcal{P}_{h}}
\newcommand{\wspace}{\mathcal{W}_{h}}
\newcommand{\wspacepar}{\mathcal{W}_{\parallel,h}}
\newcommand{\wspaceperp}{\mathcal{W}_{\perp,h}}
\newcommand{\uspace}{\mathcal{U}_{h}}
\newcommand{\uspacepar}{\mathcal{U}_{\parallel,h}}
\newcommand{\uspaceperp}{\mathcal{U}_{\perp,h}}
\newcommand{\qspace}{\mathcal{Q}_{h}}

\newcommand{\pcardinal}{d_{\mathcal{P}}}
\newcommand{\wcardinal}{d_{\mathcal{W}}}
\newcommand{\wcardinalpar}{d_{\mathcal{W}_{\parallel}}}
\newcommand{\wcardinalperp}{d_{\mathcal{W}_{\perp}}}
\newcommand{\ucardinal}{d_{\mathcal{U}}}
\newcommand{\ucardinalpar}{d_{\mathcal{U}_{\parallel}}}
\newcommand{\ucardinalperp}{d_{\mathcal{U}_{\perp}}}
\newcommand{\qcardinal}{d_{\mathcal{Q}}}

\newcommand{\pbasis}{\epsilon^{\mathcal{P}}}
\newcommand{\wbasis}{\boldsymbol{\epsilon}^{\mathcal{W}}}
\newcommand{\wbasispar}{\boldsymbol{\epsilon}^{\mathcal{W}_{\parallel}}}
\newcommand{\wbasisperp}{\boldsymbol{\epsilon}^{\mathcal{W}_{\perp}}}
\newcommand{\ubasis}{\boldsymbol{\epsilon}^{\mathcal{U}}}
\newcommand{\ubasispar}{\boldsymbol{\epsilon}^{\mathcal{U}_{\parallel}}}
\newcommand{\ubasisperp}{\boldsymbol{\epsilon}^{\mathcal{U}_{\perp}}}
\newcommand{\qbasis}{\epsilon^{\mathcal{Q}}}

\newcommand{\ualgebraicpar}{\boldsymbol{\mathsf{u}}^{\parallel}}
\newcommand{\ualgebraicperp}{\boldsymbol{\mathsf{u}}^{\perp}}

\newcommand{\massu}{\boldsymbol{\mathsf{M}}^{\mathcal{U}}}
\newcommand{\massupar}{\boldsymbol{\mathsf{M}}^{\mathcal{U}_{\parallel}}}
\newcommand{\nupar}{\boldsymbol{\mathsf{N}}^{\mathcal{U}_{\parallel}}}
\newcommand{\massuperp}{\boldsymbol{\mathsf{M}}^{\mathcal{U}_{\perp}}}
\newcommand{\nuperp}{\boldsymbol{\mathsf{N}}^{\mathcal{U}_{\perp}}}

\newcommand{\masswpar}{\boldsymbol{\mathsf{M}}^{\mathcal{W}_{\parallel}}}
\newcommand{\masswperp}{\boldsymbol{\mathsf{M}}^{\mathcal{W}_{\perp}}}
\newcommand{\massq}{\boldsymbol{\mathsf{M}}^{\mathcal{Q}}}
\newcommand{\luq}{\boldsymbol{\mathsf{L}}^{\mathcal{U}_{\perp},\mathcal{Q}}}

\newcommand{\hdiv}[1]{H(\mathrm{div}, #1)}
\newcommand{\hcurl}[1]{H(\mathrm{curl}, #1)}
\newcommand{\ltwo}[1]{L^{2}(#1)}

\journal{Journal of Computational Physics}
\begin{document}
\begin{frontmatter}

\title{A Mixed Mimetic Spectral Element Model of the 3D Compressible Euler Equations on the Cubed Sphere}
\author[MON]{D.~Lee\corref{cor}}
\ead{davelee2804@gmail.com}
\author[TUD]{A.~Palha}

\address[MON]{Department of Mechanical and Aerospace Engineering, Monash University, Melbourne 3800, Australia}
\address[TUD]{Faculty of Aerospace Engineering, Delft University of Technology, Kluyverweg 2, 2629 HS Delft, The Netherlands}
\cortext[cor]{Corresponding author. Tel. +61 452 262 804.}

\begin{abstract}
A model of the three-dimensional rotating compressible Euler equations on the cubed sphere is presented.
The model uses a mixed mimetic spectral element discretization which allows for the exact exchanges
of kinetic, internal and potential energy via the compatibility properties of the chosen function
spaces. A Strang carryover dimensional splitting procedure is used, with
the horizontal dynamics solved explicitly and the vertical dynamics solved implicitly so as to avoid
the CFL restriction of the vertical sound waves. The function spaces used to represent the horizontal
dynamics are discontinuous across vertical element boundaries, such that each horizontal layer is solved
independently so as to avoid the need to invert a global 3D mass matrix, while the function spaces used
to represent the vertical dynamics are similarly discontinuous across horizontal element boundaries,
allowing for the serial solution of the vertical dynamics independently for each horizontal element. The
model is validated against standard test cases for baroclinic instability within an otherwise
hydrostatically and geostrophically balanced atmosphere, and a non-hydrostatic gravity wave as driven by
a temperature perturbation.
\end{abstract}


\end{frontmatter}

\section{Introduction}

Mimetic finite element families are an appealing choice for the discretization of
geophysical flow problems. This is on account of their capacity to preserve both conservation laws
and leading order balance relations in the discrete form \cite{Cotter12,McRae14,LPG18,Eldred18},
due to the compatibility properties of the chosen function spaces, as well as 
their ability to represent complex geometries such as the surface of the sphere \cite{LP18,Bauer18,Shipton18}.

The use of mimetic discretizations to represent the solution variables, and the adjoint properties
of the differential operators implied by those spaces, allows for the conservation of energy via the
exact balance of energetic exchanges, as well as the orthogonality of vorticity
evolution to those exchanges \cite{McRae14,LPG18,Eldred18}. By satisfying exactly
the balance between energetic exchanges, it is hoped that these methods may improve the
statistical behaviour of climate simulations over long time scales by mitigating against
internal biases in the representation of dynamical processes.

Several mimetic finite volume models of the compressible Euler equations on the sphere have previously
been presented \cite{Skamarock12,Gassmann13,Dubos15}, and the Raviart-Thomas family of compatible finite
elements has been chosen to form the basis of the LFric atmospheric model \cite{Natale16,Melvin19}.
Collocated spectral element methods, which in contrast to the above methods represent all solution
variables on the same function space, are a popular choice for the simulation of non-hydrostatic
atmospheric flows using both continuous \cite{Giraldo13} and \blue{discontinuous Galerkin } \cite{KG12,Bao15}
formulations. These models typically satisfy some mimetic properties, such as the divergence
theorem and the adjoint property between the gradient and divergence operators \cite{Gassner13}, however
the HOMME spectral element model \cite{TF10,Taylor11} is capable of preserving a more complete set of
mimetic properties so as to satisfy both the divergence and circulation theorems through the careful
consideration and construction of the metric terms for the various differential operators.

In the present formulation we make use of the mixed mimetic spectral element method \cite{Gerritsma11,Kreeft13},
a compatible family of function spaces with spectral error convergence. We use this method to
build on previous work on the rotating shallow water equations \cite{LPG18,LP18} in order to develop
a solver for the three-dimensional rotating compressible Euler equations on the cubed sphere.
In contrast to
collocated methods, weak formulations using mixed function spaces provide a clearer mathematical
formalism for the analysis and construction of schemes which preserve internal properties of the
governing equations. Like other weak form mimetic discretizations such as mixed finite elements
\cite{Cotter12,McRae14,Shipton18} and mimetic Galerkin differences \cite{Eldred18}, the mixed
mimetic spectral element method satisfies by construction an exact mapping between function
spaces via the various differential operators, as well as adjoint relationships between strong
form mappings to higher function spaces, and weak form mappings to lower spaces, as defined
within an appropriate Hilbert space.

The remainder of this article proceeds as follows: In Section 2 the rotating compressible Euler
equations, and their energetic properties will be introduced in the continuous form.
Section 3 will provide a brief introduction to the mixed mimetic spectral element method. Readers
are referred to references therein for more detailed discussions. Section 4 will discuss the
construction of discrete function spaces built off the mixed mimetic spectral element method
required to solve the compressible Euler equations, as well as the use of those spaces to ensure
consistent energetic exchanges in the discrete form and the associated metric transformations for
these spaces. The details of the time stepping
scheme, including the implicit vertical solver will be discussed in Section 5, and the results
for a standard baroclinic instability test case and a high resolution gravity wave
will be presented in Section 6. Section 7 will
discuss the conclusions of this work and the future directions we intend to pursue.

\section{The rotating compressible Euler equations}

The compressible Euler equations for a shallow atmosphere may be expressed as \cite{Nevir09,Gassmann13}
\begin{subequations}
	\label{eq:compressible_euler}
	\begin{align}
		\frac{\partial\bm u}{\partial t} + (\bm\omega + \bm f)\times\bm u + \nabla\Bigg(\frac{1}{2}\|u\|^2 + gz\Bigg) + 	\theta\nabla\Pi &= 0,\label{eq::mom}\\
		\frac{\partial\rho}{\partial t} + \nabla\cdot(\rho\bm u) &= 0,\label{eq::dens}\\
		\frac{\partial\rho\theta}{\partial t} + \nabla\cdot(\rho\theta \bm u) &= 0,\label{eq::temp}
	\end{align}
\end{subequations}
where $\bm u = u\ex + v\ey + w\ez$ are the zonal, meridional and vertical velocity components respectively,
$\rho$ is the density, $\bm f=f\ez$ is the Coriolis term, $g$ is the acceleration due to gravity, $\theta$
is the potential temperature, and $\Pi$ is the Exner pressure (including the specific heat at constant pressure).
The last two are defined with respect to the standard thermodynamic variables of temperature, $T$, and pressure, $p$, as
\begin{subequations}
	\label{eq:definition_exner_potential_temperature}
	\begin{align}
		\Pi &:= c_p\Bigg(\frac{p}{p_0}\Bigg)^{\frac{R}{c_p}}, \label{eq:exener_definition} \\
	 	\theta &:= \frac{c_{p}T}{\Pi}, \label{eq:potential_temperature_definition} \\
	 	p &= \rho RT. \label{eq:pressure_definition}
	\end{align}
\end{subequations}
For these identities we used $c_p$ for the specific heat at constant pressure, $p_0$ for the reference pressure, and $R$ for the ideal gas constant.
We may remove the direct dependence on pressure from the system by simply substituting expression \eqref{eq:pressure_definition} for pressure into \eqref{eq:exener_definition} and \eqref{eq:potential_temperature_definition}, obtaining
\begin{subequations}
	\label{eq::eos}
	\begin{align}
		\Pi &= c_p\Bigg(\frac{\rho R\theta\Pi}{p_0c_p}\Bigg)^{\frac{R}{c_p}} =
c_p\Bigg(\frac{\rho R\theta}{p_0}\Bigg)^\frac{R}{c_v}, \label{eq:pi_from_theta_rho} \\
		\theta &= \frac{\rho R T^{-\frac{c_{v}}{R}}}{p_{0}},
	\end{align}
\end{subequations}
where $c_v = c_p-R$ is the specific heat at constant volume.

The potential temperature/Exner pressure form of the pressure gradient term, $\theta\nabla\Pi$, in \eqref{eq::mom} is equivalent to the
standard density/pressure form since
\begin{equation}
\frac{1}{\rho}\nabla p \stackrel{\eqref{eq:pressure_definition}}{=} \frac{RT}{p}\nabla p \stackrel{\eqref{eq:exener_definition} + \eqref{eq:potential_temperature_definition}}{=} \frac{R\theta}{p}\Bigg(\frac{p}{p_0}\Bigg)^{\frac{R}{c_p}}\nabla p =
\frac{R\theta}{p_0^{R/c_p}}p^{R/c_p - 1}\nabla p = \frac{R\theta}{p_0^{R/c_p}}\frac{c_p}{R}\nabla p^{R/c_p} =
c_p\theta\nabla\Bigg(\frac{p}{p_0}\Bigg)^{\frac{R}{c_p}} \stackrel{\eqref{eq:exener_definition} + \eqref{eq:potential_temperature_definition}}{=} \theta\nabla\Pi.
\end{equation}
One advantage of the Exner pressure/potential temperature representation of the thermodynamics is the formulation
of the temperature equation in flux form \eqref{eq::temp}, which allows us to exploit the adjoint relationship
between gradient and divergence in the mimetic framework in order to preserve energetic exchanges.

To obtain a closed system for the solution of the compressible Euler equations, \eqref{eq:compressible_euler} and
\eqref{eq::eos} must be supplemented by Dirichlet and Neumann boundary conditions. The following identities are
imposed as Dirichlet boundary conditions on the $z$-component of velocity, $w$
\begin{equation}
	\left.w\right|_{z = 0} = \left.w\right|_{z = z^{\mathrm{top}}} = 0, \label{eq:dirichlet_bc_velocity}
\end{equation}
where $z^{\mathrm{top}}$ corresponds to the $z$-coordinates of the top boundary of the domain.
For Neumann boundary conditions the following identities are imposed
\begin{equation}\label{eq::bcs_neumann}
	\frac{\partial\Pi}{\partial z}\Bigg|_{z=0} = \frac{\partial\Pi}{\partial z}\Bigg|_{z=z^{top}} = 0.
\end{equation}

Note that in this formulation we have invoked the  shallow atmosphere approximation, for which gravity is constant throughout the fluid
column, the height of the fluid column is negligible with respect to the earth's radius, and the
horizontal components of the Coriolis term are omitted \cite{White05}.

\subsection{Energetics}

Before introducing the discrete form of the Euler equations, we analyse the energetics of the continuous system.
This will help to guide our choice of function spaces for the various solution variables for the discrete form.

\subsubsection{Kinetic, potential, and internal energy}
	The kinetic energy, $K$, is defined as
	\begin{equation}
		K := \frac{1}{2}\langle\boldsymbol{u}, \rho\boldsymbol{u} \rangle = \frac{1}{2}\int_{\Omega}\rho \|u\|^{2}, \label{eq:kinetic_energy_definition}
	\end{equation}
	where $\|u\| := \langle \boldsymbol{u}, \boldsymbol{u} \rangle$, and $\langle \cdot, \cdot \rangle$ is the $L^{2}$ inner product given as usual as
	\begin{equation}
		\langle f, g \rangle := \int_{\Omega} f g\,\mathrm{d}\Omega\,,
	\end{equation}
	for scalar fields, and as
	\begin{equation}
		\langle \boldsymbol{u}, \boldsymbol{v} \rangle := \int_{\Omega} \boldsymbol{u} \cdot \boldsymbol{v}\,\mathrm{d}\Omega\,,
	\end{equation}
	for vector fields.

The time variation of kinetic energy is obtained by summing the $L^{2}$ inner product, between the momentum equation, \eqref{eq::mom},
and $\rho\bm u$, and between the continuity equation, \eqref{eq::dens}, and $\frac{1}{2}\|u\|^2$
\begin{equation}
	\frac{\partial K}{\partial t}  = -\langle g ,\rho \uz\rangle - \langle\rho\bm u,\theta\nabla\Pi\rangle\,,\label{eq::ke_evol}
\end{equation}
where again $\uz$ is the $z$-component of the velocity field, $\boldsymbol{u}$.

	The potential energy, $P$, is given by
	\begin{equation}
		P := \innerprod{gz}{\rho} = \int_{\Omega}\rho gz \,\mathrm{d}\Omega, \label{eq:potential_energy_definition}
	\end{equation}
	and its time derivative follows directly
\begin{equation}
\frac{\partial P}{\partial t} = \innerprod{\frac{\partial\rho}{\partial t}}{gz} \stackrel{\eqref{eq::dens}}{=}
-\langle gz,\nabla\cdot(\rho\bm u)\rangle
= \langle g,\rho \uz\rangle\,,\label{eq:time_variation_potential_energy}
\end{equation}
where we have used integration by parts on the last identity and assumed periodic boundary conditions in the horizontal directions, together with homogeneous boundary conditions for the vertical component of the velocity field, \eqref{eq:dirichlet_bc_velocity}.

	The internal energy, $I$,  is defined as
	\begin{equation}
		I := \int_{\Omega}c_{v}\rho T \,\mathrm{d}\Omega\stackrel{\eqref{eq:potential_temperature_definition}}{=} \int_{\Omega}\frac{c_{v}}{c_{p}}\rho\theta\Pi \,\mathrm{d}\Omega\stackrel{\eqref{eq:pi_from_theta_rho}}{=} \int_{\Omega}c_{v}\rho\theta\left(\frac{R\rho\theta}{p_{0}}\right)^{\frac{R}{c_{v}}} \,\mathrm{d}\Omega= \int_{\Omega} c_{v} \left(\frac{R}{p_{0}}\right)^{\frac{R}{c_{v}}}\left(\rho\theta\right)^{\frac{c_{p}}{c_{v}}}\,\mathrm{d}\Omega\,. \label{eq:internal_energy_definition}
	\end{equation}
After some manipulation, the time variation of internal energy is given by
	\begin{equation}
		\frac{\partial I}{\partial t} = - \innerprod{\Pi}{\nabla\cdot(\rho\theta \bm u)} = \innerprod{\theta\nabla\Pi}{\rho\bm u} \,. \label{eq:time_variation_internal_energy}
	\end{equation}
	where integration by parts was used on the last identity, together with homogeneous boundary conditions for $\boldsymbol{u}$ and periodic boundary conditions on the horizontal directions.

\subsubsection{Conservation of total energy}
	 Following \cite{Gassmann13}, the total energy of the system, $\mathcal{H}$, is given as the sum of kinetic, $K$, potential, $P$, and internal, $I$, energy
	 \begin{equation}
	 	\mathcal{H} := K + P + I = \int_{\Omega}\frac{1}{2}\rho u^2\,\mathrm{d}\Omega +  \int_{\Omega}\rho gz\,\mathrm{d}\Omega +  \int_{\Omega}\frac{c_v}{c_p}\Theta\Pi\,\mathrm{d}\Omega\,. \label{eq::tot_en}
	 \end{equation}
	 where we used $\Theta := \rho\theta$.

For the proof of conservation of total energy $\mathcal{H}$, \eqref{eq::tot_en}, first consider the column vectors
\begin{equation}\label{eq:en_con_1}
	\boldsymbol{\mathsf{a}} :=
	\left[
		\begin{array}{ccc}
			\boldsymbol{u} & \rho & \Theta
		\end{array}
	\right]^{\top}\,,
\end{equation}
and
\begin{equation}\label{eq:en_con_2}
	\boldsymbol{\mathsf{h}} :=
	\left[
		\begin{array}{ccc}
			\boldsymbol{U} & \Phi & \Pi
		\end{array}
	\right]^{\top}\,,
\end{equation}
where $\boldsymbol{U} := \rho \boldsymbol{u}$, and $\Phi := \frac{1}{2}u^{2} + gz$.
Introducing the skew-symmetric operator
\begin{equation}\label{eq:en_con_3}
	\boldsymbol{\mathsf{B}} :=
	\left[
		\begin{array}{ccc}
			-\bm q\times(\cdot) & -\nabla(\cdot) & -\theta\nabla(\cdot)\\
             -\nabla\cdot(\cdot) & 0 & 0 \\
             -\nabla\cdot(\theta\cdot) & 0 & 0\\
		\end{array}
	\right]\,,
\end{equation}
where $\bm q = (\bm\omega + \bm f)/\rho$ is the potential vorticity,
the original  prognostic equations, \eqref{eq::mom}-\eqref{eq::temp}, may be rewritten as
\begin{equation}
	\frac{\partial\boldsymbol{\mathsf{a}}}{\partial t} = \boldsymbol{\mathsf{B}} \,\boldsymbol{\mathsf{h}}\,. \label{eq:da_dt_B_h}
\end{equation}
Note now that the variational derivatives of $\mathcal{H}$ with respect to the prognostic variables $\boldsymbol{u}$, $\rho$, and $\Theta$, are
\begin{equation}\label{eq::variations}
	\frac{\delta\mathcal{H}}{\delta\bm u} = \rho\bm u = \bm U,\qquad
	\frac{\delta\mathcal{H}}{\delta\rho} = \frac{1}{2}u^2 + gz = \Phi,\qquad
	\frac{\delta\mathcal{H}}{\delta\Theta} = \Pi\,,
\end{equation}
and therefore
\begin{equation}
	\boldsymbol{\mathsf{h}} = \frac{\delta\mathcal{H}}{\delta\boldsymbol{\mathsf{a}}}\,. \label{eq:h_equals_H}
\end{equation}
Substituting \eqref{eq:h_equals_H} into \eqref{eq:da_dt_B_h} yields
\begin{equation}
	\frac{\partial\boldsymbol{\mathsf{a}}}{\partial t} = \boldsymbol{\mathsf{B}} \,\frac{\delta\mathcal{H}}{\delta\boldsymbol{\mathsf{a}}}\,. \label{eq:da_dt_B_H}
\end{equation}

Conservation of total energy follows directly since \cite{Hairer06}
\begin{equation}
\frac{\partial\mathcal{H}}{\partial t} =
\frac{\delta\mathcal{H}}{\delta\boldsymbol{\mathsf{a}}}
\cdot\frac{\partial\boldsymbol{\mathsf{a}}}{\partial t} \stackrel{\eqref{eq:da_dt_B_H}}{=} \frac{\delta\mathcal{H}}{\delta\boldsymbol{\mathsf{a}}}\cdot\left(
 \boldsymbol{\mathsf{B}}\frac{\delta\mathcal{H}}{\delta\boldsymbol{\mathsf{a}}}\right) = 0,
\end{equation}
where the last identity follows from the skew-symmetry of $\boldsymbol{\mathsf{B}}$. Note that here the dot product involves not only a summation over the elements of the vectors but also an integration over $\Omega$, e.g.
\[
	\frac{\delta\mathcal{H}}{\delta\bm a} \cdot\frac{\partial\bm a}{\partial t} = \int_{\Omega} \sum_{i=1}^{3} \,\frac{\delta\mathcal{H}}{\delta\mathsf{a}_{i}}\frac{\partial\mathsf{a}_{i}}{\partial t} \mathrm{d}\Omega\,.
\]

\section{Mimetic polynomial basis functions}

	\subsection{1D mimetic polynomial function spaces} \label{eq:1d_basis_functions}
The mixed mimetic spectral element method is built off two types of one-dimensional polynomials: one
associated with nodal interpolation, and the other with integral interpolation (histopolation)
\cite{robidoux-polynomial,Gerritsma11}. Subsequently, these two types of polynomials will be combined
to generate the family of three-dimensional basis functions used to discretize the system.

		\subsubsection{Nodal polynomial basis functions}
Consider the canonical interval $I=[-1,1]\subset\mathbb{R}$ and the Legendre polynomials, $L_{\degreepoly}(\xi)$ of degree $\degreepoly$ with
$\xi\in I$. The $\degreepoly+1$ roots, $\xi_{i}$, of the polynomial $(1-\xi^{2})\frac{\mathrm{d}L_{\degreepoly}}{\mathrm{d}\xi}$ are called
Gauss-Lobatto-Legendre (GLL) nodes and satisfy $-1 = \xi_{0} < \xi_{1} < \dots < \xi_{\degreepoly-1} < \xi_{\degreepoly} = 1$. Let $l^{\degreepoly}_{i}(\xi)$
be the Lagrange polynomial of degree $\degreepoly$ through the GLL $\xi_{i}$, such that
            \begin{equation}
                l^{\degreepoly}_{i}(\xi_{j}) = \delta_{i,j},\quad i,j = 0,\dots,\degreepoly\,,\label{eq::nodal_basis_polynomials}
            \end{equation}
            where $\delta_{i,j}$ is the Kronecker delta.
            The explicit form of these Lagrange polynomials is given by
            \begin{equation}
                l^{\degreepoly}_{i}(\xi) = \prod_{\substack{k=0\\k\neq i}}^{\degreepoly}\frac{\xi-\xi_{k}}{\xi_{i}-\xi_{k}}\,. \label{eq::lagrange_interpolants}
            \end{equation}

Let $q_h(\xi)$ be a polynomial of degree $\degreepoly$ defined on $I=[-1,1]$ and $q_{i} = q_h(\xi_{i})$, then the expansion of $q_h(\xi)$ in terms of Lagrange
polynomials is given by
            \begin{equation}
                q_h(\xi) := \sum_{i=0}^{\degreepoly}q_{i}l^{\degreepoly}_{i}(\xi)\,. \label{eq::nodal_polynomial_expansion}
            \end{equation}
Because the expansion coefficients in \eqref{eq::nodal_polynomial_expansion} are given by the value of $q_h$ in the nodes $\xi_i$,
we refer to this interpolation as a {\em nodal interpolation} and we will denote the Lagrange polynomials in \eqref{eq::lagrange_interpolants}
by \emph{nodal polynomials}.

		\subsubsection{Histopolant polynomial basis functions}
\blue{In addition to the nodal basis functions defined above, a second set of basis functions is required in order to discretize 
integral based quantities across edges. Together these nodal and histopolant basis functions will then be used to construct a compatible 
family of finite element function spaces in multiple dimensions. }
Using the nodal polynomials we can define \blue{this } set of basis polynomials, $e^{\degreepoly}_{i}(\xi)$, as
\begin{equation}
e^{\degreepoly}_{i}(\xi) := - \sum_{k=0}^{i-1}\frac{\mathrm{d}l^{\degreepoly}_{k}(\xi)}{\mathrm{d}\xi}\,, \qquad i=1,\dots,\degreepoly\,. \label{eq::histopolant_polynomials_definition}
\end{equation}
These polynomials $e^{\degreepoly}_{i}(\xi)$ have polynomial degree $\degreepoly-1$ and satisfy,
\begin{equation}
\int_{\xi_{j-1}}^{\xi_{j}}e^{\degreepoly}_{i}(\xi)\,\mathrm{d}\xi = \delta_{i,j},\quad i,j = 1,\dots,p\,. \label{eq::histopolant_polynomials_properties}
\end{equation}
Using \eqref{eq::histopolant_polynomials_definition} the integral of $e^{\degreepoly}_{i}(\xi)$ becomes
\cite{robidoux-polynomial,Gerritsma11}
            \[
\int_{\xi_{j-1}}^{\xi_{j}}e^{\degreepoly}_{i}(\xi)\,\mathrm{d}\xi =
- \int_{\xi_{j-1}}^{\xi_{j}}\sum_{k=0}^{i-1}\frac{\mathrm{d}l^{\degreepoly}_{k}(\xi)}{\mathrm{d}\xi} =
- \sum_{k=0}^{i-1}\int_{\xi_{j-1}}^{\xi_{j}}\frac{\mathrm{d}l^{\degreepoly}_{k}(\xi)}{\mathrm{d}\xi} =
- \sum_{k=0}^{i-1} \left(l^{\degreepoly}_{k}(\xi_{j}) - l^{\degreepoly}_{k}(\xi_{j-1})\right) = - \sum_{k=0}^{i-1} \left(\delta_{k,j} -\delta_{k,j-1}\right) =
\delta_{i,j}\,.
            \]

\begin{figure}
    \begin{center}
         \includegraphics[width=0.4\textwidth]{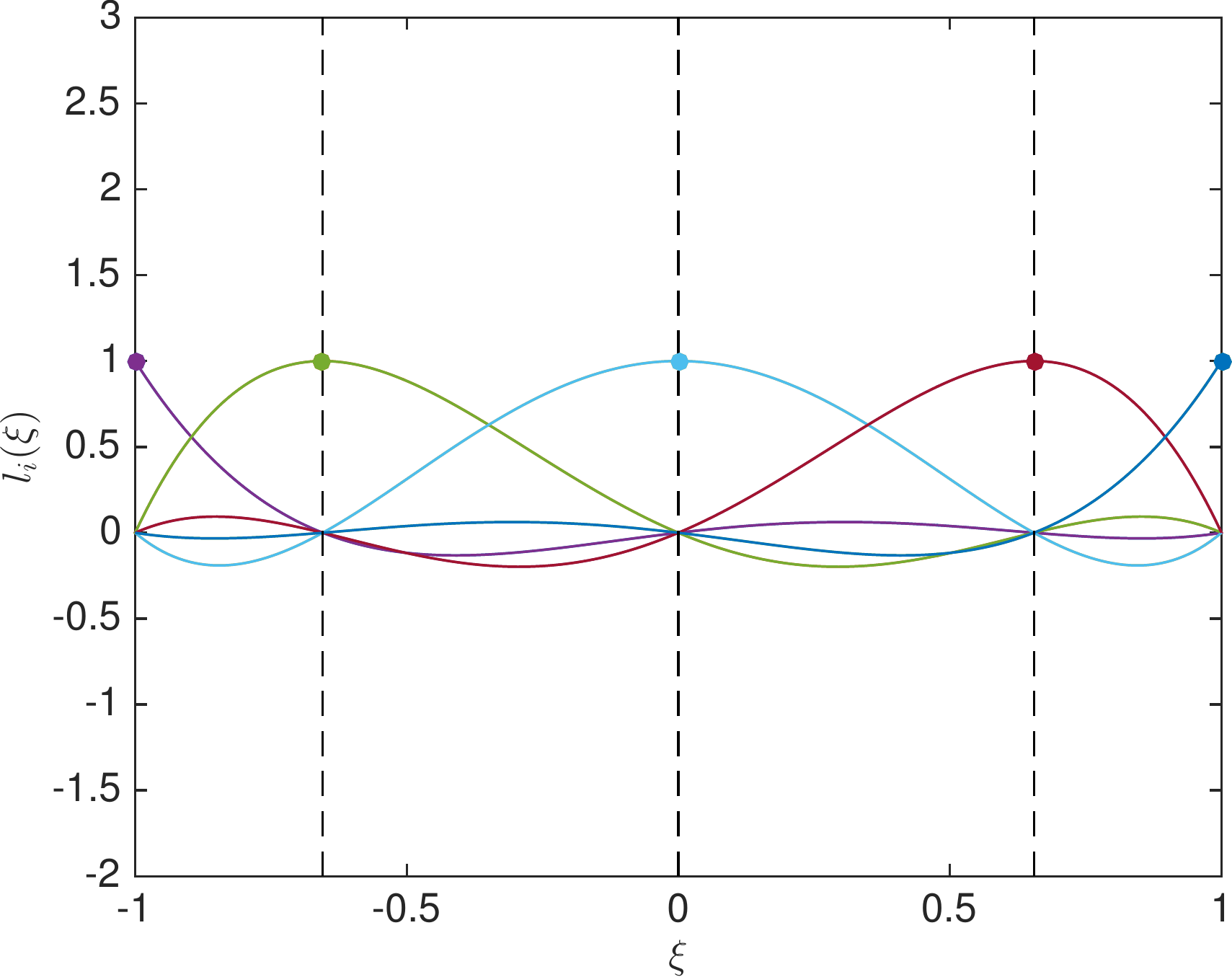} \hspace{1cm}
         \includegraphics[width=0.4\textwidth]{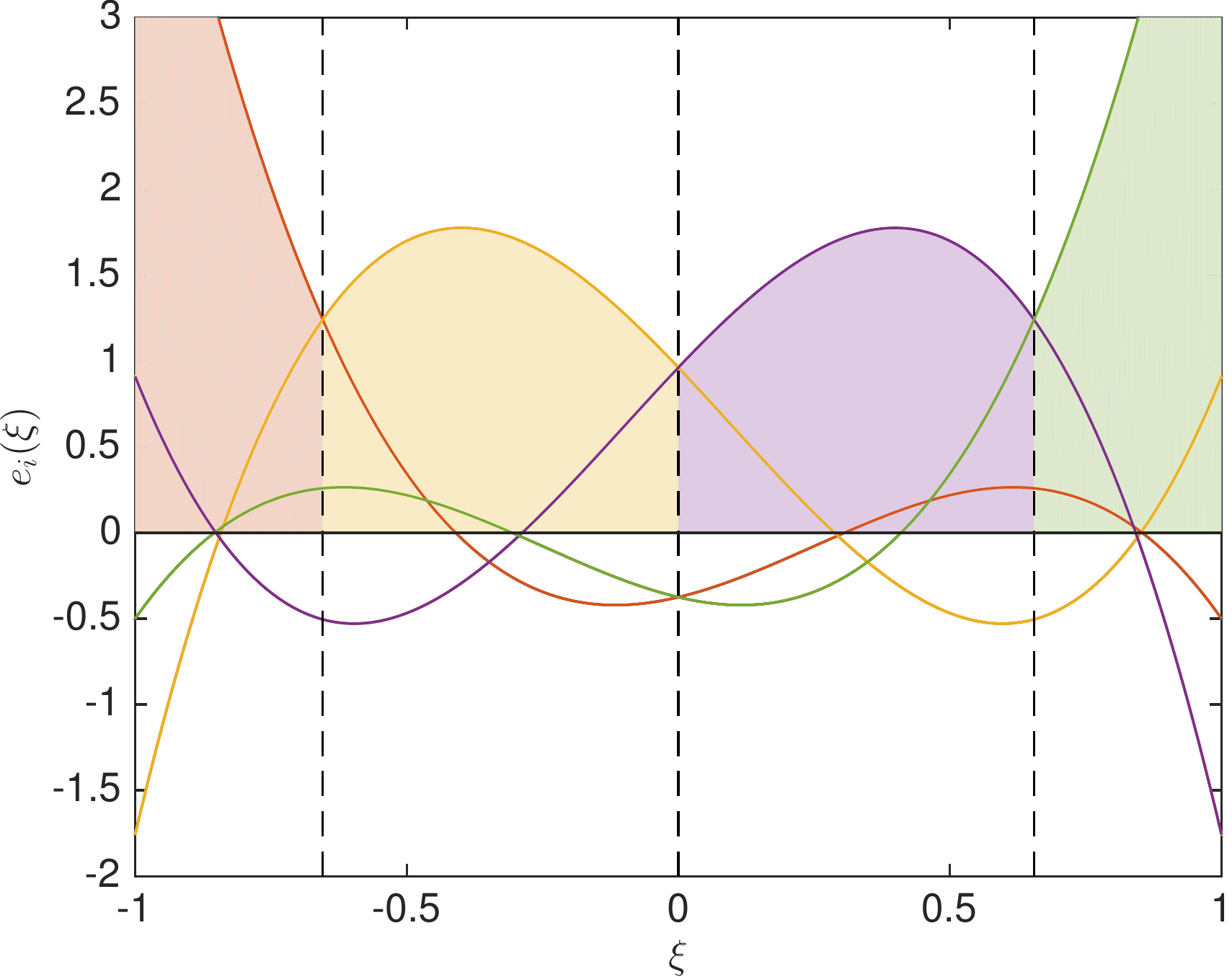}
\caption{Basis polynomials for nodal polynomials of degree $p=4$. Left: nodal polynomials, $l^p_i(\xi_j)$,
right: edge polynomials, $e^p_i(\xi_j)$.}\label{fig::basis_polynomials}
    \end{center}
\end{figure}

Let $g_h(\xi)$ be a polynomial of degree $(\degreepoly-1)$ defined on $I=[-1,1]$ and $g_{i} = \int_{\xi_{i-1}}^{\xi_{i}}g_h(\xi)\,\mathrm{d}\xi$,
then its expansion in terms of the polynomials $e_i^{\degreepoly}(\xi)$ is given by
\begin{equation}
g_{h}(\xi) := \sum_{i=1}^{\degreepoly}g_{i}e^{\degreepoly}_{i}(\xi)\;.\label{eq:1D_expansion_in_edge_polynomials}
\end{equation}
Because the expansion coefficients in \eqref{eq:1D_expansion_in_edge_polynomials} are the integral values of $g_h(\xi)$,
we denote the polynomials in \eqref{eq::histopolant_polynomials_definition} by \emph{histopolant polynomials} and refer
to \eqref{eq:1D_expansion_in_edge_polynomials} as {\em histopolation}. It can be shown, \cite{robidoux-polynomial,Gerritsma11},
that if $q_h(\xi)$ is expanded in terms of nodal polynomials, as in \eqref{eq::nodal_polynomial_expansion}, then the
expansion of its derivative $\frac{\mathrm{d}q_h(\xi)}{\mathrm{d}\xi}$ in terms of histopolant, or edge polynomials is
\begin{align}
     \left(\frac{\mathrm{d}q_h(\xi)}{\mathrm{d}\xi}\right)_{h} & =
    \sum_{i=1}^{\degreepoly} \left(\int_{\xi_{i-1}}^{\xi_{i}}\frac{\mathrm{d}q_h(\xi)}{\mathrm{d}\xi}\mathrm{d}\xi\right)e^{\degreepoly}_{i}(\xi) =
    \sum_{i=1}^{\degreepoly} \left(q_h(\xi_{i}) - q_h(\xi_{i-1})\right)e^{\degreepoly}_{i}(\xi) \nonumber \\
 &= \sum_{i=1}^{\degreepoly} \left(q_{i} - q_{i-1}\right)e^{\degreepoly}_{i}(\xi) = \sum_{i=1,j=0}^{p}\mathsf{E}^{1,0}_{i,j}q_{j}e^{\degreepoly}_{i}(\xi)\;,
\end{align}
where $\mathsf{E}^{1,0}_{i,j}$ are the coefficients of the
$\degreepoly\times(\degreepoly+1)$ matrix $\boldsymbol{\mathsf{E}}^{1,0}$ for the one dimensional case, hereafter
referred to as an \emph{incidence} matrix (see \appendixref{appendix:incidence_matrix} for details and \cite{LPG18, palha2014, gerritsmaSEMASIMAI2018} 
for an extensive discussion). The following identity holds (Commuting property)
\begin{equation}
     \left(\frac{\mathrm{d}q(\xi)}{\mathrm{d}\xi}\right)_{h} = \frac{\mathrm{d}q_{h}(\xi)}{\mathrm{d}\xi}\,.
\end{equation}
For an example of the one-dimensional basis polynomials corresponding to $\degreepoly=4$, see \figref{fig::basis_polynomials}.

	\subsection{3D mimetic polynomial function spaces} \label{sec:three_d_mimetic_polynomial_basis_functions}
		A fundamental element in the proposed discretization for the compressible Euler equations, \eqref{eq:compressible_euler}, is the de Rham sequence of function spaces in the domain $\Omega \subset \mathbb{R}^{3}$
		\begin{equation}
			\mathbb{R}\longrightarrow H^{1}(\Omega) \stackrel{\nabla}{\longrightarrow} H(\mathrm{curl}, \Omega) \stackrel{\nabla\times}{\longrightarrow} H(\mathrm{div}, \Omega) \stackrel{\nabla\cdot}{\longrightarrow} L^{2} (\Omega) \longrightarrow 0\,,
		\end{equation}
		where, as usual, the space $H^{1}(\Omega)$ represents square integrable functions over $\Omega$ whose gradient is also square integrable, the function spaces $H(\mathrm{curl}, \Omega)$ and $H(\mathrm{div}, \Omega)$ contain square integrable vector fields over $\Omega$ with square integrable curl and divergence, respectively, and the function space $L^{2}(\Omega)$ contains square integrable functions.

		More specifically, this work relies on approximating polynomial spaces $\pspace(\Omega) \subset H^{1}(\Omega)$, $\wspace(\Omega)\subset H(\mathrm{curl}, \Omega)$, $\uspace(\Omega) \subset H(\mathrm{div}, \Omega)$, and $\qspace(\Omega) \subset L^{2}(\Omega)$ such that
		\begin{equation}
			\mathbb{R}\longrightarrow\pspace(\Omega) \stackrel{\nabla}{\longrightarrow} \wspace(\Omega) \stackrel{\nabla\times}{\longrightarrow} \uspace(\Omega)  \stackrel{\nabla\cdot}{\longrightarrow} \qspace(\Omega) \longrightarrow 0\,. \label{eq:de_rham_sequence_discrete}
		\end{equation}
		These 3-dimensional polynomial function spaces may be constructed from tensor products of the 1-dimensional function spaces presented in \secref{eq:1d_basis_functions}. Moreover, each of these polynomial function spaces has an associated finite set of basis functions $\pbasis_{i}$, $\wbasis_{i}$, $\ubasis_{i}$, and $\qbasis_{i}$, such that
		\begin{equation}
			\pspace = \mathrm{span}\{\pbasis_{1}, \dots, \pbasis_{\pcardinal}\}, \quad \wspace = \mathrm{span}\{\wbasis_{1}, \dots, \wbasis_{\wcardinal}\}, \quad \uspace = \mathrm{span}\{\ubasis_{1}, \dots, \ubasis_{\ucardinal}\}, \quad\mathrm{and}\quad\qspace = \mathrm{span}\{\qbasis_{1}, \dots, \qbasis_{\qcardinal}\}. \label{eq:function_spaces_span}
		\end{equation}
		As previously discussed, e.g. \cite{palha2014}, these basis functions  are given by
		\begin{align*}
			\pbasis_{m}(\xi,\eta,\zeta; \degreepoly) &:= l_{i}^{\degreepoly}(\xi)l_{j}^{\degreepoly}(\eta)l_{k}^{\degreepoly}(\zeta), \qquad m = i + j\degreepoly_{+} + k\degreepoly_{+}^{2}, \qquad i,j,k = 0, \cdots, \degreepoly, \\ \\
			\wbasis_{m}(\xi,\eta,\zeta; \degreepoly) &:=
			\begin{dcases}
				e_{i}^{\degreepoly}(\xi)l_{j}^{\degreepoly}(\eta)l_{k}^{\degreepoly}(\zeta)\ex,  \qquad m = i + j\degreepoly + k\degreepoly\degreepoly_{+} - 1, \qquad i = 1,\dots, \degreepoly, \quad j,k = 0, \cdots, \degreepoly, \\
				l_{i}^{\degreepoly}(\xi)e_{j}^{\degreepoly}(\eta)l_{k}^{\degreepoly}(\zeta)\ey,  \qquad m = i + j_{-}\degreepoly_{+} + k\degreepoly\degreepoly_{+} + \degreepoly\degreepoly_{+}^{2}, \qquad i,k = 0, \cdots, \degreepoly, \quad j = 1, \dots, \degreepoly, \\
				l_{i}^{\degreepoly}(\xi)l_{j}^{\degreepoly}(\eta)e_{k}^{\degreepoly}(\zeta)\ez,  \qquad m = i + j\degreepoly_{+} + k_{-}\degreepoly_{+}^{2} + 2\degreepoly\degreepoly_{+}^{2} + 1, \qquad i,j = 0, \cdots, \degreepoly, \quad k = 1, \dots,\degreepoly\,
			\end{dcases}\\ \\
			\ubasis_{m}(\xi,\eta,\zeta; \degreepoly) &:=
			\begin{dcases}
				l_{i}^{\degreepoly}(\xi)e_{j}^{\degreepoly}(\eta)e_{k}^{\degreepoly}(\zeta)\ex,  \qquad m = i + j_{-}\degreepoly_{+} + k_{-}\degreepoly\degreepoly_{+}, \qquad i = 0,\dots, \degreepoly, \quad j,k = 1, \cdots, \degreepoly, \\
				e_{i}^{\degreepoly}(\xi)l_{j}^{\degreepoly}(\eta)e_{k}^{\degreepoly}(\zeta)\ey,  \qquad m = i_{-} + j\degreepoly + k_{-}\degreepoly\degreepoly_{+} + \degreepoly_{+}\degreepoly^{2}, \qquad i,k = 1, \cdots, \degreepoly, \quad j = 0, \dots, \degreepoly, \\
				e_{i}^{\degreepoly}(\xi)e_{j}^{\degreepoly}(\eta)l_{k}^{\degreepoly}(\zeta)\ez,  \qquad m = i_{-} + j_{-}\degreepoly + k\degreepoly^{2} + 2\degreepoly_{+}\degreepoly^{2} , \qquad i,j = 1, \cdots, \degreepoly, \quad k = 0, \dots,\degreepoly\,
			\end{dcases} \\ \\
			\qbasis_{m}(\xi,\eta,\zeta; \degreepoly) &:= e_{i}^{\degreepoly}(\xi)e_{j}^{\degreepoly}(\eta)e_{k}^{\degreepoly}(\zeta), \qquad m = i + j_{-}\degreepoly + k_{-}\degreepoly^{2} - 1, \qquad i,j,k = 1, \cdots, \degreepoly,
		\end{align*}
		where for compactness, the subscripts $+$ and $-$ mean addition and subtraction of 1, e.g. $\degreepoly_{+} := \degreepoly + 1$ and $i_{-} :=  i - 1$.

		Moreover, the basis functions satisfy the following identities, see for example \cite{palha2014}:
		\begin{equation}
			\nabla\pbasis_{j} = \sum_{k=0}^{\qcardinal}\mathsf{E}_{k,j}^{1,0}\wbasis_{k},\qquad\nabla\times\wbasis_{j} = \sum_{k=0}^{\ucardinal}\mathsf{E}_{k,j}^{2,1}\ubasis_{k},\qquad \mathrm{and} \qquad \nabla\cdot\ubasis_{j} = \sum_{k=0}^{\qcardinal}\mathsf{E}_{k,j}^{3,2}\qbasis_{k}, \label{eq_incidence_matrices}
		\end{equation}
		where $\boldsymbol{\mathsf{E}}^{1,0}$, $\boldsymbol{\mathsf{E}}^{2,1}$, and $\boldsymbol{\mathsf{E}}^{3,2}$ are the incidence matrices corresponding to the discrete versions of the differential operators grad, curl, and div, respectively (see Appendix for details).

\section{Numerical discretization}

	\subsection{Splitting into horizontal and vertical contributions}
		Consider the following splitting into the horizontal, $\upar$, and vertical, $\uperp$, components of the velocity field $\boldsymbol{u} = u\ex + v\ey + w\ez$
		\begin{equation}
			\upar := \ux\ex + \uy\ey\,, \qquad \uperp := \uz\ez\,. \label{eq:vector_splitting}
		\end{equation}
		Moreover, let $\gradpar$ and $\gradperp$ represent the horizontal and vertical components of the gradient operator of a scalar field $\rho$
		\begin{equation}
			\gradpar \rho := \frac{\partial\rho}{\partial\lambda} \ex + \frac{\partial\rho}{\partial\phi} \ey\,, \qquad \gradperp\rho := \frac{\partial\rho}{\partial z} \ez\,. \label{eq:grad_splitting}
		\end{equation}
		In a similar way, $\curlpar$ and $\curlperp$ represent, respectively, the horizontal and vertical components of the curl of a vector field $\boldsymbol{u} = u\ex + v\ey + w\ez$
		\begin{equation}
			\curlpar \boldsymbol{u} := \left(\frac{\partial w}{\partial\phi} - \frac{\partial v}{\partial z}\right)\ex + \left(\frac{\partial u}{\partial z} - \frac{1}{r\cos(\phi)}\frac{\partial w}{\partial\lambda}\right)\ey\,, \qquad  \curlperp \boldsymbol{u} := \frac{1}{r\cos(\phi)}\left(\frac{\partial v}{\partial\lambda} - \frac{\partial(\cos(\phi) u)}{\partial\phi}\right)\ez\,. \label{eq:curl_splitting}
		\end{equation}
Note that we have assumed the shallow atmosphere approximation of constant radius, $r$, in the above expressions.
		From \eqref{eq:grad_splitting} and \eqref{eq:curl_splitting} follows directly that
		\[
			\nabla\rho = \gradpar\rho + \gradperp\rho\,, \quad\mathrm{and}\quad \nabla\times\boldsymbol{u} = \curlpar\boldsymbol{u} +\curlperp\boldsymbol{u}\,.
		\]
		With \eqref{eq:vector_splitting} and \eqref{eq:curl_splitting} it is possible to rewrite the definition of vorticity, $\boldsymbol{\omega} := \nabla\times\boldsymbol{u}$, as
		\begin{equation}
			\boldsymbol{\omega} = \underbrace{\overbrace{\curlpar\upar}^{\omegaparpar}+ \overbrace{\curlpar\uperp}^{\omegaparperp}}_{\omegapar} + \underbrace{\curlperp\upar}_{\omegaperp}\,. \label{eq:vorticity_splitting}
		\end{equation}

		Using \eqref{eq:vector_splitting}, \eqref{eq:grad_splitting}, \eqref{eq:curl_splitting}, and \eqref{eq:vorticity_splitting}, it is possible to split the compressible Euler equations, \eqref{eq:compressible_euler}, into horizontal and vertical components
		\begin{subequations}
			\label{eq:compressible_euler_splitting}
			\begin{align}
				\frac{\partial\upar}{\partial t} +(\omegaperp + \boldsymbol{f}_{\perp})\times\upar + \omegaparpar\times\uperp + \frac{1}{2}\gradpar\|\upar\|^{2} + \theta\gradpar\Pi &= 0, \label{eq:mom_par_split}\\
				\frac{\partial\uperp}{\partial t} + \omegaparperp\times\upar + \gradperp\left(\frac{1}{2}\|\uperp\|^{2} + gz\right) + \theta\gradperp\Pi &= 0, \label{eq:mom_perp_split} \\
				\frac{\partial\rho}{\partial t} + \nabla\cdot\left(\rho\upar\right) + \nabla\cdot\left(\rho\uperp\right) &= 0, \label{eq:dens_split} \\
				\frac{\partial (\rho\theta)}{\partial t} + \nabla\cdot\left(\rho\theta\upar\right) + \nabla\cdot\left(\rho\theta\uperp\right) &= 0, \label{eq:temp_split}
			\end{align}
		\end{subequations}
		where, as in \eqref{eq:vorticity_splitting}, $\omegaparpar := \curlpar\upar$ and $\omegaparperp := \curlpar\uperp$.

		The same splitting into horizontal and vertical components may be done for the basis functions, \secref{sec:three_d_mimetic_polynomial_basis_functions},
		\begin{equation}
			\wbasis_{j} = \left(\wbasis_{j}\right)_{\parallel} +  \left(\wbasis_{j}\right)_{\perp} := \wbasispar_{j} + \wbasisperp_{j}, \label{eq:wbasis_split}
		\end{equation}
		and
		\begin{equation}
			\ubasis_{j} = \left(\ubasis_{j}\right)_{\parallel} +  \left(\ubasis_{j}\right)_{\perp} := \ubasispar_{j} + \ubasisperp_{j}. \label{eq:ubasis_split}
		\end{equation}
		Recalling the definition of $\wbasis_{j}$ and $\ubasis_{j}$, \secref{sec:three_d_mimetic_polynomial_basis_functions}, we can explicitly write their horizontal and vertical components as
		\begin{align*}
			\wbasispar_{m}(\xi,\eta,\zeta; \degreepoly) &:=
			\begin{dcases}
				e_{i}^{\degreepoly}(\xi)l_{j}^{\degreepoly}(\eta)l_{k}^{\degreepoly}(\zeta)\ex,  \qquad m = i + j\degreepoly + k\degreepoly\degreepoly_{+} - 1, \qquad i = 1,\dots, \degreepoly, \quad j,k = 0, \cdots, \degreepoly, \\
				l_{i}^{\degreepoly}(\xi)e_{j}^{\degreepoly}(\eta)l_{k}^{\degreepoly}(\zeta)\ey,  \qquad m = i + j_{-}\degreepoly_{+} + k\degreepoly\degreepoly_{+} + \degreepoly\degreepoly_{+}^{2}, \qquad i,k = 0, \cdots, \degreepoly, \quad j = 1, \dots, \degreepoly,
			\end{dcases}\\ \\
			\wbasisperp_{m}(\xi,\eta,\zeta; \degreepoly) &:=
			l_{i}^{\degreepoly}(\xi)l_{j}^{\degreepoly}(\eta)e_{k}^{\degreepoly}(\zeta)\ez,  \qquad m = i + j\degreepoly_{+} + k_{-}\degreepoly_{+}^{2}, \qquad i,j = 0, \cdots, \degreepoly, \quad k = 1, \dots,\degreepoly\,\\ \\
			\ubasispar_{m}(\xi,\eta,\zeta; \degreepoly) &:=
			\begin{dcases}
				l_{i}^{\degreepoly}(\xi)e_{j}^{\degreepoly}(\eta)e_{k}^{\degreepoly}(\zeta)\ex,  \qquad m = i + j_{-}\degreepoly_{+} + k_{-}\degreepoly\degreepoly_{+}, \qquad i = 0,\dots, \degreepoly, \quad j,k = 1, \cdots, \degreepoly, \\
				e_{i}^{\degreepoly}(\xi)l_{j}^{\degreepoly}(\eta)e_{k}^{\degreepoly}(\zeta)\ey,  \qquad m = i_{-} + j\degreepoly + k_{-}\degreepoly\degreepoly_{+} + \degreepoly_{+}\degreepoly^{2}, \qquad i,k = 1, \cdots, \degreepoly, \quad j = 0, \dots, \degreepoly,
			\end{dcases} \\ \\
			\ubasisperp_{m}(\xi,\eta,\zeta; \degreepoly) &:=
			e_{i}^{\degreepoly}(\xi)e_{j}^{\degreepoly}(\eta)l_{k}^{\degreepoly}(\zeta)\ez,  \qquad m = i_{-} + j_{-}\degreepoly + k\degreepoly^{2} , \qquad i,j = 1, \cdots, \degreepoly, \quad k = 0, \dots,\degreepoly.
		\end{align*}
		Using this splitting of the basis functions we can also split the discrete function spaces such that \eqref{eq:function_spaces_span} becomes
		\begin{align}
			\wspace = \wspacepar\oplus\wspaceperp = \mathrm{span}\{\wbasispar_{1}, \dots, \wbasispar_{\wcardinalpar}\}\oplus\mathrm{span}\{\wbasisperp_{1}, \dots, \wbasisperp_{\wcardinalperp}\}, \\
			\uspace = \uspacepar\oplus\uspaceperp = \mathrm{span}\{\ubasispar_{1}, \dots, \ubasispar_{\ucardinalpar}\}\oplus\mathrm{span}\{\ubasisperp_{1}, \dots, \ubasisperp_{\ucardinalperp}\}, \label{eq:function_spaces_span_split}
		\end{align}
		where
		\begin{equation}
			\wcardinalpar := 2pp_{+}^{2}, \qquad \wcardinalperp := pp_{+}^{2}, \qquad \ucardinalpar := 2p_{+}p^{2}, \qquad \ucardinalperp = p_{+}p^{2}.
		\end{equation}

		\begin{remark}
			An important point relevant in the derivations that follow is that
			\begin{align*}
				\wbasispar_{j} &= \wbasis_{j},  \qquad j = 0, \dots, \wcardinalpar-1, \\
				\wbasisperp_{j} &= \wbasis_{j+\wcardinalpar}, \qquad j = 0, \dots, \wcardinalperp - 1, \\
				\ubasispar_{j} &= \ubasis_{j},  \qquad j = 0, \dots, \ucardinalpar-1, \\
				\ubasisperp_{j} &= \ubasis_{j+\ucardinalpar}, \qquad j = 0, \dots, \ucardinalperp - 1. \\
			\end{align*}
		\end{remark}
	\subsection{Unsplit discretization}
To discretize the unsplit form of the compressible Euler equations, \eqref{eq:compressible_euler},  we first introduce the weak
form: Given a domain $\Omega\subset\mathbb{R}^{3}$ and a Coriolis term $\boldsymbol{f} = f\ez \in H(\mathrm{curl},\Omega)$, find
$\boldsymbol{u}, \boldsymbol{U}, \boldsymbol{F}, \boldsymbol{P}\in\hdiv{\Omega}$,
$\boldsymbol{\omega}\in\hcurl{\Omega}$, and $\Pi,\rho,\theta,\Theta\in\ltwo{\Omega}$ for the prognostic equations
\begin{subequations}
\label{eq:euler_wf_prog}
\begin{align}
\begin{split}
\innerprod{\frac{\partial\bm u}{\partial t}}{\boldsymbol{\sigma}}_{\Omega} + \innerprod{(\bm\omega + \bm f)\times\bm u}{\boldsymbol{\sigma}}_{\Omega} +
\innerprod{\frac{1}{2}\|\boldsymbol{u}\|^2 + gz}{\nabla\cdot\boldsymbol{\sigma}}_{\Omega} + \\
\innerprod{\theta \boldsymbol{P}}{\boldsymbol{\sigma}}_{\Omega} &= 0, \qquad \forall \boldsymbol{\sigma} \in \hdiv{\Omega},
\end{split}
\label{eq::mom_weak}\\
\innerprod{\frac{\partial\rho}{\partial t}}{\alpha}_{\Omega} + \innerprod{\nabla\cdot\boldsymbol{U}}{\alpha}_{\Omega} &= 0, \qquad\forall\alpha\in\ltwo{\Omega}\,\label{eq::dens_weak}\\
\innerprod{\frac{\partial\Theta}{\partial t}}{\alpha}_{\Omega} + \innerprod{\nabla\cdot\boldsymbol{F}}{\alpha}_{\Omega} &= 0, \qquad\forall\alpha\in\ltwo{\Omega},\label{eq::temp_weak}
\end{align}
\end{subequations}
as well as the associated diagnostic equations
\begin{subequations}
\label{eq:euler_wf_diag}
\begin{align}
\innerprod{\Pi}{\nabla\cdot\boldsymbol{\sigma}}_{\Omega} - \innerprod{\boldsymbol{P}}{\boldsymbol{\sigma}}_{\Omega} &= 0, \qquad \forall\boldsymbol{\sigma}\in\hdiv{\Omega}\label{eq:Pi_P}\\
\innerprod{\rho\boldsymbol{u}}{\boldsymbol{\sigma}}_{\Omega} - \innerprod{\boldsymbol{U}}{\boldsymbol{\sigma}}_{\Omega} &= 0, \qquad \forall\boldsymbol{\sigma}\in\hdiv{\Omega}, \label{eq:rho_u_V}\\
\innerprod{\theta\boldsymbol{U}}{\boldsymbol{\sigma}}_{\Omega} - \innerprod{\boldsymbol{F}}{\boldsymbol{\sigma}}_{\Omega} &= 0, \qquad \forall\boldsymbol{\sigma}\in\hdiv{\Omega}, \label{eq:theta_V_F} \\
\innerprod{\boldsymbol{u}}{\nabla\times\boldsymbol{\beta}}_{\Omega} - \innerprod{\boldsymbol{\omega}}{\boldsymbol{\beta}}_{\Omega} & = 0, \qquad \forall\boldsymbol{\beta}\in\hcurl{\Omega}, \label{eq:curl_u_omega} \\
\innerprod{\Theta}{\alpha}_{\Omega}
- \innerprod{\rho\theta}{\alpha}_{\Omega} &= 0,
\qquad\forall\alpha\in\ltwo{\Omega},\label{eq:rho_theta_Theta}\\
\innerprod{\Pi}{\alpha}_{\Omega} - c_p\Bigg(\frac{R}{p_0}\Bigg)^{R/c_v}\innerprod{\Theta^{R/c_v}}{\alpha}_{\Omega} &= 0, \qquad\forall\alpha\in\ltwo{\Omega}.
\label{eq:eos_cont}
\end{align}
\end{subequations}

		Consider now the domain $\Omega\subset\mathbb{R}^{3}$ and its tessellation $\mathcal{T}(\Omega)$ consisting of $M$ arbitrary quadrilaterals (curved), $\Omega_{m}$, with $m = 1, \dots, M$. We assume that all quadrilateral elements $\Omega_{m}$ can be obtained from a map $\Phi_{m}: (\xi, \eta, \zeta) \in I^{3}\mapsto (\lambda,\phi,z)\in\Omega_{m}$. Then the pushforward $\Phi_{m,*}$ maps functions in the reference element $I^{3}$ to functions in
the physical element $\Omega_{m}$, see for example \cite{abraham_diff_geom, frankel}. For this reason it suffices to explore the analysis on the reference domain $I^{3}$. Additionally, the multi-element case follows the standard approach in finite elements.

		\begin{remark}
			If a differential geometry formulation was used, the physical quantities would be represented by differential k-forms and the map $\Phi_{m}: (\xi, \eta, \zeta) \in I^{3}\mapsto (\lambda,\phi,z)\in\Omega_{m}$ would generate a pullback, $\Phi_{m}^{*}$, mapping $k$-forms in physical space, $\Omega_{m}$, to $k$-forms in the reference element, $I^{3}$, \cite{palha2014}.
		\end{remark}

The discrete weak formulation can be stated as: Given $\Omega = I^{3}$, the polynomial degree $N$ and a Coriolis term
$\boldsymbol{f}_{h}\in\mathcal{W}_{h,\perp}(\Omega)$, for any time $t \in (0, t_{F}]$ find
$\boldsymbol{u}_{h}, \boldsymbol{U}_{h}, \boldsymbol{F}_{h}, \boldsymbol{P}_{h}\in\uspace(\Omega)$,
$\theta_h\in\mathcal{U}_{h,\perp}(\Omega)$,
$\boldsymbol{\omega}_{h}\in\wspace(\Omega)$, and $\Pi_{h},\rho_{h},\Theta_{h}\in\qspace(\Omega)$ such that

\begin{subequations}
\label{eq:euler_wf_discrete_prog}
\begin{align}
\begin{split}
\innerprod{\frac{\partial\boldsymbol{u}_{h}}{\partial t}}{\boldsymbol{\sigma}_{h}}_{\Omega} +
\innerprod{(\boldsymbol{\omega}_{h} + \boldsymbol{f}_{h})\times\boldsymbol{u}_{h}}{\boldsymbol{\sigma}_{h}}_{\Omega} +
\innerprod{\frac{1}{2}\|\boldsymbol{u}_h\|^2 + gz}{\nabla\cdot\boldsymbol{\sigma}_{h}}_{\Omega} + &\\
\innerprod{\theta_{h} \boldsymbol{P}_{h}}{\boldsymbol{\sigma}_{h}}_{\Omega} &= 0, \qquad \forall\boldsymbol{\sigma}_{h}\in\uspace(\Omega),
\end{split}
\label{eq::mom_weak_discrete}\\
\innerprod{\frac{\partial\rho_{h}}{\partial t}}{\alpha_{h}}_{\Omega} + \innerprod{\nabla\cdot\boldsymbol{U}_{h}}{\alpha_{h}}_{\Omega} &= 0, \qquad\forall\alpha_{h}\in\qspace(\Omega)\,\label{eq::dens_weak_discrete}\\
\innerprod{\frac{\partial\Theta_{h}}{\partial t}}{\alpha_{h}}_{\Omega} + \innerprod{\nabla\cdot\boldsymbol{F}_{h}}{\alpha_{h}}_{\Omega} &= 0, \qquad\forall\alpha_{h}\in\qspace(\Omega),\label{eq::temp_weak_discrete}
\end{align}
\end{subequations}
and
\begin{subequations}
\label{eq:euler_wf_discrete_diag}
\begin{align}
\innerprod{\Pi_{h}}{\nabla\cdot\boldsymbol{\sigma}_{h}}_{\Omega} - \innerprod{\boldsymbol{P}_{h}}{\boldsymbol{\sigma}_{h}}_{\Omega} &= 0, \qquad \forall\boldsymbol{\sigma}_{h}\in\uspace(\Omega)\label{eq:Pi_P_discrete}\\
\innerprod{\rho_{h}\boldsymbol{u}_{h}}{\boldsymbol{\sigma}_{h}}_{\Omega} - \innerprod{\boldsymbol{U}_{h}}{\boldsymbol{\sigma}_{h}}_{\Omega} &= 0, \qquad \forall\boldsymbol{\sigma}_{h}\in\uspace(\Omega), \label{eq:rho_u_V_discrete}\\
\innerprod{\theta_{h}\boldsymbol{U}_{h}}{\boldsymbol{\sigma}_{h}}_{\Omega} - \innerprod{\boldsymbol{F}_{h}}{\boldsymbol{\sigma}_{h}}_{\Omega} &= 0, \qquad \forall\boldsymbol{\sigma}_{h}\in\uspace(\Omega), \label{eq:theta_V_F_discrete} \\
\innerprod{\boldsymbol{u}_{h}}{\nabla\times\boldsymbol{\beta}_{h}}_{\Omega} - \innerprod{\boldsymbol{\omega}_{h}}{\boldsymbol{\beta}_{h}}_{\Omega} & = 0, \qquad \forall\boldsymbol{\beta}_{h}\in\wspace(\Omega), \label{eq:curl_u_omega_discrete} \\
\innerprod{\Theta_{h}}{\sigma_{h,\perp}}_{\Omega}
- \innerprod{\rho_{h}\theta_{h}}{\sigma_{h,\perp}}_{\Omega} &= 0,
\qquad\forall\alpha_{h}\in\mathcal{U}_{h,\perp}(\Omega),\label{eq:rho_theta_Theta_discrete}\\
\innerprod{\Pi_h}{\alpha_h}_{\Omega} - c_p\Bigg(\frac{R}{p_0}\Bigg)^{R/c_v}\innerprod{\Theta_h^{R/c_v}}{\alpha_h}_{\Omega} &= 0, \qquad\forall\alpha_h\in\qspace(\Omega).
\end{align}
\end{subequations}

Using the expansions for all unknowns, \eqref{eq:euler_wf_discrete_prog}, \eqref{eq:euler_wf_discrete_diag}
may be written as: Find $\boldsymbol{\mathsf{u}},
\boldsymbol{\mathsf{U}}, \boldsymbol{\mathsf{F}}, \boldsymbol{\mathsf{P}}\in\mathbb{R}^{\ucardinal}$,
$\mathsf{\theta}\in\mathbb{R}^{d_{\mathcal{U}_{\perp}}}$,
$\boldsymbol{\mathsf{\omega}}\in\mathbb{R}^{\wcardinal}$, and
$\mathsf{\Pi},\mathsf{\rho},\mathsf{\Theta}\in\mathbb{R}^{\qcardinal}$ such that
\begin{subequations}
\label{eq:euler_wf_discrete_algebraic_prog}
\begin{align}
\begin{split}
\sum_{i=0}^{\ucardinal-1}\innerprod{\ubasis_{i}}{\ubasis_{j}}_{\Omega}\frac{\mathrm{d}\mathsf{u}_{i}}{\mathrm{d} t} + \sum_{i=0}^{\ucardinal-1}\innerprod{(\boldsymbol{\omega}_{h} + \boldsymbol{f}_{h})\times\ubasis_{i}}{\ubasis_{j}}_{\Omega}\mathsf{u}_{i} + &\\
\sum_{i=0}^{\qcardinal-1}\innerprod{\frac{1}{2}\boldsymbol{u}_{h} \cdot \ubasis_{i}}{\nabla\cdot\ubasis_{j}}_{\Omega}\mathsf{u}_{i} +
\sum_{i=0}^{\qcardinal-1}\innerprod{g}{\nabla\cdot\ubasis_{j}}_{\Omega}\mathsf{z}_{i} + &\\
\sum_{i=0}^{\ucardinal-1}\innerprod{\theta_{h} \ubasis_{i}}{\ubasis_{j}}_{\Omega}\mathsf{P}_{i} &= 0, \qquad j=0,\dots,\ucardinal - 1,
					\end{split} \label{eq::mom_weak_discrete_algebraic}\\
				\sum_{i=0}^{\qcardinal-1}\innerprod{\qbasis_{i}}{\qbasis_{j}}_{\Omega}\frac{\mathrm{d}\mathsf{\rho}_{i}}{\mathrm{d}t} + \sum_{i=0}^{\ucardinal-1}\innerprod{\nabla\cdot\ubasis_{i}}{\qbasis_{j}}_{\Omega}\mathsf{U}_{i} &= 0, \qquad j=0,\dots,\qcardinal-1\,\label{eq::dens_weak_discrete_algebraic}\\
				\sum_{i=0}^{\qcardinal-1}\innerprod{\qbasis_{i}}{\qbasis_{j}}_{\Omega}\frac{\mathrm{d}\mathsf{\Theta}_{i}}{\mathrm{d}t} + \sum_{j=0}^{\ucardinal-1}\innerprod{\nabla\cdot\ubasis_{i}}{\qbasis_{j}}_{\Omega}\mathsf{F}_{i} &= 0, \qquad j = 0, \dots, \qcardinal-1,\label{eq::temp_weak_discrete_algebraic}
\end{align}
\end{subequations}
and
\begin{subequations}
\label{eq:euler_wf_discrete_algebraic_diag}
\begin{align}
\sum_{i=0}^{\qcardinal-1}\innerprod{\qbasis_{i}}{\nabla\cdot\ubasis_{j}}_{\Omega} \mathsf{\Pi}_{i} - \sum_{i=0}^{\ucardinal-1}\innerprod{\ubasis_{i}}{\ubasis_{j}}_{\Omega}\mathsf{P}_{i} &= 0, \qquad j = 0, \dots, \ucardinal-1\label{eq:Pi_P_discrete_algebraic}\\
\sum_{i=0}^{\ucardinal}\innerprod{\rho_{h}\ubasis_{i}}{\ubasis_{j}}_{\Omega}\mathsf{u}_{i} - \sum_{i=0}^{\ucardinal}\innerprod{\ubasis_{i}}{\ubasis_{j}}_{\Omega}\mathsf{U}_{i} &= 0, \qquad j = 0, \dots, \ucardinal-1, \label{eq:rho_u_V_discrete_algebraic}\\
\sum_{i=0}^{\ucardinal-1}\innerprod{\theta_{h}\ubasis_{i}}{\ubasis_{j}}_{\Omega}\mathsf{U}_{i} - \sum_{i=0}^{\ucardinal-1}\innerprod{\ubasis_{i}}{\ubasis_{j}}_{\Omega}\mathsf{F}_{i} &= 0, \qquad j = 0, \dots, \ucardinal-1, \label{eq:theta_V_F_discrete_algebraic} \\
\sum_{i=0}^{\ucardinal-1}\innerprod{\ubasis_{i}}{\nabla\times\wbasis_{j}}_{\Omega}\mathsf{u}_{i} - \sum_{i=0}^{\wcardinal-1}\innerprod{\wbasis_{i}}{\wbasis_{j}}_{\Omega}\mathsf{\omega}_{i} & = 0, \qquad j = 0, \dots, \wcardinal-1, \label{eq:curl_u_omega_discrete_algebraic} \\
\begin{split}
\sum_{i=0}^{\qcardinal-1} \innerprod{\qbasis_{i}}{\epsilon^{\mathcal{U}_{\perp}}_{j}}_{\Omega}\mathsf{\Theta}_{i}
- \sum_{i=0}^{d_{\mathsf{U}_{\perp}}-1}\innerprod{\rho_{h}\epsilon^{\mathcal{U}_{\perp}}_{i}}{\epsilon^{\mathcal{U}_{\perp}}_{j}}_{\Omega}\mathsf{\theta}_{i}
&= 0, \qquad j = 0, \dots, d_{\mathcal{U}_{\perp}}-1,
\end{split}\label{eq:rho_theta_Theta_discrete_algebraic} \\
c_p\Big(\frac{R}{p_0}\Big)^{R/c_v}\sum_{i=0}^{\ucardinal-1}\innerprod{(\qbasis_{i}\mathsf{\Theta}_{i})^{R/c_v}}{\qbasis_{j}}_{\Omega} -
\sum_{i=0}^{\qcardinal-1}\innerprod{\qbasis_{i}}{\qbasis_{j}}_{\Omega}\mathsf{\Pi}_{i} &= 0, \qquad j = 0, \dots,
\qcardinal-1.\label{eq:Pi_Theta_discrete_algebraic}
\end{align}
\end{subequations}

	\subsection{Split discretization}

It is important to note that functions in $\mathcal{W}_{\perp,h}$ and $\mathcal{U}_{\parallel,h}$ are
discontinuous across vertical element boundaries, while those in $\mathcal{U}_{\perp,h}$ are discontinuous
across horizontal boundaries. Similarly functions in $\mathcal{Q}_h$ are discontinuous across both vertical
and horizontal boundaries. These properties allow us to split the three dimensional problem presented in
\eqref{eq:euler_wf_discrete_algebraic_prog}, \eqref{eq:euler_wf_discrete_algebraic_diag} into separate horizontal and vertical problems,
and in doing so avoid solving any global three dimensional implicit systems.

The split discretization is obtained by using \eqref{eq:vector_splitting}, \eqref{eq:grad_splitting}, \eqref{eq:curl_splitting}, and \eqref{eq:vorticity_splitting}, together with \eqref{eq:wbasis_split} and \eqref{eq:ubasis_split} in \eqref{eq:euler_wf_discrete_algebraic_prog}, \eqref{eq:euler_wf_discrete_algebraic_diag} in order to obtain a discrete version of the split Euler equations \eqref{eq:compressible_euler_splitting}. The horizontal discrete equations are given by
		\begin{subequations}
			\label{eq:euler_wf_discrete_algebraic_par_split_incidence}
			\begin{align}
	\begin{split}
\sum_{i=0}^{\ucardinalpar-1}\innerprod{\ubasispar_{i}}{\ubasispar_{j}}_{\Omega}\frac{\mathrm{d}\mathsf{u}_{i,\parallel}}{\mathrm{d} t} + \sum_{i=0}^{\ucardinalpar-1}\innerprod{(\boldsymbol{\omega}_{h,\perp} + \boldsymbol{f}_{h,\perp})\times\ubasispar_{i}}{\ubasispar_{j}}_{\Omega}\mathsf{u}_{i,\parallel} + & \\
\sum_{i=0}^{\ucardinalperp-1}\innerprod{\boldsymbol{\omega}_{h,\parallel,\parallel}\times\ubasisperp_{i}}{\ubasispar_{j}}_{\Omega}\mathsf{u}_{i,\perp} + & \\
\sum_{i,k=0}^{\qcardinal-1,\qcardinal-1}\left(\mathsf{E}^{3,2}_{\parallel}\right)_{j,k}^{\top}\innerprod{\frac{1}{2}\boldsymbol{u}_{h,\parallel} \cdot \ubasispar_{i}}{\qbasis_{k}}_{\Omega}\mathsf{u}_{i,\parallel}  + &\\
\sum_{i=0}^{\ucardinalpar-1}\innerprod{\theta_{h} \ubasispar_{i}}{\ubasispar_{j}}_{\Omega}\mathsf{P}_{i,\parallel} &= 0, \qquad j=0,\dots,\ucardinalpar-1,
	\end{split} \label{eq::mom_weak_discrete_algebraic_par_split_incidence}\\
\sum_{i,k=0}^{\qcardinal-1,\qcardinal-1}\left(\mathsf{E}^{3,2}_{\parallel}\right)_{j,k}^{\top}\innerprod{\qbasis_{i}}{\qbasis_{k}}_{\Omega} \mathsf{\Pi}_{i}  - \sum_{i=0}^{\ucardinalpar-1}\innerprod{\ubasispar_{i}}{\ubasispar_{j}}_{\Omega}\mathsf{P}_{i,\parallel} &= 0, \qquad j = 0, \dots, \ucardinalpar-1\label{eq:Pi_P_discrete_algebraic_par_split_incidence}\\
				\sum_{i=0}^{\ucardinalpar-1}\innerprod{\rho_{h}\ubasispar_{i}}{\ubasispar_{j}}_{\Omega}\mathsf{u}_{i,\parallel} - \sum_{i=0}^{\ucardinalpar-1}\innerprod{\ubasispar_{i}}{\ubasispar_{j}}_{\Omega}\mathsf{U}_{i,\parallel} &= 0, \qquad j = 0, \dots, \ucardinalpar-1, \label{eq:rho_u_V_discrete_algebraic_par_split_incidence}\\
				\sum_{i=0}^{\ucardinalpar-1}\innerprod{\theta_{h}\ubasispar_{i}}{\ubasispar_{j}}_{\Omega}\mathsf{U}_{i,\parallel} - \sum_{i=0}^{\ucardinalpar-1}\innerprod{\ubasispar_{i}}{\ubasispar_{j}}_{\Omega}\mathsf{F}_{i,\parallel} &= 0, \qquad j = 0, \dots, \ucardinalpar-1, \label{eq:theta_V_F_discrete_algebraic_par_split_incidence} \\
\sum_{i,k=0}^{\ucardinalpar-1,\ucardinalpar-1}\left(\mathsf{E}^{2,1}_{\parallel,\parallel}\right)_{j,k}^{\top}\innerprod{\ubasispar_{i}}{\ubasispar_{k}}_{\Omega}\mathsf{u}_{i,\parallel} - \sum_{i=0}^{\wcardinalpar-1}\innerprod{\wbasispar_{i}}{\wbasispar_{j}}_{\Omega}\mathsf{\omega}_{i,\parallel,\parallel} & = 0, \qquad j = 0, \dots, \wcardinalpar-1, \label{eq:curl_u_omega_discrete_algebraic_par_par_split_incidence} \\
\sum_{i,k=0}^{\ucardinalpar-1,\ucardinalpar-1}\left(\mathsf{E}^{2,1}_{\parallel,\perp}\right)_{j,k}^{\top}\innerprod{\ubasispar_{i}}{\ubasispar_{k}}_{\Omega}\mathsf{u}_{i,\parallel} - \sum_{i=0}^{\wcardinalperp-1}\innerprod{\wbasisperp_{i}}{\wbasisperp_{j}}_{\Omega}\mathsf{\omega}_{i,\perp} & = 0, \qquad j = 0, \dots, \wcardinalperp-1, \label{eq:curl_u_omega_discrete_algebraic_perp_split_incidence}
			\end{align}
		\end{subequations}
		where we have introduced $\omega_{h,\parallel,\parallel} := \sum_{i=0}^{\wcardinalpar-1}\mathsf{\omega}_{i,\parallel,\parallel} \wbasispar_{i}$, $\omega_{h,\parallel,\perp} := \sum_{i=0}^{\wcardinalpar-1} \mathsf{\omega}_{i,\parallel,\perp}\wbasispar_{i}$. In the same way, the vertical discrete equations are
		\begin{subequations}
			\label{eq:euler_wf_discrete_algebraic_perp_split_incidence}
			\begin{align}
					\begin{split}
\sum_{i=0}^{\ucardinalperp-1}\innerprod{\ubasisperp_{i}}{\ubasisperp_{j}}_{\Omega}\frac{\mathrm{d}\mathsf{u}_{i,\perp}}{\mathrm{d} t} + \sum_{i=0}^{\ucardinalpar-1}\innerprod{\boldsymbol{\omega}_{h,\parallel,\perp}\times\ubasispar_{i}}{\ubasisperp_{j}}_{\Omega}\mathsf{u}_{i,\parallel} + & \\
\sum_{i,k=0}^{\qcardinal-1,\qcardinal-1}\left(\mathsf{E}^{3,2}_{\perp}\right)_{j,k}^{\top}\innerprod{\frac{1}{2}{\boldsymbol{u}}_{h,\perp}\cdot\ubasisperp_{i}}{\qbasis_{k}}_{\Omega}\mathsf{u}_{i,\perp} + &\\
\sum_{i,k=0}^{\qcardinal-1,\qcardinal-1}g\left(\mathsf{E}^{3,2}_{\perp}\right)_{j,k}^{\top}\innerprod{\qbasis_i}{\qbasis_{k}}_{\Omega}\mathsf{z}_{i} +
\sum_{i=0}^{\ucardinalperp-1}\innerprod{\theta_{h} \ubasisperp_{i}}{\ubasisperp_{j}}_{\Omega}\mathsf{P}_{i,\perp} &= 0, \qquad j=0,\dots,\ucardinalperp-1,
					\end{split} \label{eq::mom_weak_discrete_algebraic_perp_split_incidence}\\
\sum_{i,k=0}^{\qcardinal-1,\qcardinal-1}\left(\mathsf{E}^{3,2}_{\perp}\right)_{j,k}^{\top}\innerprod{\qbasis_{i}}{\qbasis_{j}}_{\Omega} \mathsf{\Pi}_{i} - \sum_{i=0}^{\ucardinalperp-1}\innerprod{\ubasisperp_{i}}{\ubasisperp_{j}}_{\Omega}\mathsf{P}_{i,\perp} &= 0, \qquad j = 0, \dots, \ucardinalperp-1\label{eq:Pi_P_discrete_algebraic_perp_split_incidence}\\
				\sum_{i=0}^{\ucardinalperp-1}\innerprod{\rho_{h}\ubasisperp_{i}}{\ubasisperp_{j}}_{\Omega}\mathsf{u}_{i,\perp} - \sum_{i=0}^{\ucardinalperp-1}\innerprod{\ubasisperp_{i}}{\ubasisperp_{j}}_{\Omega}\mathsf{U}_{i,\perp} &= 0, \qquad j = 0, \dots, \ucardinalperp-1, \label{eq:rho_u_V_discrete_algebraic_perp_split_incidence}\\
				\sum_{i=0}^{\ucardinalperp-1}\innerprod{\theta_{h}\ubasisperp_{i}}{\ubasisperp_{j}}_{\Omega}\mathsf{U}_{i,\perp} - \sum_{i=0}^{\ucardinalperp-1}\innerprod{\ubasisperp_{i}}{\ubasisperp_{j}}_{\Omega}\mathsf{F}_{i,\perp} &= 0, \qquad j = 0, \dots, \ucardinalperp-1, \label{eq:theta_V_F_discrete_algebraic_perp_split_incidence} \\
\sum_{i,k=0}^{\ucardinalperp-1,\ucardinalperp-1}\left(\mathsf{E}^{2,1}_{\perp}\right)_{j,k}^{\top}\innerprod{\ubasisperp_{i}}{\ubasisperp_{k}}_{\Omega}\mathsf{u}_{i,\perp} - \sum_{i=0}^{\wcardinalpar-1}\innerprod{\wbasispar_{i}}{\wbasispar_{j}}_{\Omega}\mathsf{\omega}_{i,\parallel,\perp} & = 0, \qquad j = 0, \dots, \wcardinalpar-1, \label{eq:curl_u_omega_discrete_algebraic_par_perp_split_incidence}
			\end{align}
		\end{subequations}
where we have introduced $\omega_{h,\perp} := \sum_{i=0}^{\wcardinalperp-1}\mathsf{\omega}_{i,\perp} \wbasisperp_{i}$.
Note that for efficiency we do not assemble the second term in \eqref{eq::mom_weak_discrete_algebraic_perp_split_incidence}
for simulations at hydrostatic resolutions
(and as such we do not apply the diagnostic equation \eqref{eq:curl_u_omega_discrete_algebraic_par_perp_split_incidence}
for $\omega_{h,\perp}$). This term represents the horizontal advection of vertical velocity, and as such is not
significant for hydrostatic motions where the vertical scales are small with respect to the horizontal scales, and so
the horizontal gradients of vertical terms have minimal impact on the dynamics. While this term is also omitted for
hydrostatic formulations \cite{White05}, our motivation here is based on simple scale analysis. Additionally the rotational terms
have no projection onto the energy of the system within the mimetic discretization \cite{LPG18} (though they do re-arrange
kinetic energy). At non-hydrostatic resolutions however we do include this term, since the horizontal and vertical
scales are closer to unity, and failure to incorporate this term means that vertical motions are not being properly
transported with the horizontal flow, such that upstream vertical oscillations may be exaggerated.

The only equations in the above system that cannot be effectively split between the horizontal and vertical systems
are \eqref{eq:curl_u_omega_discrete_algebraic_par_par_split_incidence} and \eqref{eq:curl_u_omega_discrete_algebraic_par_perp_split_incidence},
since these involve vertical gradients of a horizontally continuous field and vice versa. These terms are required for the 
vertical advection of horizontal velocity and the horizontal advection of vertical velocity respectively. However by
employing a horizontal velocity space which is piecewise constant in the vertical, and a vertical velocity that is
piecewise linear in the vertical, we can avoid the need to diagnose these terms through the use of global matrices
by employing a direct differencing at the vertical layer interfaces.

Additionally, we also have the flux form equations for density and density weighted potential temperature transport that
contain both vertical and horizontal components. While we have not included these in the split systems described in
\eqref{eq:euler_wf_discrete_algebraic_par_split_incidence} and
\eqref{eq:euler_wf_discrete_algebraic_perp_split_incidence}, since doing so incurs a temporal splitting
error, in practice these equations are also split between their horizontal and vertical components. For the Strang carryover
scheme detailed in Section 5.3 this results in a second order temporal error for the full system. These equations are
given as
		\begin{subequations}
			\label{eq:euler_wf_discrete_algebraic_no_par_no_perp_split_incidence}
			\begin{align}
			\begin{split}
\sum_{i=0}^{\qcardinal-1}\innerprod{\qbasis_{i}}{\qbasis_{j}}_{\Omega}\frac{\mathrm{d}\mathsf{\rho}_{i}}{\mathrm{d}t} + \sum_{i,k=0}^{\ucardinalpar-1,\qcardinal-1}\innerprod{\qbasis_{k}}{\qbasis_{j}}_{\Omega}\left(\mathsf{E}^{3,2}_{\parallel}\right)_{k,i}\mathsf{U}_{i,\parallel} + & \\
\sum_{i,k=0}^{\ucardinalperp-1,\qcardinal-1}\innerprod{\qbasis_{k}}{\qbasis_{j}}_{\Omega}\left(\mathsf{E}^{3,2}_{\perp}\right)_{k,i}\mathsf{U}_{i,\perp} &= 0, \qquad j=0,\dots,\qcardinal-1\,
				\end{split} \label{eq::dens_weak_discrete_algebraic_split_incidence} \\
				\begin{split}
\sum_{i=0}^{\qcardinal-1}\innerprod{\qbasis_{i}}{\qbasis_{j}}_{\Omega}\frac{\mathrm{d}\mathsf{\Theta}_{i}}{\mathrm{d}t} + \sum_{i,k=0}^{\ucardinalpar-1,\qcardinal-1}\innerprod{\qbasis_{k}}{\qbasis_{j}}_{\Omega}\left(\mathsf{E}^{3,2}_{\parallel}\right)_{k,i}\mathsf{F}_{i,\parallel} + &\\
\sum_{i,k=0}^{\ucardinalperp-1,\qcardinal-1}\innerprod{\qbasis_{k}}{\qbasis_{j}}_{\Omega}\left(\mathsf{E}^{3,2}_{\perp}\right)_{k,i}\mathsf{F}_{i,\perp} &= 0, \qquad j = 0, \dots, \qcardinal-1,
				\end{split} \label{eq::temp_weak_discrete_algebraic_split_incidence}
			\end{align}
		\end{subequations}
Note that the diagnostic equations for potential temperature \eqref{eq:rho_theta_Theta_discrete_algebraic}
and the equation of state \eqref{eq:Pi_Theta_discrete_algebraic} are also included in both the horizontal
and vertical systems. \blue{We also note that the mass matrices all cancel in \eqref{eq:euler_wf_discrete_algebraic_no_par_no_perp_split_incidence}, 
such that these flux form transport equations may effectively be solved in the strong form \cite{LPG18,LP18} 
for the point-wise conservation of mass and mass weighted potential temperature.}

		\begin{remark}
			In \eqref{eq:euler_wf_discrete_algebraic_par_split_incidence},  \eqref{eq:euler_wf_discrete_algebraic_perp_split_incidence}, and \eqref{eq:euler_wf_discrete_algebraic_no_par_no_perp_split_incidence}, we have used the fact that the incidence matrices can be written as
		\begin{equation}
		\arraycolsep=2.0pt\def\arraystretch{2.2}
			\incidencegrad =
			\left[
				\begin{array}{c}
					\incidencegradpar \\
					\incidencegradperp
				\end{array}
			\right], \qquad
			\incidencecurl =
			\left[
				\begin{array}{cc}
					\incidencecurlparpar & \incidencecurlparperp \\
					\incidencecurlperp & \boldsymbol{\mathsf{0}}
				\end{array}
			\right], \qquad \mathrm{and}\qquad
			\incidencediv =
			\left[
				\begin{array}{cc}
					\incidencedivpar & \incidencedivperp
				\end{array}
			\right], \label{eq:incidence_matrices_splitting}
		\end{equation}
		where $\incidencegradpar$ is a $\wcardinalpar\times\pcardinal$ matrix, $\incidencegradperp$ is a $\wcardinalperp\times\pcardinal$ matrix, $\incidencecurlparpar$ is a $\ucardinalpar\times\wcardinalpar$ matrix, $\incidencecurlparperp$ is a $\ucardinalpar\times\wcardinalperp$ matrix, $\incidencecurlperp$ is a $\ucardinalperp\times\wcardinalpar$ matrix, $\incidencedivpar$ is a $\qcardinal\times\ucardinalpar$ matrix, and $\incidencedivperp$ is a $\qcardinal\times\ucardinalperp$ matrix.
		\end{remark}

		\begin{remark}
			In \eqref{eq:euler_wf_discrete_algebraic_par_split_incidence},  \eqref{eq:euler_wf_discrete_algebraic_perp_split_incidence}, and \eqref{eq:euler_wf_discrete_algebraic_no_par_no_perp_split_incidence}, two important points to note are
			\[
				\mathsf{u}_{i} = \mathsf{u}_{i,\parallel}, \quad i = 0, \dots, \ucardinalpar - 1 \qquad \mathrm{and}\qquad \mathsf{u}_{i + \ucardinalpar} = \mathsf{u}_{i,\perp}, \quad i = 0, \dots, \ucardinalperp - 1
			\]
			and
			\[
				\boldsymbol{\omega}_{h} = \omega_{h,\parallel,\parallel} + \omega_{h,\parallel,\perp} + \omega_{h,\perp}  =   \sum_{i=0}^{\wcardinalpar-1}\mathsf{\omega}_{i,\parallel,\parallel} \wbasispar_{i} + \sum_{i=0}^{\wcardinalpar-1} \mathsf{\omega}_{i,\parallel,\perp}\wbasispar_{i} + \sum_{i=0}^{\wcardinalperp-1}\mathsf{\omega}_{i,\perp} \wbasisperp_{i}.
			\]
		\end{remark}

		In compact matrix notation we can write \eqref{eq:euler_wf_discrete_algebraic_par_split_incidence} as
		\begin{subequations}
			\label{eq:euler_wf_discrete_algebraic_par_split_incidence_matrix}
			\begin{align}
\massupar\frac{\mathrm{d}\ualgebraicpar}{\mathrm{d}t} + \boldsymbol{\mathsf{R}}^{\parallel,\parallel}\ualgebraicpar +
\boldsymbol{\mathsf{R}}^{\parallel, \perp}\ualgebraicperp +
\left(\boldsymbol{\mathsf{E}}^{3,2}_{\parallel}\right)^{\top}\boldsymbol{\mathsf{T}}^{\mathcal{U}_{\parallel}}\ualgebraicpar +
\boldsymbol{\mathsf{S}}^{\mathcal{U}_{\parallel}}\boldsymbol{\mathsf{P}}^{\parallel} &= \boldsymbol{\mathsf{0}},\label{eq::mom_weak_discrete_algebraic_par_split_incidence_matrix}\\
				\left(\boldsymbol{\mathsf{E}}^{3,2}_{\parallel}\right)^{\top}\massq\boldsymbol{\mathsf{\Pi}} - \massupar\boldsymbol{\mathsf{P}}^{\parallel} &= \boldsymbol{\mathsf{0}}, \label{eq:Pi_P_discrete_algebraic_par_split_incidence_matrix}\\
				\nupar\ualgebraicpar - \massupar\boldsymbol{\mathsf{U}}^{\parallel} & = \boldsymbol{\mathsf{0}}, \label{eq:rho_u_V_discrete_algebraic_par_split_incidence_matrix}\\
\boldsymbol{\mathsf{S}}^{\mathcal{U}_{\parallel}}\boldsymbol{\mathsf{U}}^{\parallel} - \massupar\boldsymbol{\mathsf{F}}^{\parallel} &= \boldsymbol{\mathsf{0}}, \label{eq:theta_V_F_discrete_algebraic_par_split_incidence_matrix}\\
				\left(\boldsymbol{\mathsf{E}}^{2,1}_{\parallel,\parallel}\right)^{\top}\massupar\ualgebraicpar - \masswpar\boldsymbol{\mathsf{\omega}}^{\parallel,\parallel} &= \boldsymbol{\mathsf{0}}, \label{eq:curl_u_omega_discrete_algebraic_par_split_incidence_matrix}\\
				\left(\boldsymbol{\mathsf{E}}^{2,1}_{\perp}\right)^{\top}\massuperp\ualgebraicperp - \masswpar\boldsymbol{\mathsf{\omega}}^{\parallel,\perp} &= \boldsymbol{\mathsf{0}},
			\end{align}
		\end{subequations}
		with
		\[
			\mathsf{M}^{\mathcal{U}_{\parallel}}_{ij} := \innerprod{\ubasispar_{j}}{\ubasispar_{i}}_{\Omega}, \qquad \mathsf{R}^{\parallel,\parallel}_{ij} := \innerprod{(\boldsymbol{\omega}_{h,\perp} + \boldsymbol{f}_{h,\perp})\times\ubasispar_{j}}{\ubasispar_{i}}_{\Omega}, \qquad \mathsf{R}^{\parallel,\perp}_{ij} := \innerprod{\boldsymbol{\omega}_{h,\parallel,\parallel}\times\ubasisperp_{j}}{\ubasispar_{i}}_{\Omega}, \quad
		\]
		\[
\mathsf{M}^{\mathcal{Q}}_{ij} := \innerprod{\qbasis_{j}}{\qbasis_{i}}_{\Omega},\qquad\mathsf{S}^{\mathcal{U}_{\parallel}}_{ij} := \innerprod{\theta_{h} \ubasispar_{j}}{\ubasispar_{i}}_{\Omega}, \qquad \mathsf{N}^{\mathcal{U}_{\parallel}}_{ij} := \innerprod{\rho_{h}\ubasispar_{j}}{\ubasispar_{i}}_{\Omega},
		\]
		\[
\mathsf{M}^{\mathcal{W}_{\parallel}}_{ij} := \innerprod{\wbasispar_{j}}{\wbasispar_{i}}_{\Omega},
\qquad\mathsf{M}^{\mathcal{U}_{\perp}}_{ij} := \innerprod{\ubasisperp_{j}}{\ubasisperp_{i}}_{\Omega},
\qquad\mathsf{T}^{\mathcal{U}_{\parallel}}_{ij} := \innerprod{\frac{1}{2}\boldsymbol{u}_{h,\parallel} \cdot \ubasispar_{j}}{\qbasis_{i}}_{\Omega}.
		\]
		In a similar manner, the vertical equations, \eqref{eq:euler_wf_discrete_algebraic_perp_split_incidence}, can be written in compact matrix notation as
		\begin{subequations}
			\label{eq:euler_wf_discrete_algebraic_perp_split_incidence_matrix}
			\begin{align}
\massuperp\frac{\mathrm{d}\ualgebraicperp}{\mathrm{d}t} + \boldsymbol{\mathsf{R}}^{\perp,\parallel}\ualgebraicpar +
\left(\boldsymbol{\mathsf{E}}^{3,2}_{\perp}\right)^{\top}\boldsymbol{\mathsf{T}}^{\mathcal{U}_{\perp}}\ualgebraicperp +
g\left(\boldsymbol{\mathsf{E}}^{3,2}_{\perp}\right)^{\top}\massq\boldsymbol{\mathsf{z}} +
\boldsymbol{\mathsf{S}}^{\mathcal{U}_{\perp}}\boldsymbol{\mathsf{P}}^{\perp} &= \boldsymbol{\mathsf{0}},\label{eq::mom_weak_discrete_algebraic_perp_split_incidence_matrix}\\
				\left(\boldsymbol{\mathsf{E}}_{\perp}^{3,2}\right)^{\top}\massq\boldsymbol{\mathsf{\Pi}} - \massuperp\boldsymbol{\mathsf{P}}^{\perp} &= \boldsymbol{\mathsf{0}},\label{eq:Pi_P_discrete_algebraic_perp_split_incidence_matrix}\\
				\nuperp\ualgebraicperp - \massuperp\boldsymbol{\mathsf{U}}^{\perp} &= \boldsymbol{\mathsf{0}},\label{eq:rho_u_V_discrete_algebraic_perp_split_incidence_matrix}\\
\boldsymbol{\mathsf{S}}^{\mathcal{U}_{\perp}}\boldsymbol{\mathsf{U}}^{\perp} - \massuperp\boldsymbol{\mathsf{F}}^{\perp} &= \boldsymbol{\mathsf{0}}, \label{eq:theta_V_F_discrete_algebraic_perp_split_incidence_matrix} \\
				\left(\boldsymbol{\mathsf{E}}^{2,1}_{\parallel,\perp}\right)^{\top}\massuperp\ualgebraicperp - \masswperp\boldsymbol{\mathsf{\omega}}^{\perp} &= \boldsymbol{\mathsf{0}},\label{eq:curl_u_omega_discrete_algebraic_perp_split_incidence_matrix},
			\end{align}
		\end{subequations}
		with
		\[
\mathsf{R}^{\perp,\parallel}_{ij} := \innerprod{(\boldsymbol{\omega}_{h,\perp} + \boldsymbol{f}_{h,\perp})\times\ubasispar_{j}}{\ubasisperp_{i}}_{\Omega},
\qquad
\mathsf{S}^{\mathcal{U}_{\perp}}_{ij} := \innerprod{\theta_{h} \ubasisperp_{j}}{\ubasisperp_{i}}_{\Omega}, \qquad\mathsf{N}^{\mathcal{U}_{\perp}}_{ij} := \innerprod{\rho_{h}\ubasisperp_{j}}{\ubasisperp_{i}}_{\Omega},
		\]
		\[
\mathsf{M}^{\mathcal{W}_{\perp}}_{ij} := \innerprod{\wbasisperp_{j}}{\wbasisperp_{i}}_{\Omega},\qquad
\mathsf{T}^{\mathcal{U}_{\perp}}_{ij} := \innerprod{\frac{1}{2}\boldsymbol{u}_{h,\perp} \cdot \ubasisperp_{j}}{\qbasis_{i}}_{\Omega}.
		\]
In \eqref{eq::mom_weak_discrete_algebraic_perp_split_incidence_matrix} the vector $\boldsymbol{\mathsf{z}}$ represents the discrete
vertical coordinate projected onto $\mathcal{Q}_h$. While it may seem counter-intuitive to take a discrete vertical gradient of the
vertical coordinate itself, rather than just representing this as unity, this weak representation is required in order to preserve the
exact balance of kinetic and potential energy exchanges, as shown below in Section 4.4.

		Finally, \eqref{eq:euler_wf_discrete_algebraic_no_par_no_perp_split_incidence} may be written in compact matrix notation as
		\begin{subequations}
			\label{eq:euler_wf_discrete_algebraic_no_par_no_perp_split_incidence_matrix}
			\begin{align}
				\massq\frac{\mathrm{d}\boldsymbol{\mathsf{\rho}}}{\mathrm{d}t} + \massq\boldsymbol{\mathsf{E}}^{3,2}_{\parallel}\boldsymbol{\mathsf{U}}^{\parallel} + \massq\boldsymbol{\mathsf{E}}^{3,2}_{\perp}\boldsymbol{\mathsf{U}}^{\perp} &= \boldsymbol{\mathsf{0}},\label{eq::dens_weak_discrete_algebraic_split_incidence_matrix}\\
				\massq\frac{\mathrm{d}\boldsymbol{\mathsf{\Theta}}}{\mathrm{d}t} + \massq\boldsymbol{\mathsf{E}}^{3,2}_{\parallel}\boldsymbol{\mathsf{F}}^{\parallel} + \massq\boldsymbol{\mathsf{E}}^{3,2}_{\perp}\boldsymbol{\mathsf{F}}^{\perp} &= \boldsymbol{\mathsf{0}},\label{eq::temp_weak_discrete_algebraic_split_incidence_matrix} \\
\luq\boldsymbol{\mathsf{\Theta}}
- \nuperp\boldsymbol{\mathsf{\theta}}
&= \boldsymbol{\mathsf{0}},\label{eq:rho_theta_Theta_discrete_algebraic_split_incidence_matrix}\\
c_p\Big(\frac{R}{p_0}\Big)^{R/c_v}\sum_{i=0}^{\ucardinal-1}\innerprod{(\qbasis_{i}\mathsf{\Theta}_{i})^{R/c_v}}{\qbasis_{j}}_{\Omega}
- \massq \boldsymbol{\mathsf{\Pi}} &= \boldsymbol{\mathsf{0}},
			\end{align}
		\end{subequations}
		with
\[
\mathsf{L}^{\mathcal{U}_{\perp},\mathcal{Q}}_{ij} := \innerprod{\rho_{h}\qbasis_{j}}{\ubasisperp_{i}}_{\Omega}.
\]

\subsection{Discrete energetics}

The conservation of energy for the rotating shallow water equations via balanced kinetic-potential
exchanges has previously been analysed \cite{LPG18} and experimentally verified \cite{LP18} for
the mixed mimetic spectral element method. In terms of energetics, the qualitative difference
between the rotating shallow water equations and the compressible Euler equations is the presence
of kinetic and internal energy exchanges. As such we here extend the previous analysis to demonstrate
that these exchanges may be balanced in the discrete from.

As seen before, we have that the discrete velocity field, $\boldsymbol{u}_{h}$, density, $\rho_{h}$, and density weighted potential temperature, $\Theta_{h}$, are
\begin{equation}
	\boldsymbol{u}_{h} = \sum_{i=0}^{\ucardinal}\mathsf{u}_{i}\ubasis_{i}, \qquad \rho_{h} = \sum_{i=0}^{\qcardinal}\mathsf{\rho}_{i}\qbasis_{i}, \quad\mathrm{and}\quad \Theta_{h} = \sum_{i=0}^{\qcardinal}\mathsf{\Theta}_{i}\qbasis_{i}\,.
\end{equation}
The discrete Hamiltonian $\mathcal{H}_{h} := \mathcal{H}[\boldsymbol{u}_{h}, \rho_{h}, \Theta_{h}]$ is then given by
\begin{equation}
	\mathcal{H}[\boldsymbol{u}_{h}, \rho_{h}, \Theta_{h}] = \int_{\Omega}\frac{1}{2}\rho_{h}\|\boldsymbol{u}_{h}\|^{2}\,\mathrm{d}\Omega + \int_{\Omega}\rho_{h}gz_{h}\,\mathrm{d}\Omega +  \int_{\Omega} c_{v} \left(\frac{R}{p_{0}}\right)^{\frac{R}{c_{v}}}\Theta_{h}^{\frac{c_{p}}{c_{v}}}\,\mathrm{d}\Omega\,.
\end{equation}
Using the definition of the variational derivative, see for example \cite{Celledoni12}, we can compute the variational derivative of the Hamiltonian with respect to the velocity, $\frac{\delta\mathcal{H}}{\delta\boldsymbol{u}_{h}}$
\begin{equation}
	\left.\frac{\mathrm{d}}{\mathrm{d}\epsilon}\mathcal{H}[\boldsymbol{u}_{h} + \epsilon\boldsymbol{v}_{h}, \rho_{h}, \Theta_{h}]\right|_{\epsilon=0} =:
\innerprod{\frac{\delta\mathcal{H}}{\delta\boldsymbol{u}_{h}}}{\boldsymbol{v}_{h}}, \qquad\forall\boldsymbol{v}_{h}\in\uspace\,.
\end{equation}
The left hand side of this expression may be evaluated to yield
\begin{equation}
\innerprod{\frac{\delta\mathcal{H}}{\delta\boldsymbol{u}_{h}}}{\boldsymbol{v}_{h}} =
\innerprod{\rho_{h}\boldsymbol{u}_{h}}{\boldsymbol{v}_{h}}, \qquad\forall\boldsymbol{v}_{h}\in\uspace\,.
\end{equation}
Since this expression is valid for all $\boldsymbol{v}_{h}\in\uspace$, then
\begin{equation}
\innerprod{\frac{\delta\mathcal{H}}{\delta\boldsymbol{u}_{h}}}{\ubasis_{j}} =
\innerprod{\rho_{h}\boldsymbol{u}_{h}}{\ubasis_{j}}
\stackrel{\eqref{eq::variations}}{=}
\innerprod{\boldsymbol{U}_{h}}{\ubasis_{j}}, \qquad j=0,\dots,\ucardinal\,.
\end{equation}
Following the same procedure we may obtain the weak equations for the variational derivative with respect to $\rho_{h}$
\begin{equation}
\innerprod{\frac{\delta\mathcal{H}}{\delta\rho_{h}}}{\qbasis_{j}} =
\innerprod{\frac{1}{2}\|\boldsymbol{u}_{h}\|^{2} + gz_h}{\qbasis_{j}}
\stackrel{\eqref{eq::variations}}{=}
\innerprod{\Phi_{h}}{\qbasis_{j}}, \qquad j=0,\dots,\qcardinal\,,
\end{equation}
and the variational derivative with respect to $\Theta_{h}$
\begin{equation}
\innerprod{\frac{\delta\mathcal{H}}{\delta\Theta_{h}}}{\qbasis_{j}} =
c_p\Bigg(\frac{R}{p_0}\Bigg)^{R/c_v}\innerprod{\Theta_h^{R/c_v}}{\qbasis_{j}}
\stackrel{\eqref{eq:Pi_Theta_discrete_algebraic}}{=}
\innerprod{\Pi_h}{\qbasis_{j}}
,\qquad j = 0,\dots,\qcardinal\,.
\end{equation}

The discrete compressible Euler equations \eqref{eq:euler_wf_discrete_prog}, \eqref{eq:euler_wf_discrete_diag} may then
be formulated as a skew-symmetric system for the discrete analogue of \eqref{eq:en_con_1}-\eqref{eq:da_dt_B_h} as
\begin{equation}\label{eq:euler_eqns_discrete_skew_symmetric}
\left[
\begin{array}{c}
	\massu\boldsymbol{\mathsf{u}}_{,t} \\
	\massq\mathsf{\rho}_{,t} \\
	\massq\mathsf{\Theta}_{,t}
\end{array}
\right] =
\begin{bmatrix}
-\boldsymbol{\mathsf{R}}_q & \left(\boldsymbol{\mathsf{E}}^{3,2}\right)^{\top}\massq &
\boldsymbol{\mathsf{S}}^{\mathcal{U}}\left(\massu\right)^{-1} \left(\boldsymbol{\mathsf{E}}^{3,2}\right)^{\top}\massq\\
-\massq\boldsymbol{\mathsf{E}}^{3,2} & \boldsymbol{\mathsf{0}} & \boldsymbol{\mathsf{0}} \\
-\massq\boldsymbol{\mathsf{E}}^{3,2}\left(\massu\right)^{-1} \boldsymbol{\mathsf{S}}^{\mathcal{U}} &
\boldsymbol{\mathsf{0}} & \boldsymbol{\mathsf{0}}\\
\end{bmatrix}
\left[
\begin{array}{c}
	\boldsymbol{\mathsf{U}} \\
	\mathsf{\Phi} \\
	\mathsf{\Pi}
\end{array}
\right] ,
\end{equation}
Multiplying both sides by $\begin{bmatrix}\boldsymbol{\mathsf{U}}^{\top} & \mathsf{\Phi}^{\top} & \mathsf{\Pi}^{\top}\end{bmatrix}$, gives
\begin{equation}
\boldsymbol{\mathsf{U}}^{\top}\massu\frac{\partial\boldsymbol{\mathsf{u}}}{\partial t} +
\boldsymbol{\mathsf{u}}^{\top}(\boldsymbol{\mathsf{T}}^{\mathcal{U}})^{\top}\frac{\partial\boldsymbol{\mathsf{\rho}}}{\partial t} +
g\boldsymbol{\mathsf{z}}^{\top}\massq\frac{\partial\boldsymbol{\mathsf{\rho}}}{\partial t} +
\boldsymbol{\mathsf{\Pi}}^{\top}\massq\frac{\partial\boldsymbol{\mathsf{\Theta}}}{\partial t} =
\frac{\partial K_h}{\partial t} + \frac{\partial P_h}{\partial t} + \frac{\partial I_h}{\partial t} = 0,
\end{equation}
where
\begin{equation}
K_h=\frac{1}{2}\boldsymbol{\mathsf{U}}^{\top}\massu\boldsymbol{\mathsf{u}} =
\boldsymbol{\mathsf{u}}^{\top}(\boldsymbol{\mathsf{T}}^{\mathcal{U}})^{\top}\boldsymbol{\mathsf{\rho}},\qquad
P_h=g\boldsymbol{\mathsf{z}}^{\top}\massq\boldsymbol{\mathsf{\rho}},\qquad
I_h=\frac{c_v}{c_p}\boldsymbol{\mathsf{\Pi}}^{\top}\massq\boldsymbol{\mathsf{\Theta}}
= \frac{c_v}{c_p}\Bigg(\frac{R}{p_0}\Bigg)^{R/c_v}\int_{\Omega}\Theta_h^{c_p/c_v}\mathrm{d}\Omega.
\end{equation}
Note that for
$\mathsf{R}_{q,ij} := \innerprod{\boldsymbol{q}_{h}\times\ubasis_{j}}{\ubasis_{i}}_{\Omega}$,
where $\boldsymbol{q}_h$ is the potential vorticity \cite{LPG18},
this is itself a skew-symmetric operator such that
$\boldsymbol{\mathsf{U}}^{\top}\boldsymbol{\mathsf{R}}_q\boldsymbol{\mathsf{U}} =
\boldsymbol{\mathsf{U}}^{\top}\boldsymbol{\mathsf{R}}\boldsymbol{\mathsf{u}} = \boldsymbol{\mathsf{0}}$.
As such neither $\boldsymbol{\mathsf{R}}_q$ nor $\boldsymbol{\mathsf{R}}$ projects onto the
energy in the discrete form.

Note also that the pressure gradient diagnostic equation \eqref{eq:Pi_P_discrete_algebraic} and the
temperature flux diagnostic equation \eqref{eq:theta_V_F_discrete_algebraic} appear within the
skew-symmetric operator in \eqref{eq:euler_eqns_discrete_skew_symmetric} within the top right and
bottom left blocks respectively. The discrete energy exchanges are therefore given as
\begin{align}
\frac{\partial K_h}{\partial t} &= g\boldsymbol{\mathsf{U}}^{\top}(\boldsymbol{\mathsf{E}}^{3,2})^{\top}\massq\boldsymbol{\mathsf{z}} +
\boldsymbol{\mathsf{U}}^{\top}\boldsymbol{\mathsf{S}}^{\mathcal{U}}\left(\massu\right)^{-1}
(\boldsymbol{\mathsf{E}}^{3,2})^{\top}\massq\boldsymbol{\mathsf{\Pi}},\label{eq::dKdt}\\
\frac{\partial P_h}{\partial t} &= -g\boldsymbol{\mathsf{z}}^{\top}\massq\boldsymbol{\mathsf{E}}^{3,2}\boldsymbol{\mathsf{U}},\label{eq::dPdt}\\
\frac{\partial I_h}{\partial t} &= -\boldsymbol{\mathsf{\Pi}}^{\top}\massq\boldsymbol{\mathsf{E}}^{3,2}
\left(\massu\right)^{-1}\boldsymbol{\mathsf{S}}^{\mathcal{U}}\boldsymbol{\mathsf{U}}.\label{eq::dIdt}
\end{align}

The right hand side terms of \eqref{eq::dKdt} exactly balance those of \eqref{eq::dPdt} and \eqref{eq::dIdt},
thus allowing for the exact balances of kinetic to potential and kinetic to internal energy respectively.
This holds for both the horizontal and vertical discretisations presented above, assuming periodic
boundary conditions in the horizontal, homogeneous Dirichlet boundary conditions for the vertical velocity
\eqref{eq:dirichlet_bc_velocity} and Neumann conditions for the Exner pressure \eqref{eq::bcs_neumann}
in the vertical. As an aside we note that energetic consistency is satisfied independent of the choice
of function space for $\theta_h$, since this only appears within $\boldsymbol{\mathsf{S}}^{\mathcal{U}}$,
thus justifying our choice to represent $\theta_h\in\mathcal{U}_{\perp}$.

\subsection{Metric terms}

The Jacobian matrix between local element coordinates $\bm \xi := (\xi, \eta, \zeta)$ and global coordinates
$\bm x := (\lambda, \phi, z)$ is given as
\begin{equation}
\mathsf{J} =
\begin{bmatrix}
\cos(\phi)\lambda,_{\xi} & \cos(\phi)\lambda,_{\eta} & 0    \\
\phi,_{\xi} & \phi,_{\eta} & 0    \\
0    & 0    & z,_{\zeta} \\
\end{bmatrix},
\end{equation}
where the subscripts represent derivatives with respect to local element coordinates and we have assumed
that all horizontal layers are perfectly flat, such that the projection of vertical
local coordinates onto horizontal global coordinates and horizontal local coordinates onto vertical
global coordinates are zero. The $H(\mathrm{curl},\Omega)$, $H(\mathrm{div},\Omega)$ and $L^2(\Omega)$
forms of the Piola transformation are given respectively as \cite{Natale16,Melvin19}
\begin{equation}
\mathsf{J}^{-\top},\qquad\frac{1}{J}\mathsf{J},\qquad\frac{1}{J},
\end{equation}
where $J$ is the determinant of the Jacobian matrix.
The metric transformations for the respective mass matrices are therefore
$(\mathsf{J}^{-\top})^{\top}\mathsf{J}^{-\top}$, $\frac{1}{J^2}\mathsf{J}^{\top}\mathsf{J}$
and $\frac{1}{J^2}$.
Since the horizontal and vertical components of both the $H(\mathrm{curl},\Omega)$ and
$H(\mathrm{div},\Omega)$ transformations are orthogonal,
these metric transformations are further simplified for these spaces as
\begin{subequations}
\begin{align}
\mathcal{W}_{h,\parallel}:&
\frac{1}{(\cos(\phi)\lambda,_{\xi}\phi,_{\eta} - \cos(\phi)\lambda,_{\eta}\phi,_{\xi})^2}
\begin{bmatrix}
 (\cos(\phi)\lambda,_{\eta})^2 + (\phi,_{\eta})^2 & -\cos^2(\phi)\lambda,_{\xi}\lambda,_{\eta} - \phi,_{\xi}\phi,_{\eta} \\
-\cos^2(\phi)\lambda,_{\xi}\lambda,_{\eta} - \phi,_{\xi}\phi,_{\eta} &  (\cos(\phi)\lambda,_{\xi})^2 + (\phi,_{\xi})^2   \\
\end{bmatrix},\\
\mathcal{W}_{h,\perp}:&
\frac{1}{(z,_{\zeta})^2},\\
\mathcal{U}_{h,\parallel}:&
\frac{1}{J^2}
\begin{bmatrix}
(\cos(\phi)\lambda,_{\xi})^2 + (\phi,_{\xi})^2  & \cos^2(\phi)\lambda,_{\xi}\lambda,_{\eta} + \phi,_{\xi}\phi,_{\eta}  \\
\cos^2(\phi)\lambda,_{\xi}\lambda,_{\eta} + \phi,_{\xi}\phi,_{\eta} & (\cos(\phi)\lambda,_{\eta})^2 + (\phi,_{\eta})^2 \\
\end{bmatrix},\\
\mathcal{U}_{h,\perp}:&
\frac{1}{(\cos(\phi)\lambda,_{\xi}\phi,_{\eta} - \cos(\phi)\lambda,_{\eta}\phi,_{\xi})^2},\\
\mathcal{Q}_h:&
\frac{1}{(z,_{\zeta})^2(\cos(\phi)\lambda,_{\xi}\phi,_{\eta} - \cos(\phi)\lambda,_{\eta}\phi,_{\xi})^2}.
\end{align}
\end{subequations}
As a special case, we also note that an inner product between basis functions in
$\mathcal{U}_{h,\perp}$ and functions in $\mathcal{Q}_{h}$ has the metric term
\begin{equation}
\frac{z,_{\zeta}}{J}\cdot\frac{1}{J} = \frac{1}{z,_{\zeta}(\cos(\phi)\lambda,_{\xi}\phi,_{\eta} - \cos(\phi)\lambda,_{\eta}\phi,_{\xi})^2}.
\end{equation}

As stated these metric terms do not account for any projection of horizontal vector components
onto vertical components and vice versa. As such we limit ourselves here to the case where the
horizontal and vertical degrees of freedom are strictly orthogonal, and we are unable to
account for the tilting of layers or bottom topography in the current implementation. We note
however that the tilting of horizontal layers in order to represent topography is
naturally incorporated into the finite element formulation via these cross terms, and these may
be implemented at a later date with minimal disruption to the current formulation.

Note that the vorticity term $\langle\bm\sigma_h,\bm\omega_{h,\parallel,\parallel}\times u_{h,\perp}\rangle$
is alternatively formulated as $\langle\bm\sigma_h,u_{h,\perp}\nabla_{\perp}\bm{u}_{h,\parallel}\rangle$ (with
the vertical derivative derived in the weak form). In this form the vertical velocity
derivatives may be interpreted as being oriented \emph{normal} to the edges, rather than
\emph{tangent} to them, and as such the $H(\mathrm{div},\Omega)$ form of the Piola transformation
is used to construct the metric term for this operator as well as for \eqref{eq:curl_u_omega_discrete_algebraic_par_par_split_incidence}.

\section{Time stepping}

In this section we introduce the horizontally-explicit/vertically implicit (HEVI) time stepping
scheme employed in the model. Such schemes are popular in non-hydrostatic atmospheric modelling since they
negate the explicit time step restrictions associated with the fast time scales of the vertical gravity
and acoustic waves \cite{Ullrich12,Giraldo13,Lock14,Bao15,Gardner18}. These schemes are often constructed
by perturbing the vertical dynamics around a leading order state of hydrostatic balance
\cite{Giraldo13,Bao15} and/or further splitting the vertical dynamics so as to only solve for the terms
responsible for the fast dynamics implicitly \cite{Gardner18}. Here we avoid these formulations due to
our concern that they may break the skew-symmetric structure of the discrete system \eqref{eq:euler_eqns_discrete_skew_symmetric}
that is central to the conservation of energy in the spatial discretisation.
Nevertheless we note that while our current implicit scheme detailed below preserves the exact balance
kinetic to potential energy exchanges, balance is not satisfied for the  vertical
kinetic to internal energy exchanges due to the formulation of the pressure gradient term.
Restoring this balance within the implicit vertical scheme is a subject of ongoing research.
We further emphasise that both our time splitting scheme as well as our
individual implicit vertical and explicit horizontal time integration schemes will conserve energy only
to truncation order in time.

\subsection{Directional splitting}

We use a Strang carryover splitting scheme to partition the horizontal and vertical dynamics \cite{Ullrich12,Lock14}. 
Consider the previously introduced splitting of the spatial operator (here denoted by $\boldsymbol{L}$) into a vertical ($\boldsymbol{L}^{v}$) and horizontal ($\boldsymbol{L}^{h}$)  components: equations \eqref{eq:euler_wf_discrete_algebraic_par_split_incidence_matrix} (horizontal), and equations \eqref{eq:euler_wf_discrete_algebraic_perp_split_incidence_matrix} together with equations \eqref{eq:euler_wf_discrete_algebraic_no_par_no_perp_split_incidence_matrix} (vertical). With this splitting we can write the full system of equations as
\begin{equation}
	\frac{\mathrm{d}\boldsymbol{b}}{\mathrm{d}t} = \boldsymbol{L}(\boldsymbol{b}) = \boldsymbol{L}^{h}(\boldsymbol{b}) + \boldsymbol{L}^{v}(\boldsymbol{b})\,,
\end{equation}
where $\boldsymbol{b} := [\boldsymbol{u}, \boldsymbol{U}, \boldsymbol{F}, \boldsymbol{P}, \boldsymbol{\omega}, \rho, \theta, \Theta, \Pi]$. Then, the time stepping procedure to evolve from time step $t^{n}$ to $t^{n+1}$ follows the sequence
\begin{align}
	\frac{\mathrm{d}\boldsymbol{b}'}{\mathrm{d}t} &= \boldsymbol{L}^{v}(\boldsymbol{b}'), \qquad t\in[t^{n}, t^{n}+\frac{\Delta t}{2}] \qquad \text{(vertical half step)}\nonumber\\
	\boldsymbol{b}'(t^{n}) &= \boldsymbol{b}(t^{n}) \label{eq:vertical_half_step_first}\\
	\phantom{x} \nonumber\\
	\frac{\mathrm{d}\boldsymbol{b}''}{\mathrm{d}t} &= \boldsymbol{L}^{h}(\boldsymbol{b}''), \qquad t\in[t^{n}, t^{n}+\Delta t] \qquad \text{(horizontal full step)}\nonumber\\
	\boldsymbol{b}''(t^{n}) &= \boldsymbol{b}'(t^{n} + \frac{\Delta t}{2}) \label{eq:horizontal_full_step} \\
	\phantom{x} \nonumber \\
			\frac{\mathrm{d}\boldsymbol{b}'''}{\mathrm{d}t} &= \boldsymbol{L}^{v}(\boldsymbol{b}'''), \qquad t\in[t^{n}+\frac{\Delta t}{2}, t^{n+1}] \qquad \text{(vertical half step)}\nonumber\\
	\boldsymbol{b}'''(t^{n} + \frac{\Delta t}{2}) &= \boldsymbol{b}''(t^{n+1}) \label{eq:vertical_half_step_second}\\
	\phantom{x} \nonumber \\
	\boldsymbol{b}^{n+1} &= \boldsymbol{b}'''(t^{n + 1}) \qquad \text{(update)}\,. \label{eq:update}
\end{align}

This splitted time evolution procedure is then numerically integrated. The first vertical half time step \eqref{eq:vertical_half_step_first} is explicitly integrated using a forward Euler scheme (\emph{carried over}) as
	\begin{equation}
		\bm b' = \bm b^n + \frac{\Delta t}{2}L^{v}(\bm b^n)\,.
	\end{equation}
	The horizontal full time step \eqref{eq:horizontal_full_step} is numerically integrated using an explicit 
\blue{third order } stiffly stable Runge-Kutta scheme, see \secref{sec:explicit_horizontal_solve} for more details. Finally, the second vertical half step \eqref{eq:vertical_half_step_second} is integrated using an implicit 
Euler scheme. A nonlinear Picard iteration, $m$,
is applied to solve this vertical system, which is assumed to converge once
$|\boldsymbol{b}^{m+1}(t^{n+1})-\boldsymbol{b}^m(t^{n+1})|_2/|\boldsymbol{b}^{m+1}(t^{n+1})|_2 < 10.0^{-8}$, where $|\cdot|_2$
is the $L^2$ norm. For a detailed discussion see \secref{sec:vertical_implicit_half_step}. For the 30 models levels of the baroclinic test case described below in Section 6.1, 13 Picard iterations are required to reach convergence. 

\blue{This scheme is formally second order accurate in time due to the Strang splitting, with the 
horiztonal scheme being third order accurate, and the vertical scheme being equivalent to a second 
order trapezoidal scheme \cite{Lock14}. } Our anecdotal experience is that if a
first order horizontal-vertical splitting is used, as opposed to the second order splitting described
here, then the horizontal-vertical coupling is too weak to allow for the correct transfer of potential
and internal energy to vertical kinetic energy required to properly simulate the baroclinic instability
described below in Section 6.1.

Since all the solution variables in $\bm b$ required to solve for the vertical dynamics are
discontinuous across horizontal element boundaries, the inner linear system above is solved in serial for
each horizontal element, via a direct LU solve using the PETSc library \cite{petsc-web-page,petsc-user-ref,petsc-efficient}.
Similarly, since all the solution variables involved in the vertical dynamics are discretized
on function spaces that are discontinuous across vertical element boundaries, this system may be solved
independently, in parallel for each horizontal layer. The mass matrices are solved using the PETSc GMRES
solver, with a block Jacobi preconditioner for each element, as was done for the shallow water equations
in our previous work \cite{LP18}.

\subsection{Explicit horizontal solve} \label{sec:explicit_horizontal_solve}

The horizontal dynamics are solved using an explicit stiffly stable Runge-Kutta scheme of the form \cite{Shu88,Ullrich12}
\begin{align}
\bm b^{(1)} &= \bm b'(t^{n} + \frac{\Delta t}{2}) + \Delta tL^{h}(\bm b'(t^{n} + \frac{\Delta t}{2})),\label{eq::horiz_1}\\
\bm b^{(2)} &= \frac{3}{4}\bm b'(t^{n} + \frac{\Delta t}{2}) + \frac{1}{4}\bm b^{(1)} + \frac{\Delta t}{4}L^{h}(\bm b^{(1)}),\label{eq::horiz_2}\\
\bm {b}_h''(t^{n} + \Delta t) &= \frac{1}{3}\bm b'(t^{n} + \frac{\Delta t}{2}) + \frac{2}{3}\bm b^{(2)} + \frac{2\Delta t}{3}L^{h}(\bm b^{(2)})\,,\label{eq::horiz_3}
\end{align}
where $b'(t)$ and $b''(t)$ are as defined in \eqref{eq:vertical_half_step_first} and \eqref{eq:horizontal_full_step}, respectively.

\subsection{Implicit vertical solve} \label{sec:vertical_implicit_half_step}

The implicit vertical solve involves a half step, $\textstyle\frac{1}{2}\Delta t$, from the state
at the end of the horizontal solve, $b''(t^{n} + \Delta t)$, to the end of the time
level $b^{n+1}$.
Following \cite{Gassmann13} we begin by taking the logarithm of the equation of state \eqref{eq::eos} as
\begin{equation}
\mathrm{ln}(\Pi) = \frac{R}{c_v}\Bigg(\mathrm{ln}(\rho\theta) + \mathrm{ln}\Bigg(\frac{R}{p_0}\Bigg)\Bigg) + \mathrm{ln}(c_p).
\end{equation}
Differentiating both sides with respect to time over a half step via the chain rule,
while recalling that the second and third terms on the right hand side are simply constants then gives
\begin{equation}
\frac{\Pi^{n+1} - \Pi''}{\textstyle\frac{1}{2}\Delta t\Pi''} = \frac{R}{c_v}\frac{(\rho\theta)^{n+1} - (\rho\theta)''}{\textstyle\frac{1}{2}\Delta t(\rho\theta)''}.
\end{equation}
Substituting in the potential temperature equation \eqref{eq::temp}, we
then express the evolution of the Exner pressure as
\begin{equation}\label{eq::exner}
(\rho\theta)''\Pi^{n+1} =
(\rho\theta)''\Pi'' - \frac{\Delta tR}{2c_v}\Pi''\nabla\cdot(\rho\bm u\theta)^{n+1}.
\end{equation}
Note the similarity between \eqref{eq::exner} and the internal energy evolution equation
\eqref{eq:time_variation_internal_energy}. The vertical dynamics may then be discretized in time
(with the vorticity components omitted, as discussed in Section 4.3) as

\begin{align}
u_{\perp}^{n+1} + \frac{\Delta t}{4}\nabla_{\perp}(u_{\perp}^{n+1})^2 +
\frac{\Delta t}{2}\theta\nabla_{\perp}\Pi^{n+1} &= u_{\perp}'' - \frac{\Delta t}{2}g\nabla_{\perp}z\label{eq::vert_mom_disc}\\
\rho^{n+1} + \frac{\Delta t}{2}\nabla_{\perp}\cdot(\rho u_{\perp})^{n+1} &= \rho''\\
\Theta^{n+1} + \frac{\Delta t}{2}\nabla_{\perp}\cdot(u_{\perp}\Theta)^{n+1} &= \Theta''\\
\Theta^n\Pi^{n+1} &= \Theta''\Pi'' - \frac{\Delta tR}{2c_v}\Pi''\nabla_{\perp}\cdot(u_{\perp}\Theta)^{n+1}
\end{align}
We define the additional matrix operators
$\boldsymbol{\mathsf{M}}^{\mathcal{Q}}_{\Theta}$,
$\boldsymbol{\mathsf{M}}^{\mathcal{Q}}_{\Pi}$ and
$\boldsymbol{\mathsf{M}}^{\mathcal{U}_{\perp}}_{\Theta}$,
for which
\[
\mathsf{M}^{\mathcal{Q}}_{\Theta_h,ij} := \innerprod{\Theta_h\qbasis_{j}}{\qbasis_{i}}_{\Omega},\qquad
\mathsf{M}^{\mathcal{Q}}_{\Pi_h,ij} := \innerprod{\Pi_h\qbasis_{j}}{\qbasis_{i}}_{\Omega},\qquad
\mathsf{M}^{\mathcal{U}_{\perp}}_{\Theta_h,ij} := \innerprod{\Theta_h\epsilon^{\mathcal{U}_{\perp}}_{j}}{\epsilon^{\mathcal{U}^{\perp}}_{i}}_{\Omega},
\]
Dropping the $n+1$ superscripts for variables at the end of the vertical half step (but keeping those for variables at time level $''$),
the discrete form of \eqref{eq::exner} is then given as:

\begin{equation}\label{eq::vert_exner}
\massq_{\Theta''}\boldsymbol{\mathsf{\Pi}} = \massq_{\Theta''}\boldsymbol{\mathsf{\Pi}}'' -
\frac{\Delta tR}{2c_v}\massq_{\Pi''}\boldsymbol{\mathsf{E}}^{3,2}_{\perp}
(\massuperp)^{-1}\massuperp_{\Theta}\boldsymbol{\mathsf{u}}_{\perp},
\end{equation}
and the discrete form of \eqref{eq::vert_mom_disc} as:
\begin{equation}\label{w_vert_disc}
\massuperp\boldsymbol{\mathsf{u}}_{\perp} + \frac{\Delta t}{4}(\boldsymbol{\mathsf{E}}^{3,2}_{\perp})^{\top}
\boldsymbol{\mathsf{T}}^{\mathcal{U}_{\perp}}\boldsymbol{\mathsf{u}}_{\perp} +
\frac{\Delta t}{2}\boldsymbol{\mathsf{S}}^{\mathcal{U}_{\perp}}(\massuperp)^{-1}(\boldsymbol{\mathsf{E}}^{3,2}_{\perp})^{\top}
\massq\boldsymbol{\mathsf{\Pi}} =
\massuperp\boldsymbol{\mathsf{u}}_{\perp}'' -
\frac{\Delta t g}{2}(\boldsymbol{\mathsf{E}}^{3,2}_{\perp})^{\top}\massq\boldsymbol{\mathsf{z}}.
\end{equation}
Substituting \eqref{eq::vert_exner} into \eqref{w_vert_disc} gives
\begin{multline}
\massuperp\boldsymbol{\mathsf{u}}_{\perp} + \frac{\Delta t}{4}(\boldsymbol{\mathsf{E}}^{3,2}_{\perp})^{\top}
\boldsymbol{\mathsf{T}}^{\mathcal{U}_{\perp}}\boldsymbol{\mathsf{u}}_{\perp} +
\frac{\Delta t}{2}\boldsymbol{\mathsf{S}}^{\mathcal{U}_{\perp}}(\massuperp)^{-1}(\boldsymbol{\mathsf{E}}^{3,2}_{\perp})^{\top}\massq \\
\Bigg(
\boldsymbol{\mathsf{\Pi}}'' - \frac{\Delta tR}{2c_v}(\massq_{\Theta''})^{-1}\massq_{\Pi''}\boldsymbol{\mathsf{E}}^{3,2}_{\perp}
(\massuperp)^{-1}\massuperp_{\Theta}\boldsymbol{\mathsf{u}}_{\perp}
\Bigg) =
\massuperp\boldsymbol{\mathsf{u}}_{\perp}'' - \frac{\Delta t g}{2}(\boldsymbol{\mathsf{E}_{\perp}}^{3,2})^{\top}\massq\boldsymbol{\mathsf{z}}.
\end{multline}
Finally, rearranging gives an expression for the vertical velocity at each fixed point Picard iteration as:
\begin{multline}
\Bigg[\massuperp + \frac{\Delta t}{4}(\boldsymbol{\mathsf{E}}^{3,2}_{\perp})^{\top}\boldsymbol{\mathsf{T}}^{\mathcal{U}_{\perp}} -
\frac{\Delta t^2R}{4c_v}\boldsymbol{\mathsf{S}}^{\mathcal{U}_{\perp}}(\massuperp)^{-1}(\boldsymbol{\mathsf{E}}^{3,2}_{\perp})^{\top}
\massq
(\massq_{\Theta''})^{-1}\massq_{\Pi''}\boldsymbol{\mathsf{E}}^{3,2}_{\perp}(\massuperp)^{-1}\massuperp_{\Theta}\Bigg]\boldsymbol{\mathsf{u}}_{\perp}
= \\
\massuperp\boldsymbol{\mathsf{u}}_{\perp}'' -
\frac{\Delta t}{2}\Bigg(g(\boldsymbol{\mathsf{E}}^{3,2}_{\perp})^{\top}\massq\boldsymbol{\mathsf{z}} +
\boldsymbol{\mathsf{S}}^{\mathcal{U}_{\perp}}(\massuperp)^{-1}(\boldsymbol{\mathsf{E}}^{3,2}_{\perp})^{\top}\massq\boldsymbol{\mathsf{\Pi}}''\Bigg).\label{eq::vert_w}
\end{multline}
In order to incorporate the pressure gradient term into the left hand
side in \eqref{eq::vert_w} so as to ensure the stable implicit solution of vertical motions, we have
sacrificed the energetically consistent formulation of the vertical pressure gradient term as presented in
\eqref{eq::mom_weak_discrete_algebraic_perp_split_incidence}, \eqref{eq:Pi_P_discrete_algebraic_perp_split_incidence}
and \eqref{eq:euler_eqns_discrete_skew_symmetric}.
As such we do not expect the vertical kinetic and internal exchanges to exactly
balance, as we do for the horizontal explicit discretization.
However we note that unlike the perturbation and operator splitting approaches to HEVI discretisations
discussed above, our current approach leaves open the possibility of recovering the exact balance of kinetic
and internal energy exchanges. We are actively exploring preconditioning strategies to address this problem.
We also note that our implicit vertical solve may potentially dampen the fast vertical motions for
which the time scales are not resolved, in comparison to sub-cycled vertically explicit formulations.

The Exner pressure at the current Picard iteration is then derived from \eqref{eq::vert_exner}, and the
corresponding Picard iteration solves for the other variables are then given (omitting the horizontal
terms) as
\begin{align}
\nuperp\boldsymbol{\mathsf{\theta}} &= \luq\boldsymbol{\mathsf{\Theta}}
\label{eq::vert_theta}\\
\boldsymbol{\mathsf{U}}_{\perp} &= (\massuperp)^{-1}\boldsymbol{\mathsf{N}}^{\mathcal{U}_{\perp}}\boldsymbol{\mathsf{u}}_{\perp}\label{eq::vert_W}\\
\boldsymbol{\mathsf{\rho}} &= \boldsymbol{\mathsf{\rho}}'' - \frac{\Delta t}{2}\boldsymbol{\mathsf{E}}^{3,2}_{\perp}\boldsymbol{\mathsf{U}}_{\perp}\label{eq::vert_rho}\\
\boldsymbol{\mathsf{\Theta}} &= \boldsymbol{\Theta}'' - \frac{\Delta t}{2}\boldsymbol{\mathsf{E}}^{3,2}_{\perp}(\massuperp)^{-1}
\boldsymbol{\mathsf{S}}^{\mathcal{U}_{\perp}}\boldsymbol{\mathsf{U}}_{\perp}\label{eq::vert_Theta}.
\end{align}

\subsection{Dissipative terms}

In order to stabilise the model we also include various dissipative terms. These include a biharmonic
viscosity on both the horizontal momentum and temperature equations \cite{Guba14,LP18} with a value of
$0.072\Delta x^{3.2}$, where $\Delta x$ is the average spacing between GLL nodes.
A Rayleigh friction term with a coefficient
of $0.2$ is also applied to the top layer of the vertical momentum equation, as is often used in
atmospheric models to suppress orographically forced gravity waves \cite{Klemp08}.
This term is added to \eqref{eq::mom_weak_discrete_algebraic_perp_split_incidence} as
$\sum_{i=0}^{\ucardinalperp-1}0.2\innerprod{\delta_{i,k}\ubasisperp_{i}}{\ubasisperp_{j}\delta_{j,l}}_{\Omega}\mathsf{u}_{i,\perp}$,
where $k,l$ are degrees of freedom for the trial and test functions in the top layer only,
and $\delta_{i,k}$ is the standard delta function.
While this Rayleigh
friction term is not strictly necessary for the stability of the simulation presented here, it greatly
reduces the noise in the energetic profiles as the model adjusts to a state of hydrostatic balance from
its initial conditions. No viscosity is required to stabilize the vertical solution, perhaps because
this is only second order accurate, so the internal dissipation of this low order discretization is
sufficient to prevent nonlinear instabilities.

We emphasise that these dissipative terms will necessarily remove energy from the system,
and in the case of the horizontal viscous terms are necessary to stabilise the model by arresting
the nonlinear cascades at the grid scale. While upwinding terms may be used instead of viscosity as
a means of suppressing grid scale oscillations without dissipating energy, by systematically adding
these to the skew-symmetric formulation \cite{Wimmer19}, here we limit ourselves to a consideration
only of the internal energetic processes, and not the external forcings.

\section{Results}

\subsection{Baroclinic instability}

We validate the model using a dry baroclinic instability test case \cite{Ullrich14} with the shallow
atmosphere approximation. The appeal of this test case is that the initial condition is specified
for a $z-$level vertical coordinate, whereas other such test cases that are defined on pressure level
vertical coordinates require the solution of a nonlinear problem in order to compute the corresponding
$z-$level configuration. The initial state is one of geostrophic horizontal and hydrostatic
vertical balance, overlaid with a small, $\mathcal{O}(1\mathrm{m/s})$, perturbation to the zonal and
meridional velocity components.

The model was run with $24\times 24$ elements of degree $p=3$ on each face of the cubed sphere
(and linear elements in the vertical),
for an averaged resolution of $\Delta x\approx 128\mathrm{km}$ and 30 vertical levels on 96
processors ($4^2=16$ per face of the cubed sphere) with a time step of $\Delta t = 120\mathrm{s}$.
While the vertical dynamics are all solved implicitly and so do not limit the time step size, the
explicit horizontal dynamics present both diffusive and advective CFL restrictions due to the
biharmonic viscosity and sound waves respectively.

\begin{figure}[!hbtp]
\begin{center}
\includegraphics[width=0.48\textwidth,height=0.36\textwidth]{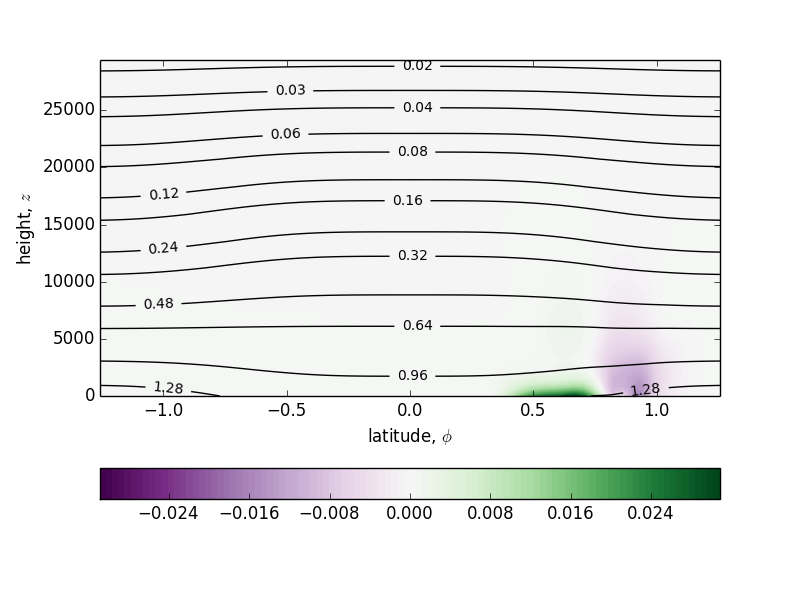}
\includegraphics[width=0.48\textwidth,height=0.36\textwidth]{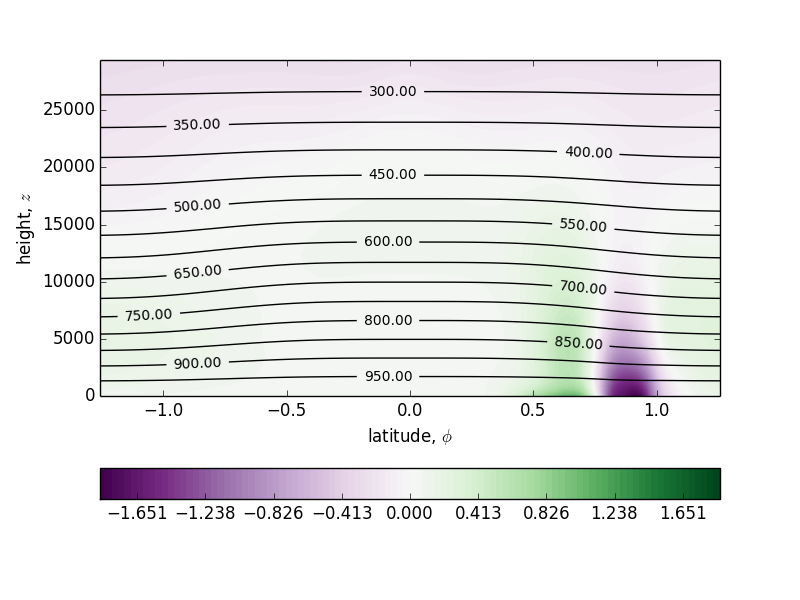}
\caption{Zonal averages of density, $\rho_h$ in $\mathrm{kg}\cdot\mathrm{m}^{-3}$ (left) and
Exner pressure, $\Pi_h$ in $\mathrm{m}^{2}\mathrm{s}^{-2}\mathrm{K}^{-1}$ (right) at day 10.
Contours represent absolute values, and shades represent differences from initial values.}
\label{fig::zonal_avg_1}
\end{center}
\end{figure}

\begin{figure}[!hbtp]
\begin{center}
\includegraphics[width=0.48\textwidth,height=0.36\textwidth]{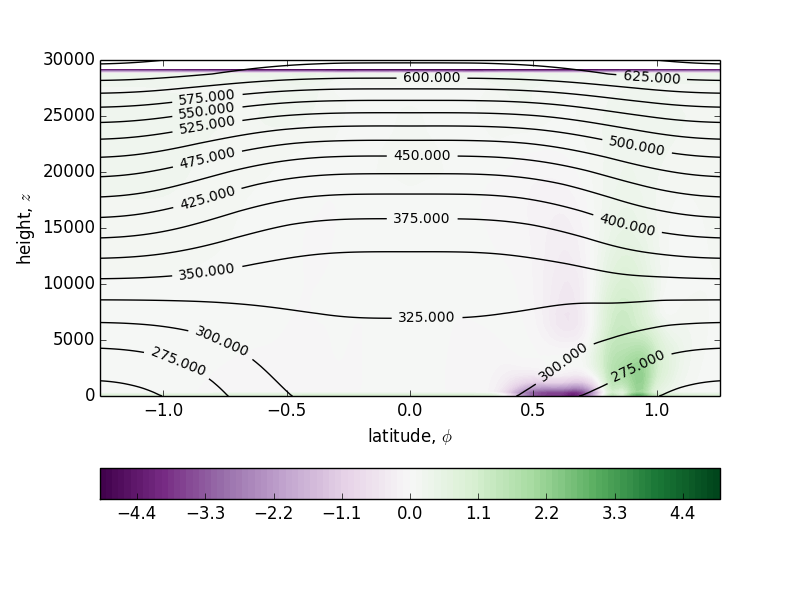}
\includegraphics[width=0.48\textwidth,height=0.36\textwidth]{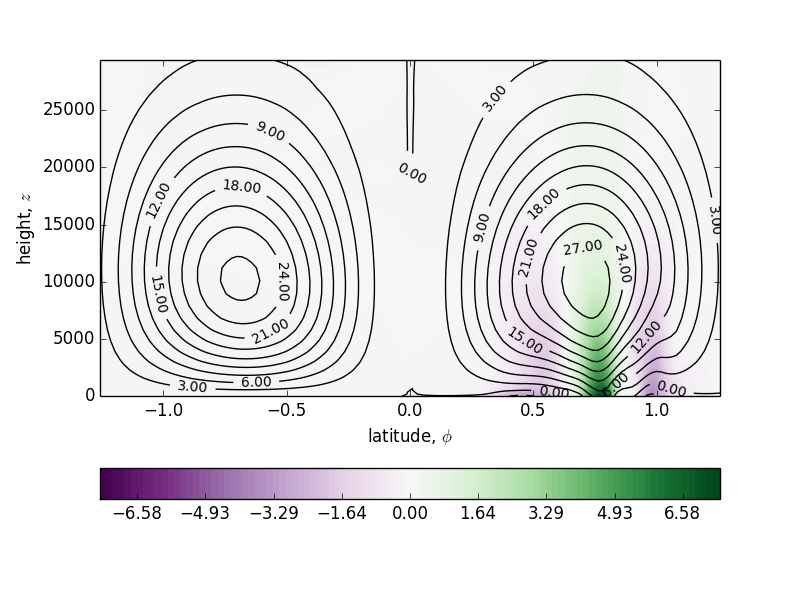}
\caption{Zonal averages of potential temperature, $\theta_h$ in $\mathrm{K}$ (left) and
zonal velocity, $u_h$ in $\mathrm{m}\cdot\mathrm{s}^{-1}$ (right) at day 10.
Contours represent absolute values, and shades represent differences from initial values.}
\label{fig::zonal_avg_2}
\end{center}
\end{figure}

Figures \ref{fig::zonal_avg_1} and \ref{fig::zonal_avg_2} show the zonal averages of density $\rho$,
Exner pressure $\Pi$, potential temperature $\theta$ and zonal velocity $u$ at day 10
(solid lines), as well as the differences between the final and initial states. These profiles
show little difference between the initial and final states, with the exception of the zonal velocity,
which exhibits a small kink near the bottom boundary where the baroclinic instability occurs,
demonstrating that the leading order geostrophic and hydrostatic balances in the horizontal and vertical
are well satisfied. The potential temperature (Fig. \ref{fig::zonal_avg_2}), is approximately
$15^{\circ}\mathrm{K}$ cooler in the top layer at day 10, which is outside the range of the
color bar. This is due to the hydrostatic adjustment of the mean state, as discussed below.

\begin{figure}[!hbtp]
\begin{center}
\includegraphics[width=0.48\textwidth,height=0.36\textwidth]{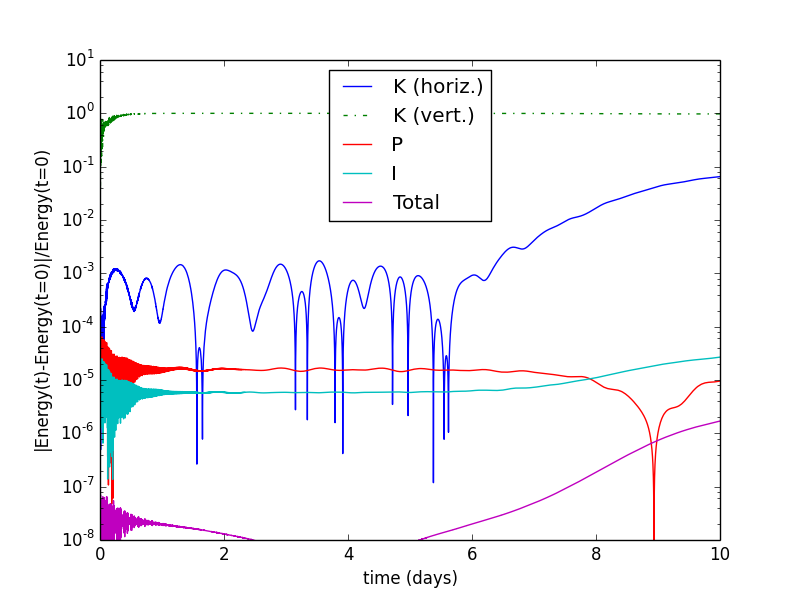}
\includegraphics[width=0.48\textwidth,height=0.36\textwidth]{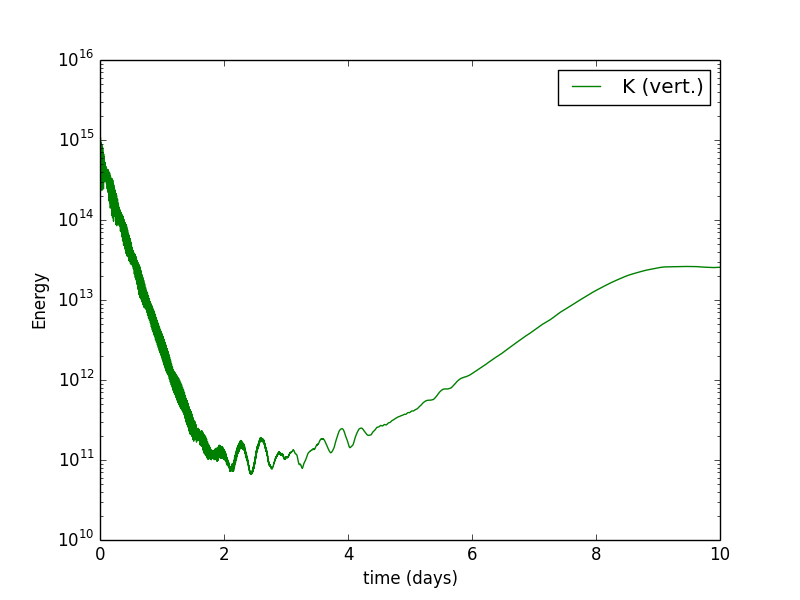}
\caption{Left: normalised difference in energy with respect to initial values. Right: vertical kinetic energy.}
\label{fig::energetics_1}
\end{center}
\end{figure}

\begin{figure}[!hbtp]
\begin{center}
\includegraphics[width=0.48\textwidth,height=0.36\textwidth]{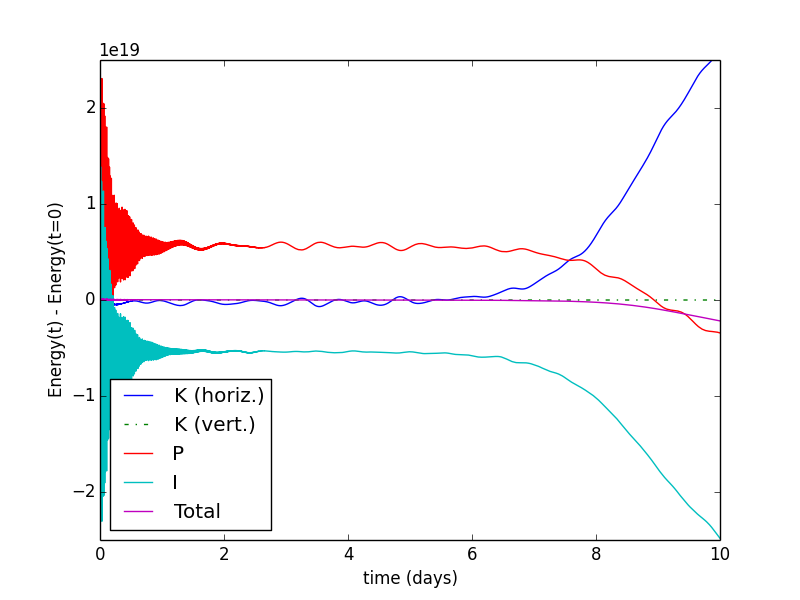}
\includegraphics[width=0.48\textwidth,height=0.36\textwidth]{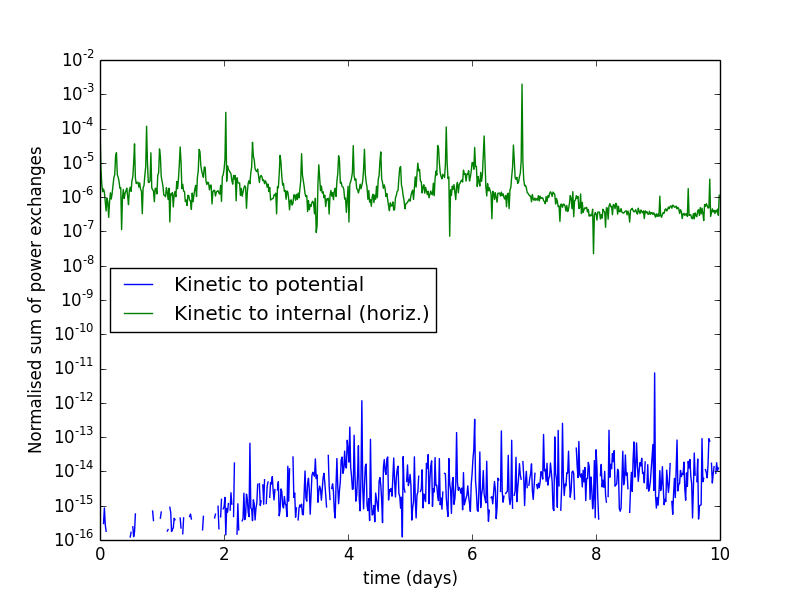}
\caption{Left: difference in energy with respect to initial values.
Right: normalised difference in ($K$)inetic to ($P$)otential energy,
$|\int K\mathrm{\ to\ }P\mathrm{d}\Omega + \int P\mathrm{\ to\ }K\mathrm{d}\Omega|/|\int K\mathrm{\ to\ }P\mathrm{d}\Omega|$,
and horizontal ($K$)inetic to ($I$)nternal energy,
$|\int K\mathrm{\ to\ }I\mathrm{d}\Omega + \int I\mathrm{\ to\ }K\mathrm{d}\Omega|/|\int K\mathrm{\ to\ }I\mathrm{d}\Omega|$.
}
\label{fig::energetics_2}
\end{center}
\end{figure}

Figures \ref{fig::energetics_1} and \ref{fig::energetics_2} show the evolution of the kinetic
(horizontal and vertical), potential and internal energy with time, and the associated exchanges.
These are shown on both logarithmic scales for their normalised absolute difference from initial
value (Fig. \ref{fig::energetics_1}), and as a direct difference between their current and
initial value (Fig. \ref{fig::energetics_2}).

The growth in the baroclinic instability is evident in the increase in kinetic energy, and the
reduction in both potential and internal energy as isopycnals flatten in the region of the
instability. Note that the total amounts of potential and internal energy are approximately
$3.6\times 10^{23}$ and $9.2\times 10^{23}\mathrm{kg\cdot m^2s^{-2}}$ respectively, and so are
several orders of magnitude greater than the amounts of horizontal and vertical kinetic energy
(approximately $4.0\times 10^{20}$ and $2.5\times 10^{13}\mathrm{kg\cdot m^2s^{-2}}$ respectively).
As such the flattening of the density contours from which the baroclinic instability draws energy
are barely evident in Fig. \ref{fig::zonal_avg_1}.

Figure \ref{fig::energetics_2} also shows the normalised sum of the globally integrated
kinetic to potential \eqref{eq::dKdt} and potential to kinetic \eqref{eq::dPdt} energy exchanges,
as well as the (horizontal) globally integrated kinetic to internal \eqref{eq::dKdt} and internal
to kinetic \eqref{eq::dIdt} exchanges. The kinetic to potential and potential to kinetic exchanges
balance to machine precision, owing to the exact skew-symmetry of these operators, as shown in
\eqref{eq:euler_eqns_discrete_skew_symmetric}. The kinetic to internal exchanges balance only
to approximately $\mathcal{O}(10^{-6})$, and are sometimes out by a factor of $\mathcal{O}(10^{-3})$.
This is due to the fact that both terms involve the inverse of a $H(\mathrm{div},\Omega)$ mass
matrix, the action of which is approximated via an iterative Krylov method. We therefore anticipate
pointwise errors in these exchanges, which may be compounded by the fact that these inverses are
approximated twice, once to determine the Exner pressure gradient \eqref{eq:Pi_P_discrete_algebraic}, and a
second time to determine the density weighted potential temperature flux \eqref{eq:theta_V_F_discrete_algebraic}.
Nevertheless, these errors exhibit no systematic drift throughout the simulation.
We do not show the vertical kinetic to internal balance errors, as our vertical implicit solve
does not satisfy this balance.

The observation that potential energy is greater than its initial value for most of the
simulation, while internal energy is smaller, as shown in Figs. \ref{fig::energetics_1} and
\ref{fig::energetics_2}, is explained by the fact that the initial condition is not precisely
in hydrostatic balance. As such the initial adjustment process leads to a slight rise in the
fluid, resulting in an increase in potential energy, and a corresponding reduction in pressure,
leading to a reduction of internal energy via the equation of state. These changes are of $\mathcal{O}(10^{-4})$
compared to the total amounts of potential and internal energy in the system. The high frequency
oscillation in the internal and potential energies observed over the first 24 hours of the simulation
during this adjustment process is reduced by the application of the Rayleigh friction to the top
layer of the vertical momentum equation.

\begin{figure}[!hbtp]
\begin{center}
\includegraphics[width=0.48\textwidth,height=0.36\textwidth]{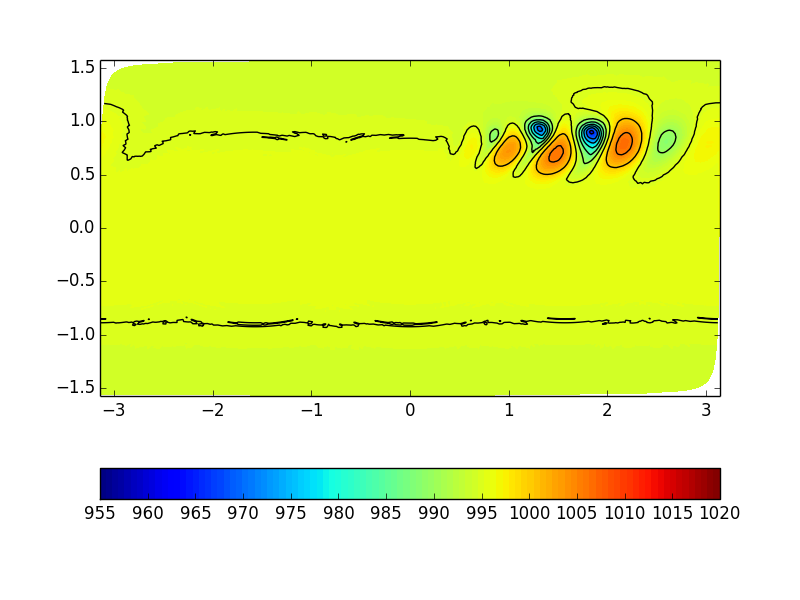}
\includegraphics[width=0.48\textwidth,height=0.36\textwidth]{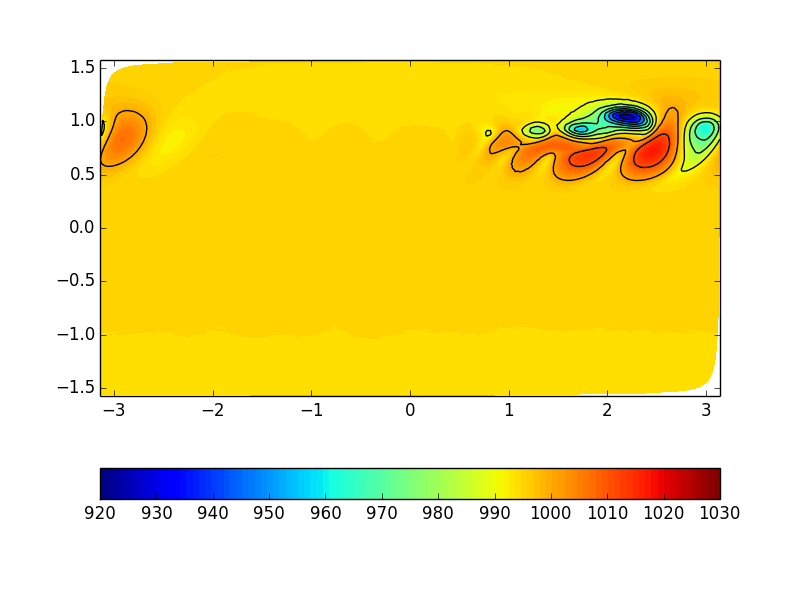}
\caption{Bottom level pressure, $p_h$ (in $\mathrm{hPa}$) day 8 (left) and 10 (right).}
\label{fig::pressure}
\end{center}
\end{figure}

\begin{figure}[!hbtp]
\begin{center}
\includegraphics[width=0.48\textwidth,height=0.36\textwidth]{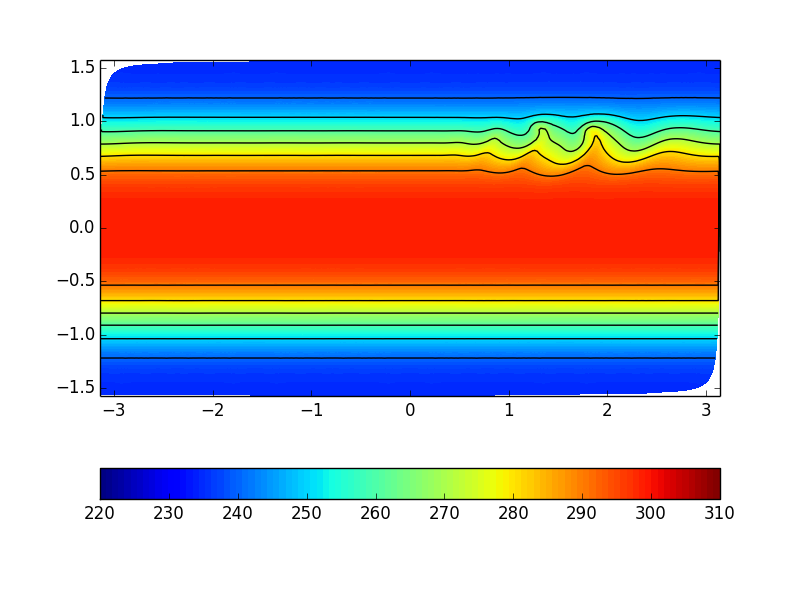}
\includegraphics[width=0.48\textwidth,height=0.36\textwidth]{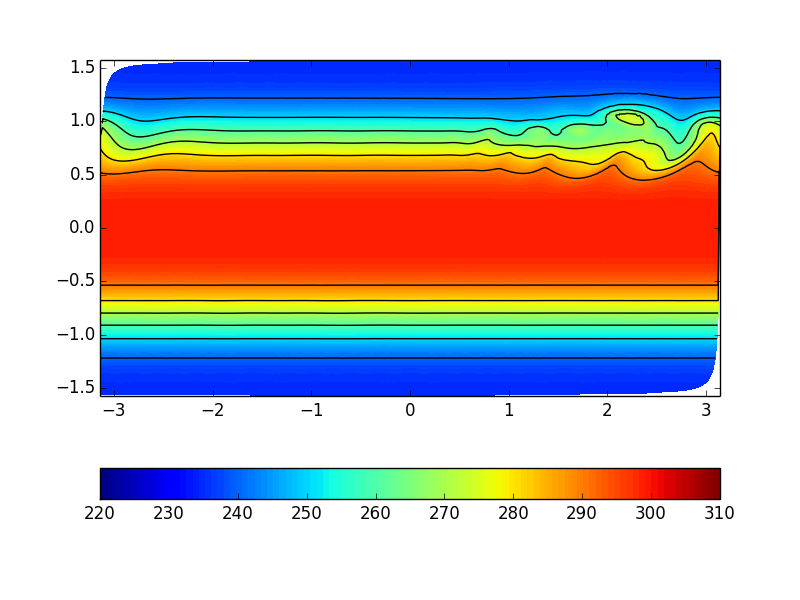}
\caption{Temperature, $T_h$ (in $^{\circ}\mathrm{K}$) at $z \approx 1.5\mathrm{km}$, day 8 (left) and 10 (right).}
\label{fig::temperature}
\end{center}
\end{figure}

\begin{figure}[!hbtp]
\begin{center}
\includegraphics[width=0.48\textwidth,height=0.36\textwidth]{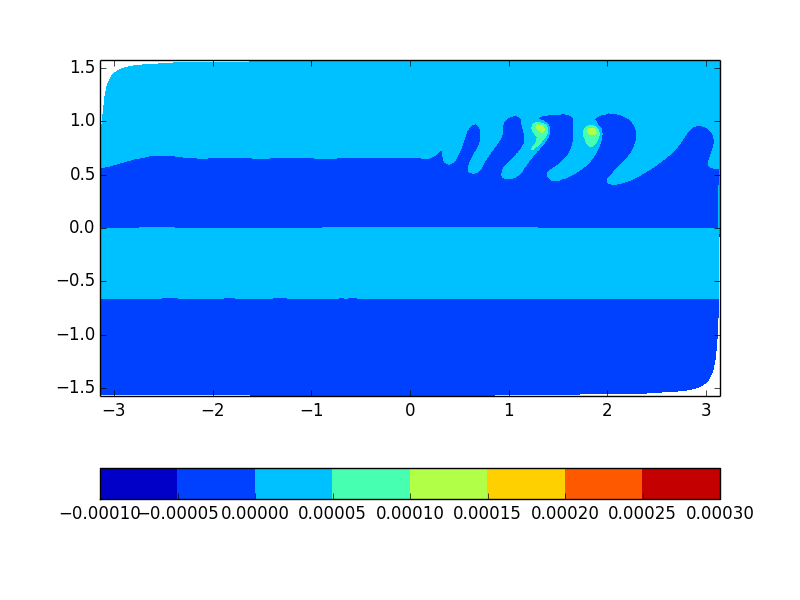}
\includegraphics[width=0.48\textwidth,height=0.36\textwidth]{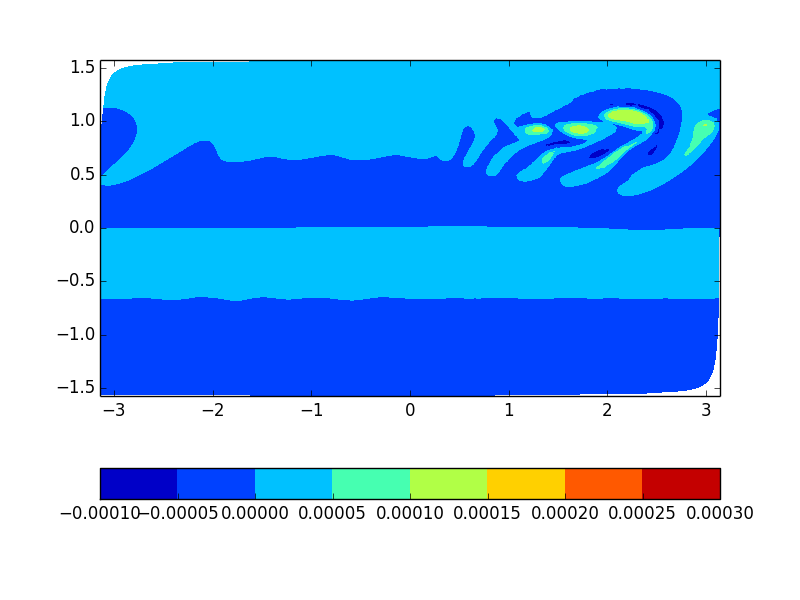}
\caption{Vertical component of the relative vorticity, $\omega_h$ (in $\mathrm{s}^{-1}$)
at $z \approx 1.5\mathrm{km}$, day 8 (left) and 10 (right).}
\label{fig::vorticity}
\end{center}
\end{figure}

\begin{figure}[!hbtp]
\begin{center}
\includegraphics[width=0.48\textwidth,height=0.36\textwidth]{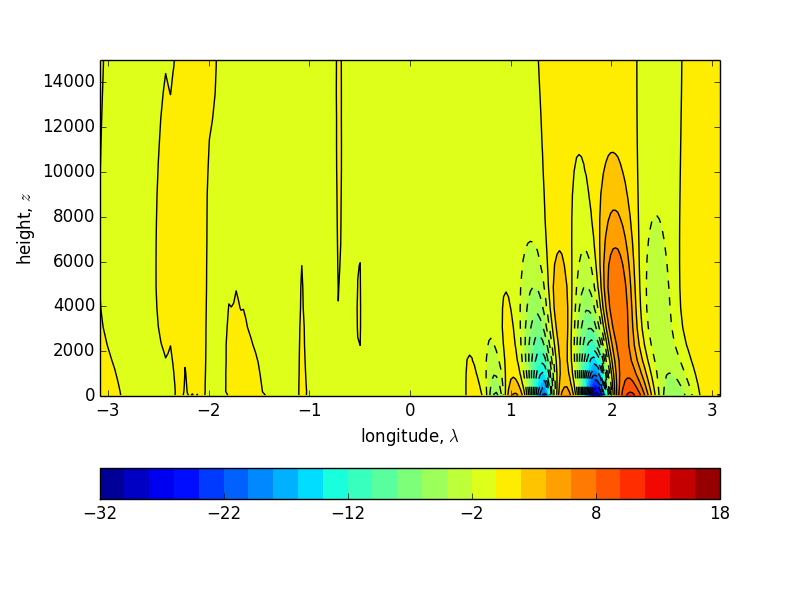}
\includegraphics[width=0.48\textwidth,height=0.36\textwidth]{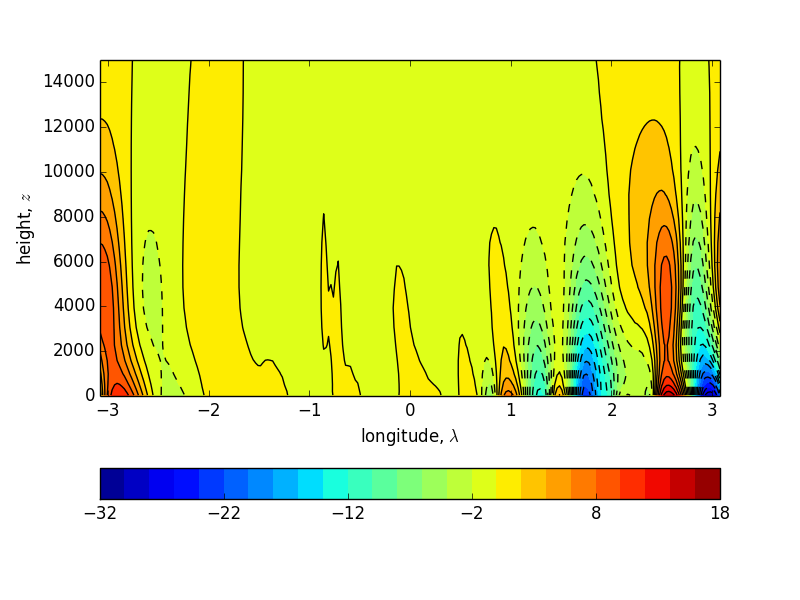}
\caption{Vertical cross section of the pressure perturbation, $p_h-\bar p_h$ (in $\mathrm{hPa}$)
at $50^{\circ}\mathrm{N}$, day 8 (left) and day 10 (right).}
\label{fig::pressure_cross_section}
\end{center}
\end{figure}

We present the bottom level pressure $p$ and temperature, $T$, and the vertical component of
the relative vorticity, $\omega$, at $z\approx 1.5\mathrm{km}$ at days $8$ and $10$ in Figs.
\ref{fig::pressure}, \ref{fig::temperature} and \ref{fig::vorticity}, as well as a meridional
cross section of the pressure perturbation at $50^{\circ}\mathrm{N}$ in Fig.
\ref{fig::pressure_cross_section}. In the cases of the pressure and temperature, these are
reconstructed from the model variables as $p=p_0(\Pi/c_p)^{c_p/R}$ and $T=\theta\Pi/c_p$.
The pressure perturbation in Fig. \ref{fig::pressure_cross_section} is derived by removing
the average pressure at the corresponding vertical level at $50^{\circ}\mathrm{S}$. These results
compare well with the previously published test case results \cite{Ullrich14} in both shape and
magnitude, and clearly show the signal of the baroclinic instability.

\begin{figure}[!hbtp]
\begin{center}
\includegraphics[width=0.48\textwidth,height=0.36\textwidth]{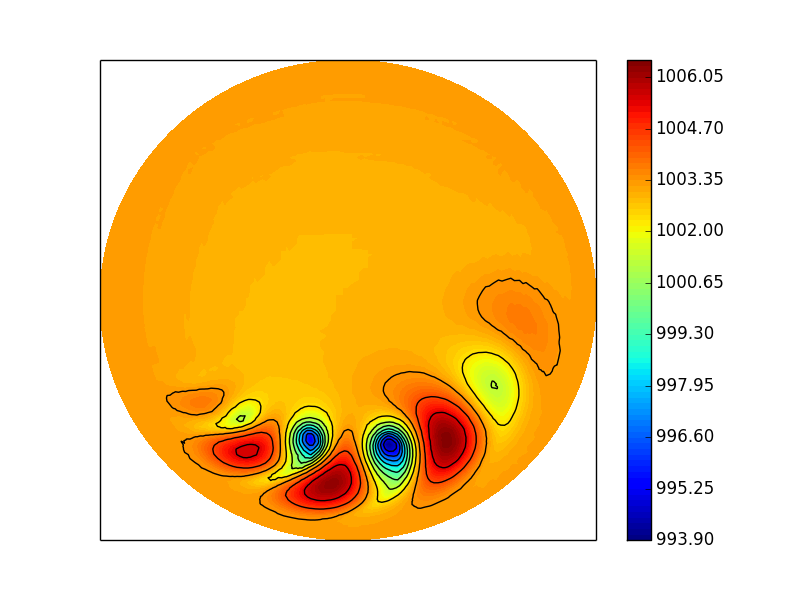}
\includegraphics[width=0.48\textwidth,height=0.36\textwidth]{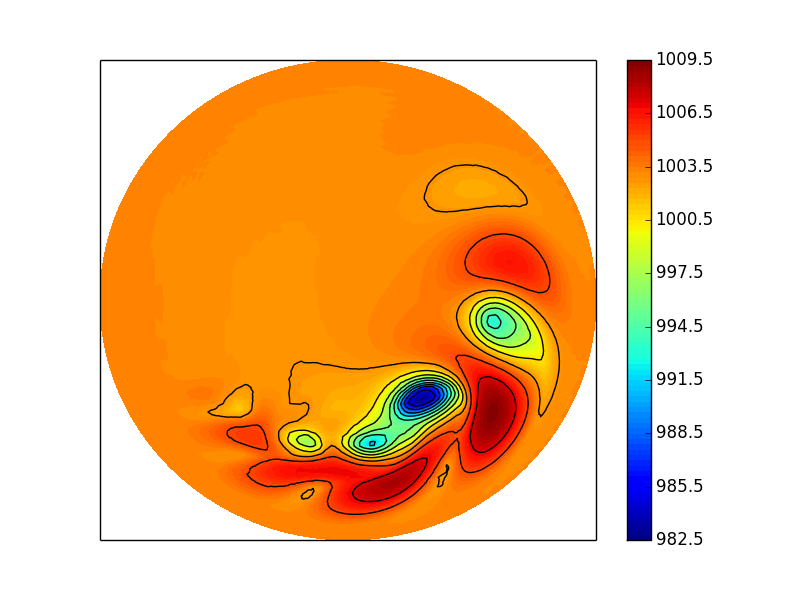}
\caption{Bottom level Exner pressure, $\Pi_h$ (in $\mathrm{m}^{2}\mathrm{s}^{-2}\mathrm{K}^{-1}$)
day 8 (left) and 10 (right).}
\label{fig::exner_pressure}
\end{center}
\end{figure}

\begin{figure}[!hbtp]
\begin{center}
\includegraphics[width=0.48\textwidth,height=0.36\textwidth]{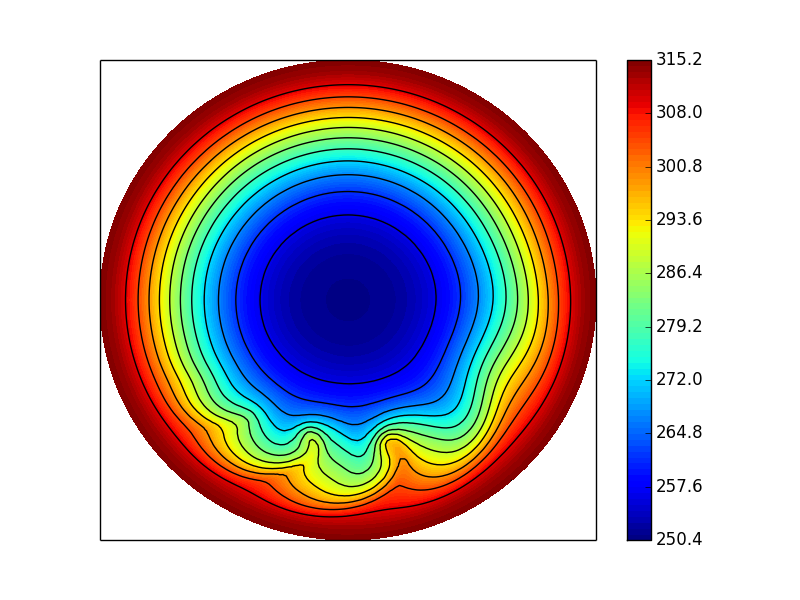}
\includegraphics[width=0.48\textwidth,height=0.36\textwidth]{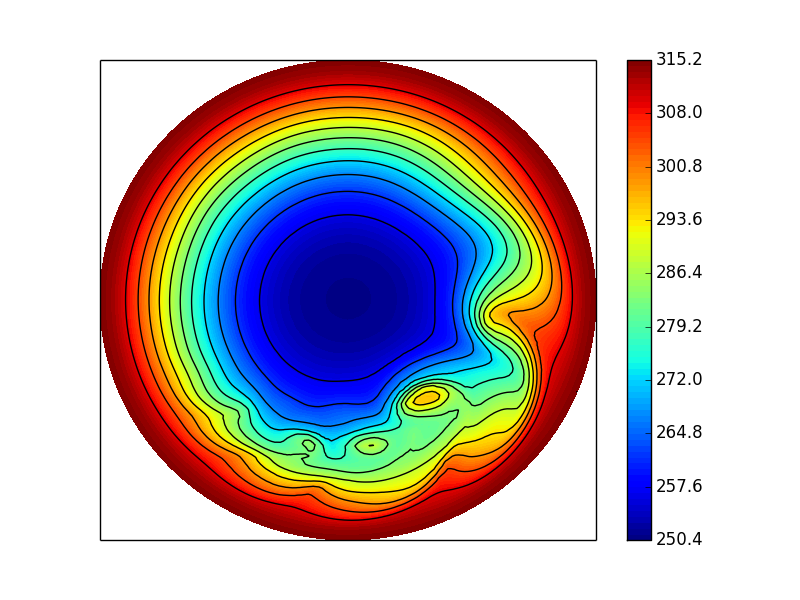}
\caption{Potential temperature, $\theta_h$ (in $^{\circ}\mathrm{K}$) at $z \approx  1.5\mathrm{km}$,
day 8 (left) and 10 (right).}
\label{fig::potential_temperature}
\end{center}
\end{figure}

\begin{figure}[!hbtp]
\begin{center}
\includegraphics[width=0.48\textwidth,height=0.36\textwidth]{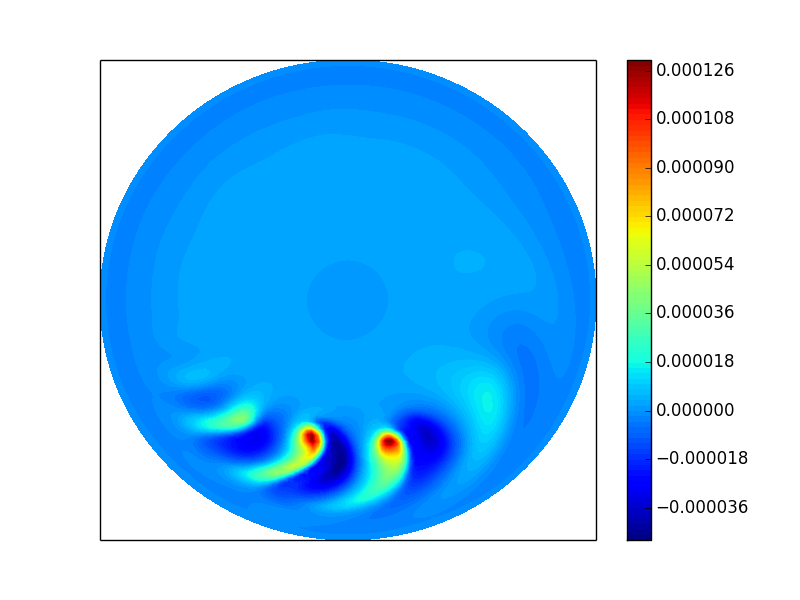}
\includegraphics[width=0.48\textwidth,height=0.36\textwidth]{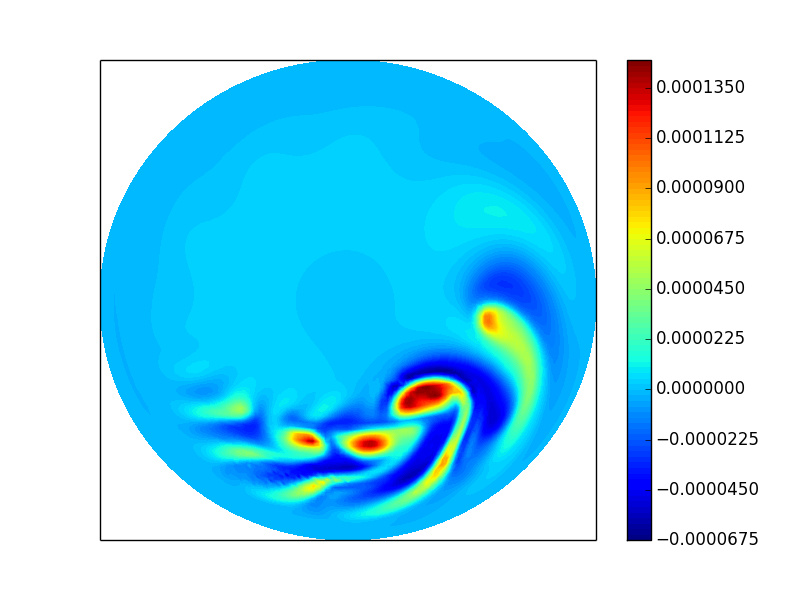}
\caption{Vertical component of the relative vorticity, $\omega_h$ (in $\mathrm{s}^{-1}$) at
$z \approx  1.5\mathrm{km}$, day 8 (left) and 10 (right). Results are the same as for Fig.
\ref{fig::vorticity}, but with color bars set to exact range of the data.}
\label{fig::vorticity_2}
\end{center}
\end{figure}

For completeness we also show the results at days 8 and 10 for the original model variables of
bottom level Exner pressure and potential temperature and relative vertical vorticity at $z\approx 1.5\mathrm{km}$
in Figs. \ref{fig::exner_pressure}, \ref{fig::potential_temperature} and \ref{fig::vorticity_2}
respectively.
These are presented for the northern hemisphere only, looking down from the north pole.
These results are perhaps slightly sharper than the reconstructed temperature and
pressure fields since they are interpolated directly from the degrees of freedom.

\subsection{Non-hydrostatic gravity wave}

The baroclinic instability test case detailed above provides an excellent means of validating
both the leading order balance relations and the horizontal and vertical coupling required to
correctly simulate a secondary nonlinear three-dimensional instability. However the limitation
of this test case is that it is configured for spatial resolutions at which the non-hydrostatic
dynamics are negligible. As such we also test the behaviour of the model for the propagation of
a non-hydrostatic gravity wave driven by a potential temperature perturbation on a planet with
a reduced radius $125$ times smaller than that of the earth. This test was originally proposed as part
of the 2012 DCMIP workshop, for which numerous hydrostatic and non-hydrostatic dynamical cores
presented results \cite{DCMIP31}. Specific details of the initial configuration can be found
within the DCMIP test case document on the web site.

As for the baroclinic instability test, we run the simulation with a resolution of
$24\times 24$ elements of degree $p=3$ in each cubed sphere panel. We use 16 evenly spaced
vertical levels over a total height of $10,000.0$m, and a time step of $\Delta t = 0.5$s for
a total simulation time of $3600$s. No Rayleigh damping or viscosity is applied in the vertical,
however we rescale the horizontal biharmonic viscosity by a factor of $2.0$ for both the momentum
and temperature equations for a value of $0.144\Delta x^{3.2}$. Since the horizontal
and vertical scales are of equal order in this configuration, we have also included the
additional vorticity term in \eqref{eq::mom_weak_discrete_algebraic_perp_split_incidence} and
the associated diagnostic term \eqref{eq:curl_u_omega_discrete_algebraic_par_perp_split_incidence},
both of which were omitted from the baroclinic test case, where the vertical scales were
$\mathcal{O}(10^{-3})$ with respect to the horizontal scales.

This test case is especially challenging for our higher order spectral element formulation, firstly
since we have not applied any sort of upwinding or monotonicity preservation method to either
the continuity or temperature equation, and secondly because unlike other models we formulate our
energy equation as a flux form equation for the density weighted potential temperature,
$\Theta=\rho\theta$, from which the potential temperature, $\theta$, is then diagnosed. This is
in contrast to a material form of the potential temperature advection, in which upwinding is
directly applied in the flux reconstruction \cite{Melvin19}. For this reason we configure our
model with slightly higher resolution than that specified \cite{DCMIP31}.

\begin{figure}[!hbtp]
\begin{center}
\includegraphics[width=0.48\textwidth,height=0.36\textwidth]{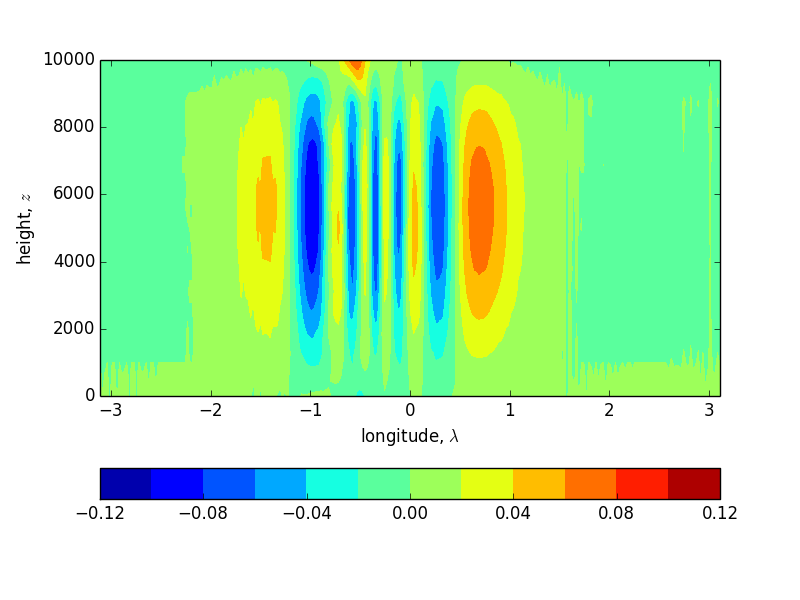}
\includegraphics[width=0.48\textwidth,height=0.36\textwidth]{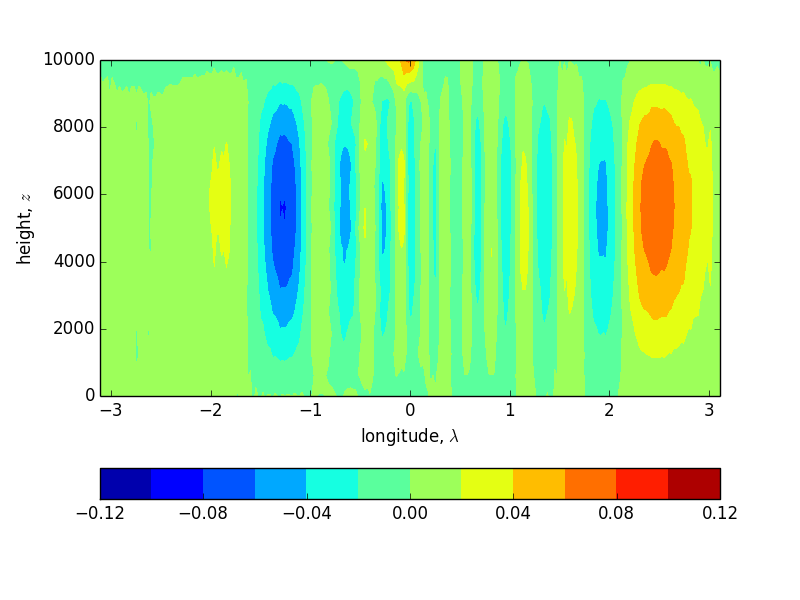}
\caption{Longitude-height equatorial ($\phi=0^{\circ}$) cross section of the potential temperature
perturbation, $\theta'(\lambda,0,z) = \theta(\lambda,0,z) - \bar{\theta}(z)$
at times $t=30$ minutes (left) and $t=60$ minutes (right).}
\label{fig::grav_wave_theta}
\end{center}
\end{figure}

Figure \ref{fig::grav_wave_theta} shows the longitude-height equatorial ($\phi=0^{\circ}$) cross section of the potential
temperature perturbation, $\theta'(\lambda,0,z) = \theta(\lambda,0,z) - \bar{\theta}(z)$, where $\bar{\theta}(z)$ is the
mean potential temperature at a given height, after 30 minutes and 1 hour. While the structure
and evolution of the perturbation qualitatively match the results presented for other non-hydrostatic
models, the absence of a monotone upwinded advection scheme means that our results are slightly
more oscillatory than those of other models. Moreover our results also show the ejection of a small
hot bubble which rises to the top of the domain, which is also most likely an artefact of our
non-monotone scheme.
In future work we will investigate energetically consistent methods for stabilising the advective
terms, as has previously been addressed for the shallow water equations \cite{Wimmer19}, and hence
recover less oscillatory solutions at non-hydrostatic scales. Fortunately the energetically
consistent flux form of the potential temperature equation presented here
(\ref{eq::temp_weak_discrete}, \ref{eq:theta_V_F_discrete}), matches with the formulation
necessary to ensure that tracer concentration fluxes are consistent with respect to the mass
flux in monotone schemes \cite{LLPR19}.

\begin{figure}[!hbtp]
\begin{center}
\includegraphics[width=0.48\textwidth,height=0.36\textwidth]{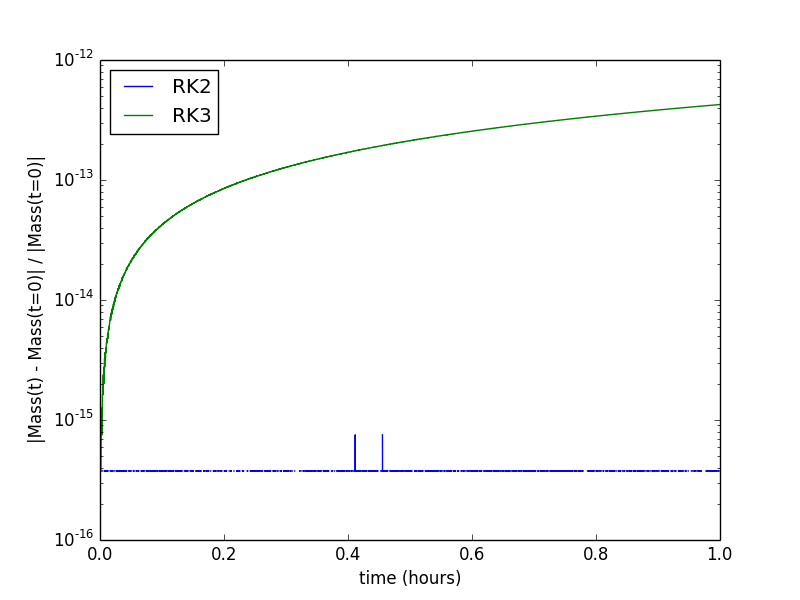}
\includegraphics[width=0.48\textwidth,height=0.36\textwidth]{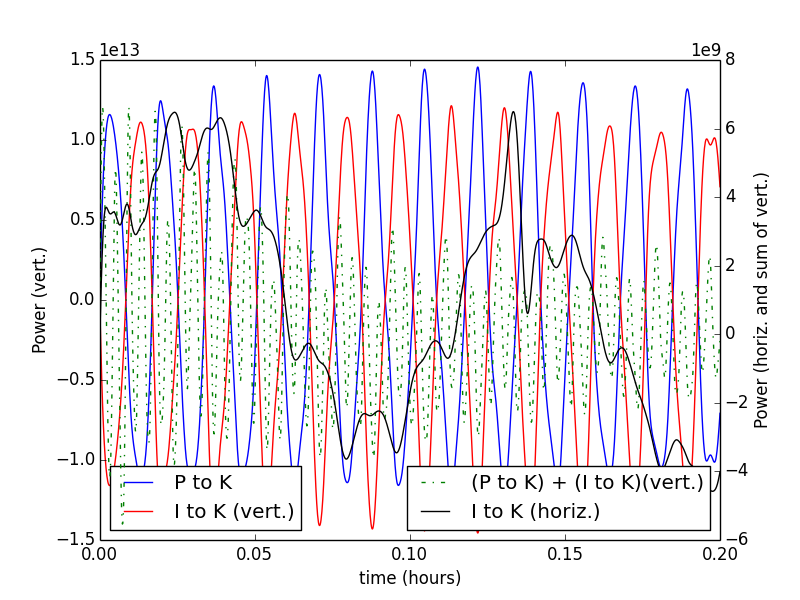}
\caption{Left: mass conservation error using second and third order Runge-Kutta integration in the horizontal.
Right: power associated with the vertical and horizontal energetic exchanges, (in $\mathrm{kg\cdot m^2s^{-3}}$).}
\label{fig::grav_wave_conservation}
\end{center}
\end{figure}

We also use this test to study the conservation properties of the model, as shown in Fig. \ref{fig::grav_wave_conservation}.
The third order Runge-Kutta scheme used for horizontal advection, as given in (\ref{eq::horiz_1}-\ref{eq::horiz_3})
involves a linear extrapolation from previous to current states. It would appear that this convex
operation leads to a small loss of mass that grows linearly with time. If we replace the third
order scheme with a second order scheme of the form $\boldsymbol{b}^{(1)} = \boldsymbol{b}' + \Delta tL^h(\boldsymbol{b}^n);
\boldsymbol{b}'' = \boldsymbol{b}' + \Delta t(L^h(\boldsymbol{b}^n) + L^h(\boldsymbol{b}^{(1)}))/2$,
which does not involve a linear extrapolation, then exact energy conservation is recovered.

Figure \ref{fig::grav_wave_conservation} also shows the power (in units of $\mathrm{kg\cdot m^2s^{-3}}$)
associated with the horizontal and vertical energetic exchanges, and the sum of vertical exchanges. Since
the vertical gravity wave oscillates at a high frequency compared to the horizontal dynamics, we only
show these results for the first 12 minutes of the simulation. The gravity wave is expressed through the
periodic exchange of potential to internal energy, via vertical motions. The power associated with
these exchanges is of $\mathcal{O}(10^{13}\mathrm{kg\cdot m^2s^{-3}})$, which is approximately
$\mathcal{O}(10^{4})$ times as large as both the sum of these exchanges, as well as the horizontal
kinetic to internal exchanges, which are observed to occur on a longer time scale than the gravity wave.

\section{Conclusions}

A model of the rotating compressible Euler equations on the cubed sphere using the mixed mimetic
spectral element method is presented. The model uses a Strang carryover directional splitting scheme
with an implicit Picard solver for the vertical dynamics in order to negate the CFL condition of the
vertical acoustic and gravity waves. The discontinuities in the discrete function spaces are exploited
so as to solve for each horizontal layer and each vertical element independently in order to avoid the
need to invert a global 3D mass matrix. The forcing and flux terms are constructed so as to take advantage
of the adjoint relations between the discrete gradient and divergence operators and balance the exchanges
of kinetic, potential and internal energy. The exception to this is the construction of the vertical
pressure gradient operator, which violates this principal in order to allow for a stable implicit solve
of the vertical dynamics, with the consequence that the vertical kinetic to internal energy
exchanges do not exactly balance.

In future work we will explore fully implicit formulations so as to avoid the horizontal vertical
splitting and restore the balance of vertical kinetic internal energy exchanges to the discrete system.
We will also explore energetically consistent flux reconstructions so as to recover less oscillatory
solutions at non-hydrostatic scales.

\blue{We also intend to optimise the parallel formulation of the code as part of additional future work. 
For the baroclinic instability simulation presented here at a resolution of $128\mathrm{km}$ on 30 vertical 
levels over 96 processors on a set of Dell PowerEdge M630 servers, the model only runs 
approximately 4.5 times faster than real time. By re-formulating the parallel decomposition of the different
function spaces on the cubed sphere we hope to better align these decompositions with the native PETSc
decompositions of vectors and matrices so as to reduce parallel communication.}

\section{Acknowledgments}

David Lee would like to thank Drs. Marcus Thatcher and John McGregor for their continued
support and access to computing resources, and Dr. Thomas Melvin for several helpful discussions.
We are also grateful to the three anonymous reviewers, each of whom contributed
thoughtful and insightful comments that helped improve the quality of this article.

	\appendix

	\section{Incident matrix} \label{appendix:incidence_matrix}

			\tikzset{
  on each segment/.style={
    decorate,
    decoration={
      show path construction,
      moveto code={},
      lineto code={
        \path [#1]
        (\tikzinputsegmentfirst) -- (\tikzinputsegmentlast);
      },
      curveto code={
        \path [#1] (\tikzinputsegmentfirst)
        .. controls
        (\tikzinputsegmentsupporta) and (\tikzinputsegmentsupportb)
        ..
        (\tikzinputsegmentlast);
      },
      closepath code={
        \path [#1]
        (\tikzinputsegmentfirst) -- (\tikzinputsegmentlast);
      },
    },
  },
  mid arrow/.style={postaction={decorate,decoration={
        markings,
        mark=at position .5 with {\arrow[#1]{stealth}}
      }}},
  end arrow/.style={postaction={decorate,decoration={
        markings,
        mark=at position 1.0 with {\arrow[#1]{stealth}}
      }}},
}
		\subsection{Discrete gradient (2D)}
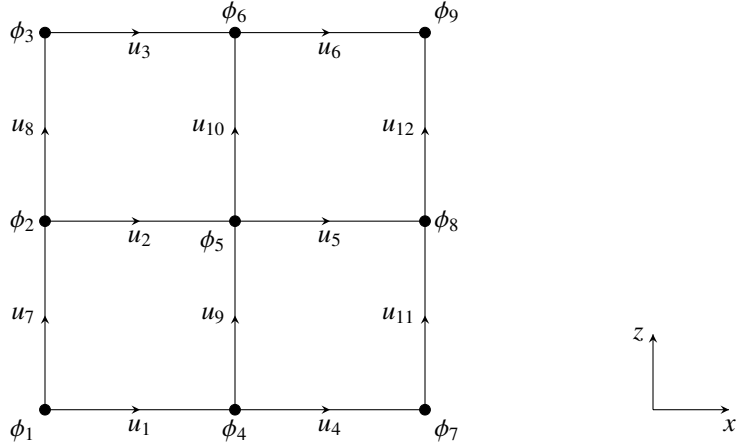
\begin{figure}[!h]
		\centering
	\begin{tikzpicture}
  		\path [draw=black,postaction={on each segment={mid arrow=black}}]
			(0.0, 0.0) -- (2.5, 0.0)
			(2.5, 0.0) -- (5.0, 0.0)
			(0.0, 2.5) -- (2.5, 2.5)
			(2.5, 2.5) -- (5.0, 2.5)
			(0.0, 5) -- (2.5, 5)
			(2.5, 5) -- (5.0, 5)
			(0.0, 0.0) -- (0.0, 2.5)
			(0.0, 2.5) -- (0.0, 5.0)
			(2.5, 0.0) -- (2.5, 2.5)
			(2.5, 2.5) -- (2.5, 5.0)
			(5, 0.0) -- (5, 2.5)
			(5, 2.5) -- (5, 5.0)
  		;
		\filldraw (0,0) circle (2pt) node[align=left,   below left] {$\phi_{1}$};
		\filldraw (0,2.5) circle (2pt) node[align=left,   left] {$\phi_{2}$};
		\filldraw (0,5.0) circle (2pt) node[align=left,   left] {$\phi_{3}$};
		\filldraw (2.5,0.0) circle (2pt) node[align=left,   below] {$\phi_{4}$};
		\filldraw (2.5,2.5) circle (2pt) node[align=left,   below left] {$\phi_{5}$};
		\filldraw (2.5,5) circle (2pt) node[align=left,   above] {$\phi_{6}$};
		\filldraw (5,0.0) circle (2pt) node[align=left,   below right] {$\phi_{7}$};
		\filldraw (5,2.5) circle (2pt) node[align=left,   right] {$\phi_{8}$};
		\filldraw (5,5) circle (2pt) node[align=left,   above right] {$\phi_{9}$};

		\draw (1.25, 0.0) node[align=left,   below] {$u_{1}$};
		\draw (1.25, 2.5) node[align=left,   below] {$u_{2}$};
		\draw (1.25, 5.0) node[align=left,   below] {$u_{3}$};

		\draw (3.75, 0.0) node[align=left,   below] {$u_{4}$};
		\draw (3.75, 2.5) node[align=left,   below] {$u_{5}$};
		\draw (3.75, 5.0) node[align=left,   below] {$u_{6}$};

		\draw (0, 1.25) node[align=left,   left] {$u_{7}$};
		\draw (0, 3.75) node[align=left,   left] {$u_{8}$};

		\draw (2.5, 1.25) node[align=left,   left] {$u_{9}$};
		\draw (2.5, 3.75) node[align=left,   left] {$u_{10}$};

		\draw (5, 1.25) node[align=left,   left] {$u_{11}$};
		\draw (5, 3.75) node[align=left,   left] {$u_{12}$};

		\path [draw=black,postaction={on each segment={end arrow=black}}]
			(8.0, 0.0) -- (9.0, 0.0)
			(8.0, 0.0) -- (8.0, 1.0)	
		;	
		\draw (9, 0.0) node[align=left,   below] {$x$};
		\draw (8, 1.0) node[align=left,   left] {$z$};
		
		\draw[white] (-4, 0.0) node[align=left,   below] {$x$};
	\end{tikzpicture}
	\caption{Mesh with orientations. All nodes oriented as positive and all edges positive in the direction of the arrows.}
	\label{fig:appendix_discrete_gradient}
	\end{figure}

		The continuous equation
		\begin{equation}
			\nabla\phi = \boldsymbol{u}\,, \label{eq:appendix_gradient}
		\end{equation}
		can be represented on the mesh depicted on \figref{fig:appendix_discrete_gradient} by the following algebraic system of equations

		\begin{equation}
			\boldsymbol{\mathsf{E}}^{1,0}\boldsymbol{\mathsf{\phi}} = \boldsymbol{\mathsf{u}}\,, \label{eq:incidence_0_1_appendix}
		\end{equation}
		where the incidence matrix $\boldsymbol{\mathsf{E}}^{1,0}$ is given by
		\begin{equation}
			\boldsymbol{\mathsf{E}}^{1,0} =
			\left[
				\begin{array}{ccccccccc}
					-1 & 0  & 0  & 1   & 0 & 0 & 0 & 0 & 0 \\
					0  & -1 & 0  & 0   & 1 & 0 & 0 & 0 & 0 \\ 
					0  & 0  & -1 & 0   & 0 & 1 & 0 & 0 & 0 \\ 
					0  & 0  &  0  & -1 & 0 & 0 & 1 & 0 & 0 \\
					0  & 0  & 0   & 0  & -1 & 0 & 0 & 1 & 0 \\ 
					0  & 0  & 0 & 0 & 0 & -1 & 0 & 0 & 1 \\ 
					-1 &  1 & 0  & 0  &  0 &  0 &  0 &  0 &  0 \\
					 0 &  -1 & 1  & 0  &  0 &  0 &  0 &  0 &  0 \\
					 0 &   0 & 0 & -1  &  1 &  0 &  0 &  0 &  0 \\
					  0 &   0 & 0  & 0  &  -1 &  1 &  0 &  0 &  0 \\
					  0 &   0 & 0  & 0  &  0 &  0 &  -1 &  1 &  0 \\
					  0 &   0 & 0  & 0  &  0 &  0 &  0 &  -1 &  1 \\
				\end{array}
			\right]\,.
		\end{equation}
		We note that we can easily split into a vertical component ($z$-direction) and a parallel component ($x$-direction). Therefore we can rewrite \eqref{eq:appendix_gradient} as
		\begin{equation}
			\nabla\phi = \gradpar\phi + \gradperp\phi = \boldsymbol{u}_{\parallel} + \boldsymbol{u}_{\perp}\,.
		\end{equation}
		In the same way, we can rewrite the discrete counterpart of this equation, \eqref{eq:incidence_0_1_appendix}, as
		\begin{equation}
			\left[
				\begin{array}{c}
					\incidencegradpar \\
					\\
					\incidencegradperp
				\end{array}
			\right]
			\boldsymbol{\mathsf{\phi}}
			=
			\left[
				\begin{array}{c}
					\boldsymbol{\mathsf{u}}_{\parallel} \\
					\\
					\boldsymbol{\mathsf{u}}_{\perp}
				\end{array}
			\right]\,,
		\end{equation}
		where
		\begin{equation}
			\incidencegradpar =
			\left[
				\begin{array}{ccccccccc}
					-1 & 0  & 0  & 1   & 0 & 0 & 0 & 0 & 0 \\
					0  & -1 & 0  & 0   & 1 & 0 & 0 & 0 & 0 \\ 
					0  & 0  & -1 & 0   & 0 & 1 & 0 & 0 & 0 \\ 
					0  & 0  &  0  & -1 & 0 & 0 & 1 & 0 & 0 \\
					0  & 0  & 0   & 0  & -1 & 0 & 0 & 1 & 0 \\ 
					0  & 0  & 0 & 0 & 0 & -1 & 0 & 0 & 1
				\end{array}
			\right]\,, \quad
			\incidencegradperp =
			\left[
				\begin{array}{ccccccccc}
					-1 &  1 & 0  & 0  &  0 &  0 &  0 &  0 &  0 \\
					 0 &  -1 & 1  & 0  &  0 &  0 &  0 &  0 &  0 \\
					 0 &   0 & 0 & -1  &  1 &  0 &  0 &  0 &  0 \\
					  0 &   0 & 0  & 0  &  -1 &  1 &  0 &  0 &  0 \\
					  0 &   0 & 0  & 0  &  0 &  0 &  -1 &  1 &  0 \\
					  0 &   0 & 0  & 0  &  0 &  0 &  0 &  -1 &  1
				\end{array}
			\right] \,,
		\end{equation}
		and
		\begin{equation}
			\boldsymbol{\mathsf{E}}^{1,0} =
			\left[
				\begin{array}{c}
					\incidencegradpar \\
					\\
					\incidencegradperp
				\end{array}
			\right] \,.
		\end{equation}

		\subsection{Discrete curl (2D)}
			The continuous equation
			\begin{equation}
				\nabla\times\boldsymbol{u} = \boldsymbol{\omega}\,, \label{eq:appendix_curl} 
			\end{equation}
			can be represented on the mesh depicted on \figref{fig:appendix_discrete_curl} by the following system of equations
			\begin{figure}[!h]
		\centering
	\begin{tikzpicture}
  		\path [draw=black,postaction={on each segment={mid arrow=black}}]
			(0.0, 0.0) -- (2.5, 0.0)
			(2.5, 0.0) -- (5.0, 0.0)
			(0.0, 2.5) -- (2.5, 2.5)
			(2.5, 2.5) -- (5.0, 2.5)
			(0.0, 5) -- (2.5, 5)
			(2.5, 5) -- (5.0, 5)
			(0.0, 0.0) -- (0.0, 2.5)
			(0.0, 2.5) -- (0.0, 5.0)
			(2.5, 0.0) -- (2.5, 2.5)
			(2.5, 2.5) -- (2.5, 5.0)
			(5, 0.0) -- (5, 2.5)
			(5, 2.5) -- (5, 5.0)
  		;

		\draw (1.25, 0.0) node[align=left,   below] {$u_{1}$};
		\draw (1.25, 2.5) node[align=left,   below] {$u_{2}$};
		\draw (1.25, 5.0) node[align=left,   below] {$u_{3}$};

		\draw (3.75, 0.0) node[align=left,   below] {$u_{4}$};
		\draw (3.75, 2.5) node[align=left,   below] {$u_{5}$};
		\draw (3.75, 5.0) node[align=left,   below] {$u_{6}$};

		\draw (0, 1.25) node[align=left,   left] {$u_{7}$};
		\draw (0, 3.75) node[align=left,   left] {$u_{8}$};

		\draw (2.5, 1.25) node[align=left,   left] {$u_{9}$};
		\draw (2.5, 3.75) node[align=left,   left] {$u_{10}$};

		\draw (5, 1.25) node[align=left,   left] {$u_{11}$};
		\draw (5, 3.75) node[align=left,   left] {$u_{12}$};

		\path [draw=black,postaction={on each segment={mid arrow=black}}]
        		(1.75,1.25) arc (0:360:0.5 and 0.5);
		\path [draw=black,postaction={on each segment={mid arrow=black}}]
        		(1.75,3.75) arc (0:360:0.5 and 0.5);
		\path [draw=black,postaction={on each segment={mid arrow=black}}]
        		(4.25,1.25) arc (0:360:0.5 and 0.5);
		\path [draw=black,postaction={on each segment={mid arrow=black}}]
        		(4.25,3.75) arc (0:360:0.5 and 0.5);
		
		\draw (1.25, 1.25) node[align=left] {$w_{1}$};
		\draw (1.25, 3.75) node[align=left] {$w_{2}$};
		\draw (3.75, 1.25) node[align=left] {$w_{3}$};
		\draw (3.75, 3.75) node[align=left] {$w_{4}$};
		
		\path [draw=black,postaction={on each segment={end arrow=black}}]
			(8.0, 0.0) -- (9.0, 0.0)
			(8.0, 0.0) -- (8.0, 1.0)	
		;	
		\draw (9, 0.0) node[align=left,   below] {$x$};
		\draw (8, 1.0) node[align=left,   left] {$z$};
		
		\draw[white] (-4, 0.0) node[align=left,   below] {$x$};
	\end{tikzpicture}
	\caption{Mesh with orientations. All edges positive in the direction of the arrows and surfaces positive in the direction of the rotations indicated.}
	\label{fig:appendix_discrete_curl}
	\end{figure}
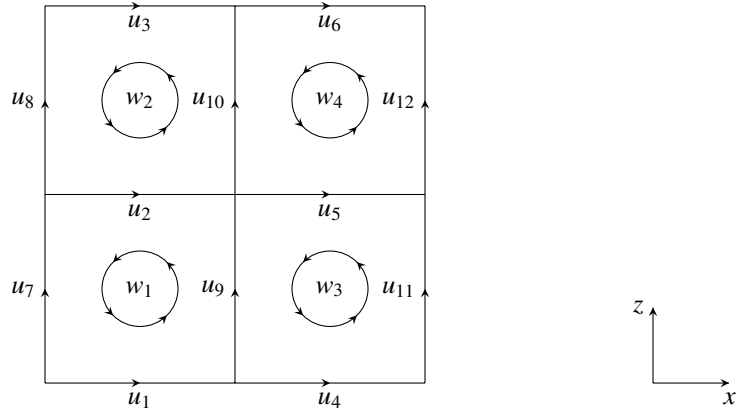
			\begin{equation}
				\incidencecurl \boldsymbol{\mathsf{u}} = \boldsymbol{\mathsf{w}}\,, \label{eq:incidence_1_2_appendix}
			\end{equation}
			where the incidence matrix $\incidencecurl$ is given by
			\begin{equation}
				\incidencecurl = 
				\left[
					\begin{array}{cccccccccccc}
						1   &   -1   &  0   &  0  &  0  &  0  &  -1  &  0   &  1  &  0  &  0  &  0  \\
						0   &    1   &  -1  &  0  &  0  &  0  &  0   &  -1  &  0  &  1  &  0  &  0  \\
						0   &    0   &  0  &  1  &  -1  &  0  &  0   &  0  &  -1  &  0  &  1  &  0  \\
						0   &    0   &  0  &  0  &  1   &  -1  &  0   &  0  &  0  &  -1  &  0  &  1
					\end{array}
				\right]\,.
			\end{equation}
			We note that, as done before, it is possible to split the curl operator into a vertical component ($z$-direction) and a parallel component ($x$-direction). Therefore we can rewrite \eqref{eq:appendix_curl}  as
			\begin{equation}
				\boldsymbol{\omega} = \underbrace{\overbrace{\curlpar\upar}^{\omegaparpar}+ \overbrace{\curlpar\uperp}^{\omegaparperp}}_{\omegapar} + \underbrace{\curlperp\upar}_{\omegaperp}\,. \label{eq:vorticity_splitting_appendix}
			\end{equation}
			Similarly, we can rewrite the discrete counterpart of this equation, \eqref{eq:incidence_1_2_appendix}, as
			\begin{equation}
			\left[
				\begin{array}{cc}
					\incidencecurlparpar & \incidencecurlparperp \\
					 & \\
					\incidencecurlperp & \boldsymbol{\mathsf{0}}
				\end{array}
			\right]
			\left[
				\begin{array}{c}
					\boldsymbol{\mathsf{u}}_{\parallel} \\
					\boldsymbol{\mathsf{u}}_{\perp}
				\end{array}
			\right]
			= 
			\left[
				\begin{array}{c}
					\boldsymbol{\mathsf{w}}_{\parallel} \\
					\boldsymbol{\mathsf{w}}_{\perp}
				\end{array}
			\right]
			\,.
			\end{equation}
			Note that since the mesh under consideration, \figref{fig:appendix_discrete_curl}, is a two dimensional mesh, $\boldsymbol{\omega}$ will only have a component in the $y$-direction. For this reason, $\incidencecurlperp = \boldsymbol{\mathsf{0}}$, and for compactness we suppress this term (for three dimensional meshes this term must be included)
			\begin{equation}
				\left[
				\begin{array}{cc}
					\incidencecurlparpar & \incidencecurlparperp
				\end{array}
			\right]
			\left[
				\begin{array}{c}
					\boldsymbol{\mathsf{u}}_{\parallel} \\
					\boldsymbol{\mathsf{u}}_{\perp}
				\end{array}
			\right]
			 = \boldsymbol{\mathsf{w}}_{\parallel}\,,
			\end{equation} 
			where
			\begin{equation}
				\incidencecurlparpar = 
				\left[
					\begin{array}{cccccc}
						1   &   -1   &  0   &  0  &  0  &  0   \\
						0   &    1   &  -1  &  0  &  0  &  0   \\
						0   &    0   &  0  &  1  &  -1  &  0    \\
						0   &    0   &  0  &  0  &  1   &  -1
					\end{array}
				\right]\,,
				\qquad
				\incidencecurlparperp =
				\left[
					\begin{array}{cccccc}
						  -1  &  0   &  1  &  0  &  0  &  0  \\
						   0   &  -1  &  0  &  1  &  0  &  0  \\
						  0   &  0  &  -1  &  0  &  1  &  0  \\
					      0   &  0  &  0  &  -1  &  0  &  1
					\end{array}
				\right]\,,
			\end{equation}
			and
			\begin{equation}
				\incidencecurl = 
				\left[
				\begin{array}{cc}
					\incidencecurlparpar & \incidencecurlparperp
				\end{array}
			\right]\,.
			\end{equation}
			
		\subsection{Discrete divergence (2D)}
			For the divergence operator, the continuous equation
			\begin{equation}
				\nabla\cdot\boldsymbol{u} = \sigma\,, \label{eq:divergence_appendix}
			\end{equation}
			can be expressed on the mesh depicted on \figref{fig:appendix_discrete_divergence} by the following system of equations
			
			\begin{figure}[!h]
		\centering
	\begin{tikzpicture}
  		\path [draw=black]
			(0.0, 0.0) -- (2.5, 0.0)
			(2.5, 0.0) -- (5.0, 0.0)
			(0.0, 2.5) -- (2.5, 2.5)
			(2.5, 2.5) -- (5.0, 2.5)
			(0.0, 5) -- (2.5, 5)
			(2.5, 5) -- (5.0, 5)
			(0.0, 0.0) -- (0.0, 2.5)
			(0.0, 2.5) -- (0.0, 5.0)
			(2.5, 0.0) -- (2.5, 2.5)
			(2.5, 2.5) -- (2.5, 5.0)
			(5, 0.0) -- (5, 2.5)
			(5, 2.5) -- (5, 5.0)
  		;

		\draw (1.25, 0.0) node[align=left,   below left] {$u_{7}$};
		\draw (1.25, 2.5) node[align=left,   below left] {$u_{8}$};
		\draw (1.25, 5.0) node[align=left,   below left] {$u_{9}$};

		\draw (3.75, 0.0) node[align=left,   below left] {$u_{10}$};
		\draw (3.75, 2.5) node[align=left,   below left] {$u_{11}$};
		\draw (3.75, 5.0) node[align=left,   below left] {$u_{12}$};

		\draw (0, 1.25) node[align=left,   below left] {$u_{1}$};
		\draw (0, 3.75) node[align=left,   below left] {$u_{2}$};

		\draw (2.5, 1.25) node[align=left,   below left] {$u_{3}$};
		\draw (2.5, 3.75) node[align=left,   below left] {$u_{4}$};

		\draw (5, 1.25) node[align=left,   below left] {$u_{5}$};
		\draw (5, 3.75) node[align=left,   below left] {$u_{6}$};
		
		\path [draw=black,postaction={on each segment={end arrow=black}}]
			(-0.25, 1.25) -- (0.25, 1.25)
			(-0.25, 3.75) -- (0.25, 3.75)
			(2.25, 1.25) -- (2.75, 1.25)
			(2.25, 3.75) -- (2.75, 3.75)
			(4.75, 1.25) -- (5.25, 1.25)
			(4.75, 3.75) -- (5.25, 3.75)
			(1.25, -0.25) -- (1.25, 0.25)
			(3.75, -0.25) -- (3.75, 0.25)
			(1.25, 2.25) -- (1.25, 2.75)
			(3.75, 2.25) -- (3.75, 2.75)
			(1.25, 4.75) -- (1.25, 5.25)
			(3.75, 4.75) -- (3.75, 5.25)
			;	
		
		\draw (1.25, 1.25) node[align=left] {$\sigma_{1}$};
		\draw (1.25, 3.75) node[align=left] {$\sigma_{2}$};
		\draw (3.75, 1.25) node[align=left] {$\sigma_{3}$};
		\draw (3.75, 3.75) node[align=left] {$\sigma_{4}$};
		
		\path [draw=black,postaction={on each segment={end arrow=black}}]
			(8.0, 0.0) -- (9.0, 0.0)
			(8.0, 0.0) -- (8.0, 1.0)	
		;	
		\draw (9, 0.0) node[align=left,   below] {$x$};
		\draw (8, 1.0) node[align=left,   left] {$z$};
		
		\draw[white] (-4, 0.0) node[align=left,   below] {$x$};
	\end{tikzpicture}
	\caption{Mesh with orientations. All fluxes positive in the direction of the arrows and volumes positive for sources (net outflux).}
	\label{fig:appendix_discrete_divergence}
	\end{figure}
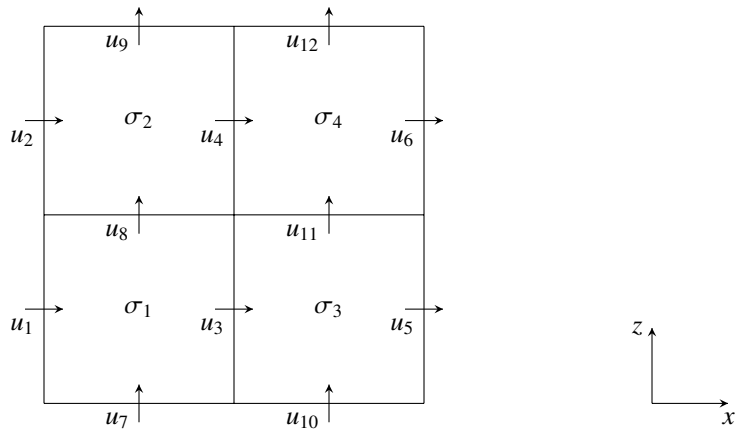
	\begin{equation}
		\incidencediv\boldsymbol{\mathsf{u}} = \boldsymbol{\mathsf{\sigma}}\,, \label{eq:divergence_discrete_appendix}
	\end{equation}
	where the incidence matrix is given by
	\begin{equation}
		\incidencediv = 
		\left[
			\begin{array}{cccccccccccc}
				-1   &   0   &  1   &  0  &  0  &  0  &  -1  &  1   &  0  &  0  &  0  &  0  \\
				0   &   -1   &  0   &  1  &  0  &  0  &  0  &  -1   &  1  &  0  &  0  &  0  \\
				0   &   0   &  -1   &  0  &  1  &  0  &  0  &  0   &  0  &  -1  &  1  &  0  \\
				0   &   0   &  0   &  -1  &  0  &  1  &  0  &  0   &  0  &  0  &  -1  &  1 
			\end{array}
		\right]\,.
	\end{equation}
	As discussed previously, it is possible to split the divergence operator into a vertical component ($z$-direction) and a parallel component ($x$-direction). Therefore we can rewrite \eqref{eq:divergence_appendix} as
	\begin{equation}
		\nabla\cdot\boldsymbol{u}_{\parallel} + \nabla\cdot\boldsymbol{u}_{\perp} = \sigma\,.
	\end{equation}
	In the same fashion, we can rewrite the discrete version of this equation, \eqref{eq:divergence_discrete_appendix}, as
	\begin{equation}
		\left[
			\begin{array}{cc}
				\incidencedivpar & \incidencedivperp
			\end{array}
		\right]
		\left[
				\begin{array}{c}
					\boldsymbol{\mathsf{u}}_{\parallel} \\
					\boldsymbol{\mathsf{u}}_{\perp}
				\end{array}
			\right]
			=
			\boldsymbol{\mathsf{\sigma}}\,,
	\end{equation}
	where
	\begin{equation}
		\incidencedivpar = \left[
			\begin{array}{cccccc}
				-1   &   0   &  1   &  0  &  0  &  0    \\
				0   &   -1   &  0   &  1  &  0  &  0   \\
				0   &   0   &  -1   &  0  &  1  &  0    \\
				0   &   0   &  0   &  -1  &  0  &  1   
			\end{array}
		\right]\,,
		\qquad
		\incidencedivperp = \left[
			\begin{array}{cccccc}
				-1  &  1   &  0  &  0  &  0  &  0  \\
				 0  &  -1   &  1  &  0  &  0  &  0  \\
				 0  &  0   &  0  &  -1  &  1  &  0  \\
				 0  &  0   &  0  &  0  &  -1  &  1 
			\end{array}
		\right]\,,
	\end{equation}
	and
	\begin{equation}
		\incidencediv = \left[
			\begin{array}{cc}
				\incidencedivpar & \incidencedivperp
			\end{array}
		\right]\,.
	\end{equation}

\vspace{1cm}

\end{document}